\documentclass{amsart}
\usepackage{mathptmx}      
\usepackage{latexsym}
\usepackage{amssymb}
\usepackage{bm}
\usepackage[usenames,dvipsnames]{color}
\usepackage{graphicx,pifont}
\usepackage{multicol,multirow}
\usepackage{ulem,cancel}
\usepackage{tikz}
\usetikzlibrary[patterns]
\usepackage{comment}
\usepackage{hyperref}
\usepackage{array}
\usepackage{caption}
\usepackage{subcaption}
\usepackage{txfonts}
\usepackage{enumerate}

\newtheorem{theorem}{Theorem}[section]

\theoremstyle{definition}
\newtheorem{definition}[theorem]{Definition}
\newtheorem{example}[theorem]{Example}
\theoremstyle{remark}
\newtheorem{remark}[theorem]{Remark}
\numberwithin{equation}{section}
\numberwithin{figure}{section}
\newtheorem{conjecture}[theorem]{Conjecture}

\newcommand{\bs}[1]{\boldsymbol{#1}}
\newcommand{\bsb}[1]{{\mbox{\sffamily #1}}}
\newcommand{\bb}[1]{{\mbox{\sffamily \bfseries #1}}}
\newcommand{\oname}[1]{\textrm{#1}}

\newcommand{\nrm}[1]{\left|\left|#1\right|\right|}
\newcommand{\eqdef}{\stackrel{\mathrm{def}}{=\joinrel=}}

\newcommand{\pp}[2]{\frac{\partial {#1}}{\partial {#2}}}

\newcommand{\termc}[1]{{\mbox{\sffamily {#1}}}}

\newcommand{\phf}[1]{{#1}+1/2}
\newcommand{\mhf}[1]{{#1}-1/2}
\newcommand{\powth}{{\textrm{\scriptsize th}}}
\newcommand{\powrd}{{\textrm{\scriptsize rd}}}

\newcommand{\pownd}{{\textrm{\scriptsize nd}}}

\begin{document}

\title[FD-FV for first-order hyperbolic conservation laws]{A hybrid FD-FV method for first-order hyperbolic conservation laws on Cartesian grids: the smooth problem case}
\author[X.~Zeng]{Xianyi Zeng}
\address{Department of Mathematical Sciences, University of Texas at El Paso, El Paso TX 79902, United States.\\
         Tel.: +1-915-747-6759}
\email[X.~Zeng]{xzeng@utep.edu}
\date{\today}



\begin{abstract}
We present a class of hybrid FD-FV (finite difference and finite volume) methods for solving general hyperbolic conservation laws written in first-order form.
The presentation focuses on one- and two-dimensional Cartesian grids; however, the generalization to higher dimensions is straightforward.
These methods use both cell-averaged values and nodal values as dependent variables to discretize the governing partial differential equation (PDE) in space, and they are combined with method of lines for integration in time.
This framework is absent of any Riemann solvers while it achieves numerical conservation naturally.
This paper focuses on the accuracy and linear stability of the proposed FD-FV methods, thus we suppose in addition that the solutions are sufficiently smooth.
In particular, we prove that the spatial-order of the FD-FV method is typically one-order higher than that of the discrete differential operator, which is involved in the construction of the method.
In addition, the methods are linearly stable subjected to a Courant-Friedrich-Lewy condition when appropriate time-integrators are used.
The numerical performance of the methods is assessed by a number of benchmark problems in one and two dimensions.
These examples include the linear advection equation, nonlinear Euler equations, the solid dynamics problem for linear elastic orthotropic materials, and the Buckley-Leverett equation.
\end{abstract}

\maketitle

\section{Introduction}
\label{sec:intro}
Hyperbolic conservation laws~\cite{CMDafermos:2010a} are used to describe many problems in the areas of fluid dynamics such as compressible flows, aeroacoustics, magnetohydrodynamics.
More recently, they are also used to formulate the problems of solid dynamics in both infinitesimal strain regime~\cite{GILemoine:2013a} and finite strain regime~\cite{MAguirre:2014a}.

Conventional numerical methods for conservation laws are typically constructed in a way such that the same discretization concept is used to construct the discrete formula of each dependent variable.
In the case of many conventional finite difference (FD) methods, such as~\cite{AJameson:1995a}, this discrete formula represents a discretization of governing equation in its strong form; 
whereas in the case of classical finite volume (FV) methods, such as those in~\cite{BvanLeer:1979a} and~\cite{RLeVeque:2002a}, the discrete formula represents a numerical realization of the integral form of the governing equation to achieve natural conservation.

One disadvantage of these methods, however, is that they typically require large (and overlapping) stencils to achieve higher-order accuracy.
For example, to achieve second-order accuracy in space, an explicit one-dimensional finite volume method requires information of five cells to update one particular cell average (i.e., two additional cells to each side); 
and a similar situation holds for the finite difference methods.
A more compact formulation can be achieved by the discontinuous Galerkin (DG) formulation (see, for example~\cite{BCockburn:1998a}), however, at the cost of having duplicate degrees of freedom on the interface between elements.

In recent years, there is continuous interest in developing numerical methods that are more compact than the conventional FD or FV methods, but still require fewer unknowns than the DG-type methods, for example, by using different discretization formula for different unknowns.
To the best knowledge of the author, numerical methods of this type date back to the 1970s: van Leer~\cite{BvanLeer:1977b} constructed a third-order accurate method (scheme V of the reference) that used both cell averages and nodal values at cell faces for 1D conservation laws.
This work, however, is not further explored.
A non-exhausted list of similar works that employ variables with different meanings includes the multi-moment methods (\cite{FXiao:2005a} and \cite{YImai:2008a}), the space-time control element and solution element methods (\cite{SCChang:1995a}, \cite{SChang:1999a}, and \cite{ZCZhang:2002a}), the staggered mesh method~\cite{RSanders:1992a} and dual-mesh method~\cite{HTHuynh:1995a}, the $\textrm{P}_N\textrm{P}_M$ methods~\cite{MDumbser:2008a}, and more recently the active flux scheme~\cite{TAEymann:2011a}.

The present work also falls into this category in that we choose the discrete variables to approximate both cell averages and nodal values that collocate at cell interfaces.
In particular, we focus on discretization in space and use the method of lines for integration in time. 
In our computational framework, on the one hand, the semi-discretization formula for the cell averages is constructed using the integral form of the governing equation, just like the FV methods, but with the numerical fluxes replaced by physical fluxes evaluated at nodal values.
On the other hand, the semi-discretization formula for the nodal values is constructed by local linearization and applying discrete differential operators that involve both nodal values and cell averages.
In the one-dimensional case, the latter construction has direct connection to the Hermite interpolation polynomials.
Due to these reasons, we call the present method the hybrid FD-FV methods (or simply FD-FV methods).
The FD-FV methods are more compact than the conventional finite difference and finite volume methods.
For example in one space dimension, we have in average two variables for each mesh cell, which enable us to construct higher-order discrete differential operators than the finite difference methods given the same stencil.
Furthermore, we will show that using a discrete differential operator that is constructed to $p^\powth\textrm{\sc-}$order operator in space, the resulting FD-FV method is formally $(p+1)^\powth\textrm{\sc-}$order accurate.
Thus we are able to construct a fourth-order accurate FD-FV method using a stencil of three consecutive cells (see the method given by $[\mathcal{D}_x^{3,up-biased}]$ in Section~\ref{sec:appr_stab}), whereas in the case of conventional FD or FV methods the spatial-order is only first-order.

In this paper, we only consider explicit time-integrators for simplicity.
In addition, because we focus on the analysis of the spatial order of accuracy of the proposed FD-FV methods and their linear stability properties, we also assume that the solutions are sufficiently smooth.
Enhancement of the nonlinear stability of the methods and their performance to solve discontinuous solutions will be addressed in future work.
To this end, the remainder of the paper is organized as follows.
We describe the general formulation of the FD-FV methods for one-dimensional first-order hyperbolic system in Section~\ref{sec:appr}.
The proofs of the existence of the coefficients and the theorem concerning the spatial order of accuracy of the FD-FV methods are provided in Appendix~\ref{app:hermite} and Appendix~\ref{app:order}, respectively.
Section~\ref{sec:ext} extends the FD-FV methodology to two-dimensional Cartesian grids, with evidence of the conjectured order of spatial accuracy provided in Appendix~\ref{app:2nd}.
Numerical examples are offered in Section~\ref{sec:num}, and finally Section~\ref{sec:concl} concludes this paper.

\section{The hybrid FD-FV approach for one-dimensional problems}
\label{sec:appr}
We describe in this section the fundamental framework of the hybrid FD-FV method for the one-dimensional hyperbolic conservation law in first-order form:
\begin{equation}\label{eq:appr_1dhcl}
\bs{w}_t + \bs{f}(\bs{w})_x = 0,
\end{equation}
which is defined on the rectangular space-time domain $(x,t)\in\Omega\times[0,T]$, and supplemented by appropriate initial and boundary conditions.
We assume that $\bs{w}(x,t)$ takes values in $\mathbb{R}^d$, and $\bs{f}:\mathbb{R}^d\to\mathbb{R}^d$ is a continuously differentiable function of $\bs{w}$.
Here and throughout the remainder of the paper, we use normal font variables such as $t$ and $x$ to denote scalars, whereas we use bold font to designate vectors or matrices.

\subsection{General FD-FV formulation for scalar problems}
\label{sec:appr_gen}
We first divide the computational domain $\Omega$ into cells with non-overlapping interiors: $\Omega = \bigcup_j\mathcal{C}_j$, where $\mathcal{C}_j=[x_{\mhf{j}},x_{\phf{j}}]$.
We denote by $h_j=x_{\phf{j}}-x_{\mhf{j}}$ and $x_j = (x_{\mhf{j}}+x_{\phf{j}})/2$ the length and the center of the cell $\Omega_j$, respectively.
In this paper, we assume uniform grids ($h_j\equiv h$), but the methodology applies equally to nonuniform 1D grids.

For simplicity, we first consider scalar problems, in which case $\bs{w}$ and $\bs{f}$ of (\ref{eq:appr_1dhcl}) are replaced by scalar functions $w$ and $f$:
\begin{equation}\label{eq:appr_gen_1dhcl}
w_t + f(w)_x = 0.
\end{equation}
The target of the hybrid FD-FV methods is to compute the numerical approximations to the nodal values at cell faces:
\begin{equation}\label{eq:appr_gen_n}
w^n_{\phf{j}} \approx w(x_{\phf{j}},t^n),
\end{equation}
and those to cell-averaged values for each cell:
\begin{equation}\label{eq:appr_gen_c}
\overline{w}^n_j \approx \frac{1}{h}\int_{\mathcal{C}_j}w(x,t^n)dx.
\end{equation}
Because we use the method of lines to decouple the discretization schemes in space and in time, we only consider the semi-discretized variables in this section, namely:
\begin{displaymath}
w_{j+1/2}(t) \approx w(x_{j+1/2},t),\quad\overline{w}_{j}(t) \approx \frac{1}{h}\int_{\mathcal{C}_j}w(x,t)dx.
\end{displaymath}
Here and throughout the remainder of the paper, we use $\overline{\dagger}_j^n,\ \dagger_{j+1/2}^n$ and $\overline{\dagger}_j,\ \dagger_{j+1/2}$ to designate the approximate cell averages and nodal values to the fully discretized system and semi-discretized system, respectively.
In the latter case, we often omit expressing the dependence in $t$ for convenience.
\begin{remark}\label{rk:appr_init_c}
For the purpose of analysis, we suppose that the cell averages are initialized exactly by the initial condition.
But this is not necessary in practical computations.
In the numerical examples provided in Section~\ref{sec:num}, we use appropriate quadrature rules to initialize the cell averages $\overline{w}^0_j = \overline{w}_j(0)$.
\end{remark}

On the one hand, in order to derive the semi-discretized formula for cell averages, we integrate equation~(\ref{eq:appr_gen_1dhcl}) over each cell $\mathcal{C}_j$ to obtain:
\begin{displaymath}
\frac{d}{dt}\left(\frac{1}{h}\int_{\mathcal{C}_j}w(x,t)dx\right)+\frac{1}{h}\left[f\left(w(x_{j+1/2},t)\right)-f\left(w(x_{j-1/2},t)\right)\right] = 0,
\end{displaymath}
which then leads to the following discretization scheme in space for $\overline{w}_j$:
\begin{equation}\label{eq:appr_gen_c_eq}
\frac{d\overline{w}_j}{dt} + \frac{1}{h}\left[f\left(w_{j+1/2}\right)-f\left(w_{j-1/2}\right)\right] = 0,
\end{equation}
where $w_{j\pm1/2}$ are the nodal variables.

On the other hand, we linearize the governing equation at a node $x_{\phf{j}}$ to derive the semi-discretized formula for nodal variables.
Supposing that $J(w)\eqdef f'(w)$ is the Jacobian of the flux function, we seek a numerical scheme for constructing the ODE for $w_{\phf{j}}$ in the form:
\begin{equation}\label{eq:appr_gen_n_eq}
\frac{dw_{\phf{j}}}{dt} + J(w_{\phf{j}})[\mathcal{D}_xw]_{\phf{j}} = 0.
\end{equation}
Here we use $[\mathcal{D}_x]$ to denote a discrete differential operator that uses both adjacent cell-averaged values and nodal values, to distinguish it from the conventional finite difference operators $\mathcal{D}_x$, which only involves the nodal variables.
A general (linear) form of this operator is:
\begin{equation}\label{eq:appr_gen_ddo}
\left[\mathcal{D}_xw\right]_{j+1/2} = \frac{1}{h}\sum_{l=-q+1}^q\alpha_l\overline{w}_{j+l} + \frac{1}{h}\sum_{l=-q}^q\beta_lw_{j+1/2+l},
\end{equation}
where $q$ is a positive integer that indicates the stencil of the method, and $\alpha_l$ and $\beta_l$ are constants to be determined later.
\begin{remark}\label{rk:appr_fdm}
One may recover a conventional finite difference method, by setting all constants $\alpha_l$ to zero.
In this case, the equations for the nodal variables decouple from the cell averages.
Thus in theory, the conventional finite difference methods are a subset of the present hybrid FD-FV framework.
However, this is not encouraged, for the reason discussed in the accuracy analysis below.
\end{remark}

\subsection{Spatial accuracy of FD-FV operators}
\label{sec:appr_acc}
There are mainly two considerations in design the constants $\alpha_l$ and $\beta_l$: the first one is that we want to achieve a certain order of accuracy, and the second one is that we need at least linear stability when these operators are applied to solve an simple advection problem.
To this end, we first focus on the accuracy analysis of the resulting FD-FV method by using (\ref{eq:appr_gen_c_eq}) and (\ref{eq:appr_gen_n_eq}) with a set of constants in this section, and move on to the linear stability analysis in Section~\ref{sec:appr_stab}.

We apply the Taylor series expansion to the nodal variables and cell-averaged variables to study the accuracy of the FD-FV operators $[\mathcal{D}_x]$ and the resulting FD-FV method.
Considering performing the expansion about a node $x_{\phf{j}}$, we have for a nearby nodal variable $w(x_{\phf{j+l}})$ (here we omit the dependence on $t$ for simplicity):
\begin{equation}\label{eq:appr_acc_taylor_n}
w(x_{\phf{j+l}}) = w(x_{\phf{j}}) + \sum_{m=1}^\infty\frac{(lh)^m}{m!}\partial_x^mw(x_{\phf{j}}),
\end{equation}
and for the average of the cell $\mathcal{C}_{j+l}$:
\begin{equation}
\label{eq:appr_acc_taylor_c}
\begin{array}{>{\displaystyle}r>{\displaystyle}c>{\displaystyle}l}
\frac{1}{h}\int_{x_{\mhf{j+l}}}^{x_{\phf{j+l}}}w(x)dx & = & \frac{1}{h}\int_{x_{\mhf{j+l}}}^{x_{\phf{j+l}}}\left(w(x_{j+1/2})+\sum_{m=1}^\infty\frac{(x-x_{j+1/2})^m}{m!}\partial_x^mw(x_{j+1/2})\right)dx \\
& = & w(x_{j+1/2}) + \sum_{m=1}^\infty\frac{[l^{m+1}-(l-1)^{m+1}]h^{m}}{(m+1)!}\partial_x^mw(x_{j+1/2}),
\end{array}
\end{equation}
where we use $\partial_x^m$ to denote the $m^\powth$ derivative in $x$.
To this end, we may define the order of the operator $[\mathcal{D}_x]$ by Definition~\ref{def:appr_acc_ddo}, and relate it to the coefficients by Theorem~\ref{thm:appr_acc_ddo}.
\begin{definition}\label{def:appr_acc_ddo}
The FD-FV operator $[\mathcal{D}_x]$ defined by~(\ref{eq:appr_gen_ddo}) is $p^\powth\textrm{\sc-}$order accurate for some integer $p\ge1$ if:
\begin{equation}\label{eq:appr_acc_ddo}
[\mathcal{D}_xw]_{\phf{j}} = \partial_xw_{\phf{j}}+c_p\partial_x^{p+1}w_{\phf{j}}h^p + O(h^{p+1})
\end{equation}
for any sufficiently smooth function $w(x)$, where each subscript $_{\phf{j}}$ denotes the exact data at $x_{\phf{j}}$.
In this formula, $c_p\ne0$ is a constant that is independent of the function $w$ and the cell size $h$.
\end{definition}
\begin{example}\label{ex:appr_acc}
The operator $[\mathcal{D}_xw]_{\phf{j}}=(w_{\phf{j}}-\overline{w}_j)/h$ is a first-order operator with $c_1=-1/3$;
and the operator $[\mathcal{D}_xw]_{\phf{j}}=2(-w_{j+3/2}+3\overline{w}_{j+1}-2w_{\phf{j}})/h$ is a second-order operator with $c_2=-1/12$.
\end{example}
\begin{theorem}\label{thm:appr_acc_ddo}
The FD-FV operator~(\ref{eq:appr_gen_ddo}) is $p^\powth\textrm{\sc-}$order accurate for some integer $p\ge1$ if and if the constants $\alpha_l, -q+1\le l\le q$ and $\beta_l, -q\le l\le q$ satisfy:
\begin{equation}\label{eq:appr_acc_ddo_iff}
\begin{array}{>{\displaystyle}r>{\displaystyle}c>{\displaystyle}l>{\displaystyle}l}
\sum_{l=-q+1}^q\frac{l^{m+1}-(l-1)^{m+1}}{(m+1)!}\alpha_l + 
\sum_{l=-q}^q\frac{l^m}{m!}\beta_l &=& \delta_{m1}, &\quad m=0,1,\cdots,p, \\
\sum_{l=-q+1}^q\frac{l^{p+2}-(l-1)^{p+2}}{(p+2)!}\alpha_l+\sum_{l=-q}^q\frac{l^{p+1}}{(p+1)!}\beta_l &=& c_p.
\end{array}
\end{equation}
Here $\delta_{m1}$ is the Kronecker delta that takes the value $1$ when $m=1$ and $0$ otherwise.
\end{theorem}
The proof of Theorem~\ref{thm:appr_acc_ddo} is straightforward: one simply needs to plug (\ref{eq:appr_acc_taylor_n}--\ref{eq:appr_acc_taylor_c}) into (\ref{eq:appr_gen_ddo}) and use Definition~\ref{def:appr_acc_ddo}.
What is more interesting is the existence of the coefficients that satisfy~(\ref{eq:appr_acc_ddo_iff}).
It turns out that at least in the 1D case, these coefficients have close relation to the construction of Hermite interpolation polynomials, see Appendix~\ref{app:hermite}.
Thus in most situations, we can construct the coefficients $\alpha_l$ and $\beta_l$ that leads to an arbitrary order operator $[\mathcal{D}_x]$.
Note that the preceding results extend naturally to nonuniform grids, seeing their close relation to the Hermite interpolation theory.

Our next target is to establish the connection of the order of $[\mathcal{D}_x]$ and that of the resulting semi-discretized FD-FV scheme, between which the latter is to be defined.
Considering applying the FD-FV method~(\ref{eq:appr_gen_c_eq}--\ref{eq:appr_gen_n_eq}) to solve the Cauchy problem for the scalar advection equation:
\begin{equation}\label{eq:appr_acc_adv}
w_t + cw_x = 0,\quad (x,t)\in\mathbb{R}\times[0,T],
\end{equation}
for some constant $c\ne0$ and the simple wave initial condition $w(x,0) = e^{ikx}$, where $k$ is the wavenumber.
The exact solutions for the cell averages and nodal variables at $t=T$ are:
\begin{equation}\label{eq:appr_acc_adv_sol}
\overline{w}^\ast_j = \frac{1}{h}\int_{(j-1)h}^{jh}e^{ik(x-cT)}dx = \frac{1}{i\theta}e^{-ickT}e^{ij\theta}\left(1-e^{-i\theta}\right),\quad
w^\ast_{j+1/2} = e^{-ickT}e^{ij\theta}.
\end{equation}
Here $i$ denotes the imaginary unit, $\theta=kh$ is the numerical wavenumber, and we suppose $x_{j-1/2}=(j-1)h$ as well as $x_j=(j-1/2)h$ for the ease of expression;
we also use the superscript $^\ast$ to denote the exact solutions.

The semi-discretized solutions to (\ref{eq:appr_acc_adv}) by using (\ref{eq:appr_gen_c_eq}--\ref{eq:appr_gen_n_eq}) leads to a coupled ODE system for $\overline{w}_j(t)$ and $w_{j+1/2}(t)$:
\begin{equation}\label{eq:appr_acc_ode}
\begin{array}{>{\displaystyle}r>{\displaystyle}c>{\displaystyle}l}
\frac{d\overline{w}_j}{dt} + \frac{c}{h}(w_{\phf{j}}-w_{\mhf{j}}) &=& 0, \\
\frac{dw_{\phf{j}}}{dt} + c[\mathcal{A}w]_{\phf{j}} + c[\mathcal{N}w]_{\phf{j}} &=& 0,
\end{array}
\end{equation}
where we split the operator $[\mathcal{D}_x]$ into two components by $[\mathcal{D}_x] = [\mathcal{A}] + [\mathcal{N}]$:
\begin{equation}\label{eq:appr_acc_ddo_spl}
[\mathcal{A}w]_{\phf{j}} \eqdef \frac{1}{h}\sum_{l=-q+1}^q\alpha_l\overline{w}_{j+l},\quad
[\mathcal{N}w]_{\phf{j}} \eqdef \frac{1}{h}\sum_{l=-q}^q\beta_lw_{\phf{j+l}}.
\end{equation}
It is not difficult to see that if the initial data is exact, the solutions to (\ref{eq:appr_acc_ode}) are given by the form:
\begin{equation}\label{eq:appr_acc_ode_sol}
\overline{w}_j = \frac{1}{i\theta}A(t)e^{ij\theta}(1-e^{-i\theta}),\quad
w_{\phf{j}} = N(t)e^{ij\theta},\quad\forall j.
\end{equation}
Here $A(\cdot)$ and $N(\cdot)$ are two scalar functions that solve the ODE system:
\begin{equation}\label{eq:appr_acc_amp_ode}
\begin{array}{>{\displaystyle}r>{\displaystyle}c>{\displaystyle}l}
A'(t) + ickN(t) &=& 0, \\
N'(t) + cka(\theta)A(t) + ckb(\theta)N(t) &=& 0,
\end{array}
\end{equation}
with initial condition $A(0) = N(0) = 1$.
The two functions $a(\theta)$ and $b(\theta)$ of (\ref{eq:appr_acc_amp_ode}) are:
\begin{equation}\label{eq:appr_acc_amp_aux}
a(\theta) \eqdef \frac{1-e^{-i\theta}}{i\theta^2}\sum_{l=-q+1}^qe^{il\theta}\alpha_l,\quad\textrm{ and }\quad
b(\theta) \eqdef \frac{1}{\theta}\sum_{l=-q}^qe^{il\theta}\beta_l.
\end{equation}
Comparing (\ref{eq:appr_acc_adv_sol}) and (\ref{eq:appr_acc_ode_sol}), we have the following definition of the order of accuracy of the semi-discretized FD-FV scheme.
\begin{definition}\label{def:appr_acc_fdfv}
The semi-discretized FD-FV scheme defined by $[\mathcal{D}_x]$ is at least $p^\powth\textrm{\sc-}$order accurate, if the following is true. 
If the method is applied to solve the Cauchy problem specified by (\ref{eq:appr_acc_adv}) with initial condition $w(x,0)=e^{ikx}$, in which case the approximated solutions are in the form of (\ref{eq:appr_acc_ode_sol}), then the following identities hold:
\begin{equation}\label{eq:appr_acc_fdfv_order}
A(t) = e^{-ickt}\left(1+c_a(t)\theta^p+O(\theta^{p+1})\right),\quad
N(t) = e^{-ickt}\left(1+c_n(t)\theta^p+O(\theta^{p+1})\right),
\end{equation}
where $c_a$ and $c_n$ are functions of time and they are independent of $\theta$.
\end{definition}
\begin{remark}\label{eq:appr_acc_fdfv_amp}
Although we express the two scalars $A(\cdot)$ and $N(\cdot)$ as functions of time, they actually also depend on the spatial parameter $\theta=hk$.
However, in our proof of spatial-order of accuracy in next theorem we suppose that the mesh size is fixed; 
thus we omit the dependence of these functions on $\theta$ for simplicity.
\end{remark}
The main result concerning the spatial-order of the semi-discretized FD-FV method and the designed-order of the operator $[\mathcal{D}_x]$ is given by the next theorem.
\begin{theorem}\label{thm:appr_acc_fdfv}
Given a $p^\powth\textrm{\sc-}$order accurate $[\mathcal{D}_x]$ such that $b_0\eqdef\sum_{l=-q}^q\beta_l\ne0$, it leads to a semi-discretized FD-FV method that is at least $(p+1)^\powth\textrm{\sc-}$order accurate provided also that $cb_0>0$, where $c$ is the constant advection speed of eq.~(\ref{eq:appr_acc_adv}).
\end{theorem}
The proof of this theorem is based on solving~(\ref{eq:appr_acc_amp_ode}) exactly and performing Taylor series expansions, which is provided as Appendix~\ref{app:order}.
The important implication of Theorem~\ref{thm:appr_acc_fdfv} is that by allowing a mixed form of cell-averages and nodal values to design the differential operator $[\mathcal{D}_x]$, one achieves one-order higher accuracy compared to conventional finite difference methods.
To see why this is true, we note that conventional finite difference operators fit into the hybrid FD-FV framework by setting all $\alpha_l$ to be zero (Remark~\ref{rk:appr_fdm}), in which case $b_0=0$.
By slight modification of the proof of Appendix~\ref{app:order} to $b_0=0$, it is easy to see that the maximum possible order of accuracy is $p^{\powth}$, i.e., the designed order of the differential operator.
It is also noted that the upwind-like condition $cb_0>0$ is very important to prevent a second mode from blowing up as the mesh size goes to zero, see eq.~(\ref{eq:order_sol_bdd}) and the discussion that follows.

To conclude this section, we listed in Table~\ref{tb:appr_ddo} several operators $[\mathcal{D}_x]$ as well as their designed order and corresponding $b_0$.
\begin{table}
\caption{FD-FV operators $[\mathcal{D}_x]$ with $b_0\ne0$.}
\label{tb:appr_ddo}
\begin{tabular}{|c|ll|c|}
\hline
Order & & $\left[\mathcal{D}_xw\right]_{j+1/2}$ & $b_0$\\ \hline
1st & Forward: & $2\left(\overline{w}_{j+1}-w_{j+1/2}\right)/h$ & $-2$\\
 & Backward: & $2\left(w_{j+1/2}-\overline{w}_j\right)/h$ & $2$ \\ \hline 
2nd & Forward: & $2\left(-w_{j+3/2}+3\overline{w}_{j+1}-2w_{j+1/2}\right)/h$ & $-6$ \\
 & Backward: & $2\left(2w_{j+1/2}-3\overline{w}_j+w_{j-1/2}\right)/h$ & $6$ \\ \hline
3rd & Forward: & $\left(\overline{w}_{j+2}-8w_{j+3/2}+17\overline{w}_{j+1}-10w_{j+1/2}\right)/(2h)$ & $-9$\\ 
 & F-biased: & $\left(-2w_{j+3/2}+7\overline{w}_{j+1}-4w_{j+1/2}-\overline{w}_j\right)/(2h)$ & $-3$\\
 & B-biased: & $\left(\overline{w}_{j+1}+4w_{j+1/2}-7\overline{w}_j+2w_{j-1/2}\right)/(2h)$ & $3$\\
 & Backward: & $\left(10w_{j+1/2}-17\overline{w}_j+8w_{j-1/2}-\overline{w}_{j-1}\right)/(2h)$ & $9$\\ \hline
4th & Forward: & $\left(-2w_{j+5/2}+7\overline{w}_{j+2}-16w_{j+3/2}+23\overline{w}_{j+1}-12w_{j+1/2}\right)/(2h)$ & $-15$\\
 & F-biased: & $\left(\overline{w}_{j+2}-12w_{j+3/2}+31\overline{w}_{j+1}-18w_{j+1/2}-2\overline{w}_j\right)/(6h)$ & $-5$\\
 & B-biased: & $\left(2\overline{w}_{j+1}+18w_{j+1/2}-31\overline{w}_j+12w_{j-1/2}-\overline{w}_{j-1}\right)/(6h)$ & $5$\\
 & Backward: & $\left(12w_{j+1/2}-23\overline{w}_j+16w_{j-1/2}-7\overline{w}_{j-1}+2w_{j-3/2}\right)/(2h)$ & $15$ \\ \hline
\end{tabular}
\end{table}
In the table, we use ``F-biased'' and ``B-biased'' to denote forward biased or backward biased stencils, respectively.
Due to the upwind-like condition $cb_0>0$, we also use the terminology ``upwind'' (designated by the superscript $^{up}$) or ``upwind-biased'' (designated by the superscripts $^{up-b}$ or $^{up-biased}$) to denote an operator.
For example, $[\mathcal{D}_x^{2,up}]$ denotes the second-order backward operator if $c>0$, and the second-order forward operator if $c<0$.

\subsection{Linear stability analysis}
\label{sec:appr_stab}
We study the linear stability of both the semi-discretized and fully-discretized FD-FV methods in this section.
In the semi-discretized case, we apply the von Neumann analysis to the same advection problem as in previous section.
In the view of (\ref{eq:appr_acc_adv_sol}) and (\ref{eq:appr_acc_ode_sol}), we can plot the dissipation error and the dispersion error by comparing $A(t)$ and $N(t)$ with $e^{-ickt}$, respectively.
For simplicity, we suppose $c>0$ and only consider operators of Table~\ref{tb:appr_ddo} with positive $b_0$ (Theorem~\ref{thm:appr_acc_fdfv}).

For a simple wave problem, the numerical solutions to the semi-discretized system is obtained by solving~(\ref{eq:appr_acc_amp_ode}), which gives rise to two eigenvalues $\lambda_1$ and $\lambda_2$;
these values correspond to the two superposed traveling waves in the numerical solutions $A(t)$ and $N(t)$, as demonstrated by~(\ref{eq:order_sol}).
Thus contrary to the von Neumann analysis to conventional finite difference systems, we study the dispersion property of the FD-FV methods by plotting $\tilde{\theta}_{1,2}$ against $\theta=kh$, where $\tilde{\theta}_{1,2}\eqdef \tilde{k}_{1,2}h$ and $\tilde{k}_{1,2}$ are the wave numbers of the corresponding traveling wave, i.e., $\Im(k\lambda_{1,2})$, the imaginary part of $k\lambda_{1,2}$.
Similarly, we study the dissipation property by plotting the numerical dissipation rate $\tilde{\varepsilon}_{1,2}$ against $\varepsilon\equiv0$, where $\tilde{\varepsilon}_{1,2}$ are the real parts of $\theta\lambda_{1,2}$, $\Re(\theta\lambda_{1,2})$, and it needs to be non-negative for linear stability.
To this end, we can rewrite~(\ref{eq:order_sol}) in the form:
\begin{displaymath}
A(t) = a_1e^{-(\tilde{\varepsilon}_1+ic\tilde{\theta}_1)t/h} + a_2e^{-(\tilde{\varepsilon}_2+ic\tilde{\theta}_2)t/h}\ \textrm{ and }\ 
N(t) = n_1e^{-(\tilde{\varepsilon}_1+ic\tilde{\theta}_1)t/h} + n_2e^{-(\tilde{\varepsilon}_2+ic\tilde{\theta}_2)t/h},
\end{displaymath}
where $a_{1,2}$ and $n_{1,2}$ are constants that are independent of $t$.
These quantities are important for the accuracy analysis, as shown in Appendix~\ref{app:order}, but they do not affect the linear stability of the semi-discretized FD-FV schemes.

\begin{figure}\centering
\begin{subfigure}[b]{0.45\textwidth}\centering
    \includegraphics[width=\textwidth]{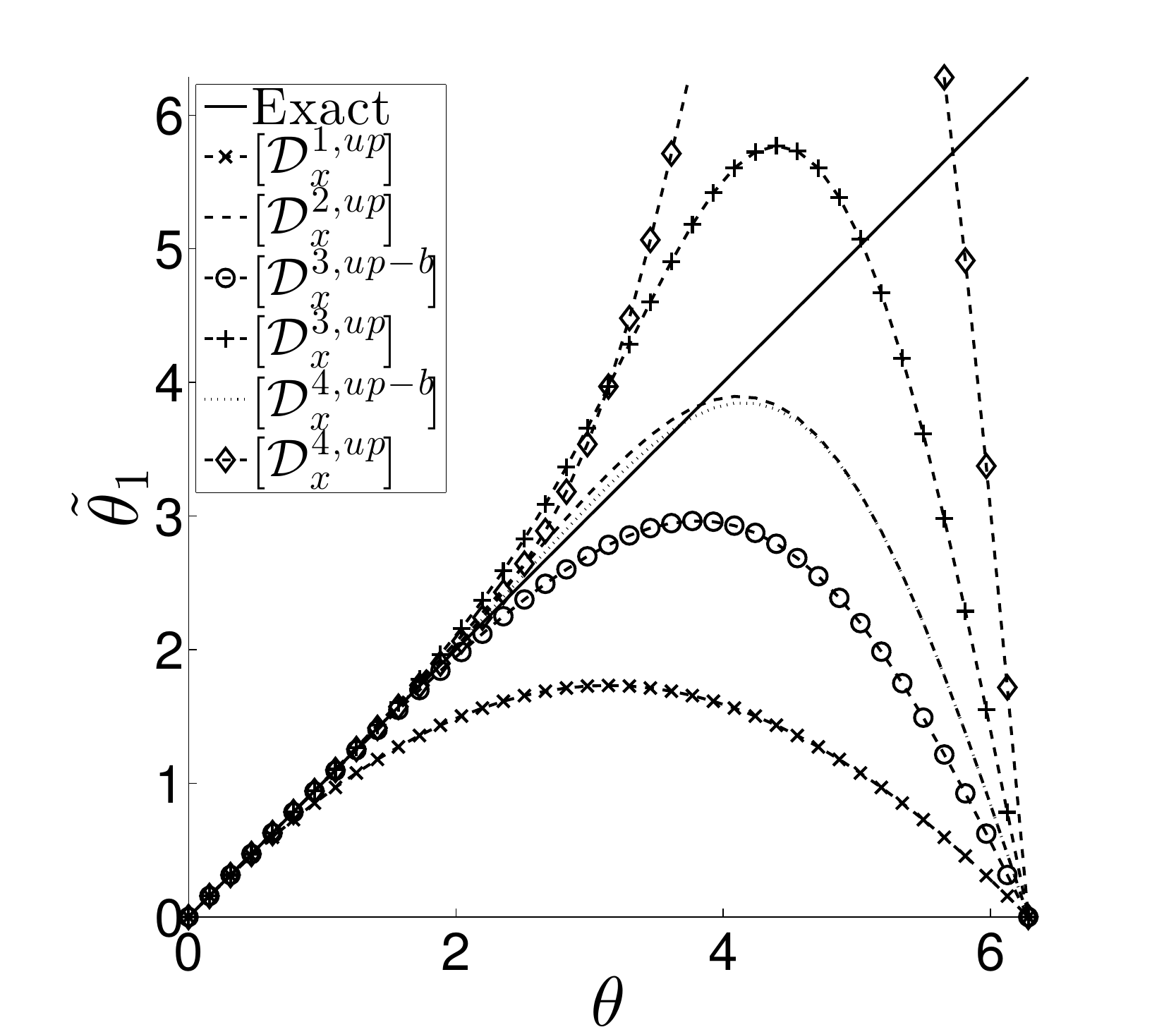}
    \caption{$\tilde{\theta}_1$ vs. $\theta$.}
    \label{fg:appr_stab_disper_pos}
\end{subfigure}
\begin{subfigure}[b]{0.45\textwidth}\centering
    \includegraphics[width=\textwidth]{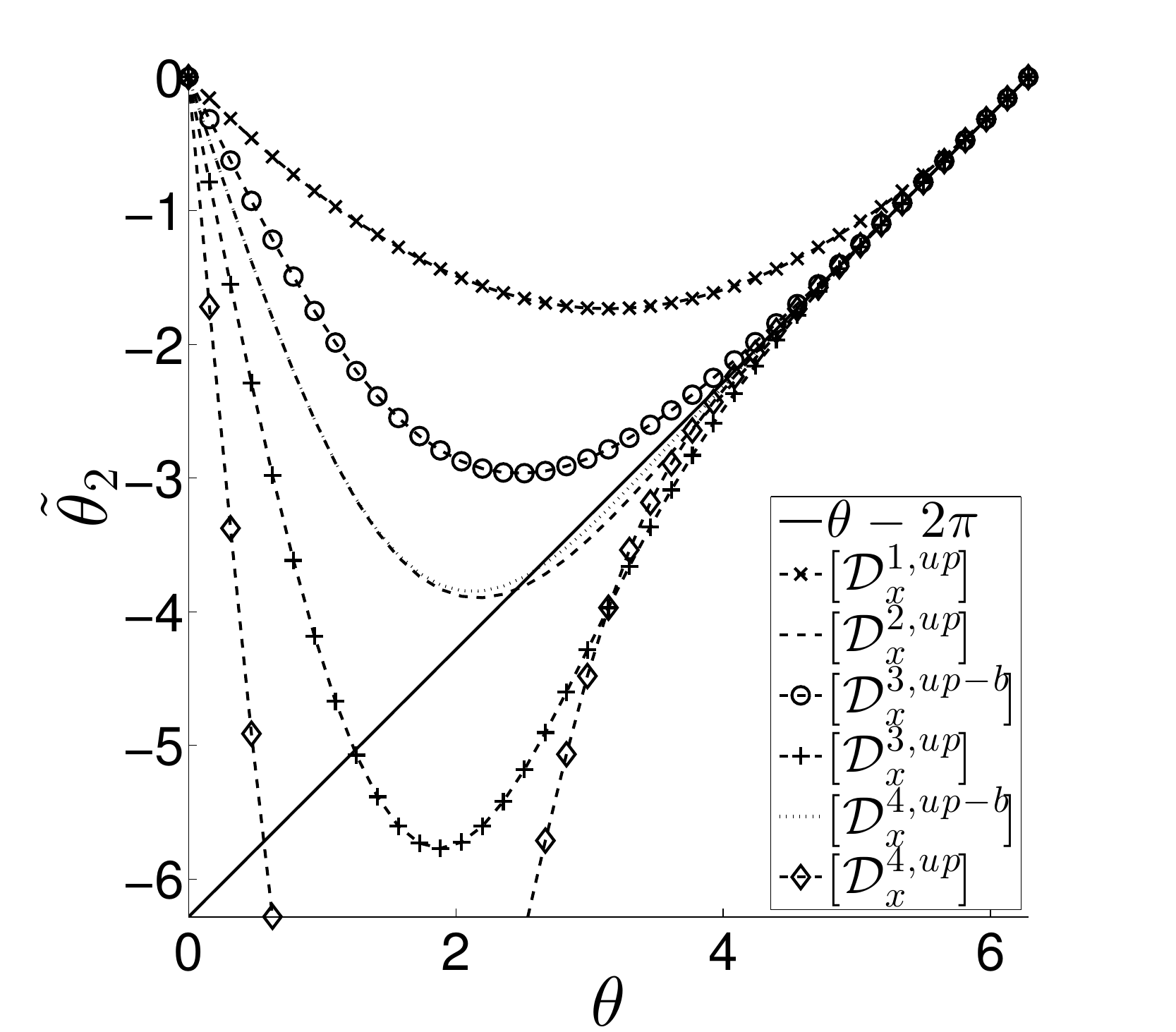}
    \caption{$\tilde{\theta}_2$ vs. $\theta$.}
    \label{fg:appr_stab_disper_neg}
\end{subfigure}
\caption{Dispersion property of various FD-FV operators $[\mathcal{D}_x]$: (\ref{fg:appr_stab_disper_pos}) $\theta\Im(\lambda_1)$ vs. $\theta$, (\ref{fg:appr_stab_disper_neg}) $\theta\Im(\lambda_2)$ vs. $\theta$.
Here ``$^{up}$'' and ``$^{up-biased}$'' denote the upwind and the upwind-biased stencils, respectively.}
\label{fg:appr_stab_disper}
\end{figure}

\begin{figure}\centering
\begin{subfigure}[b]{0.45\textwidth}\centering
    \includegraphics[width=\textwidth]{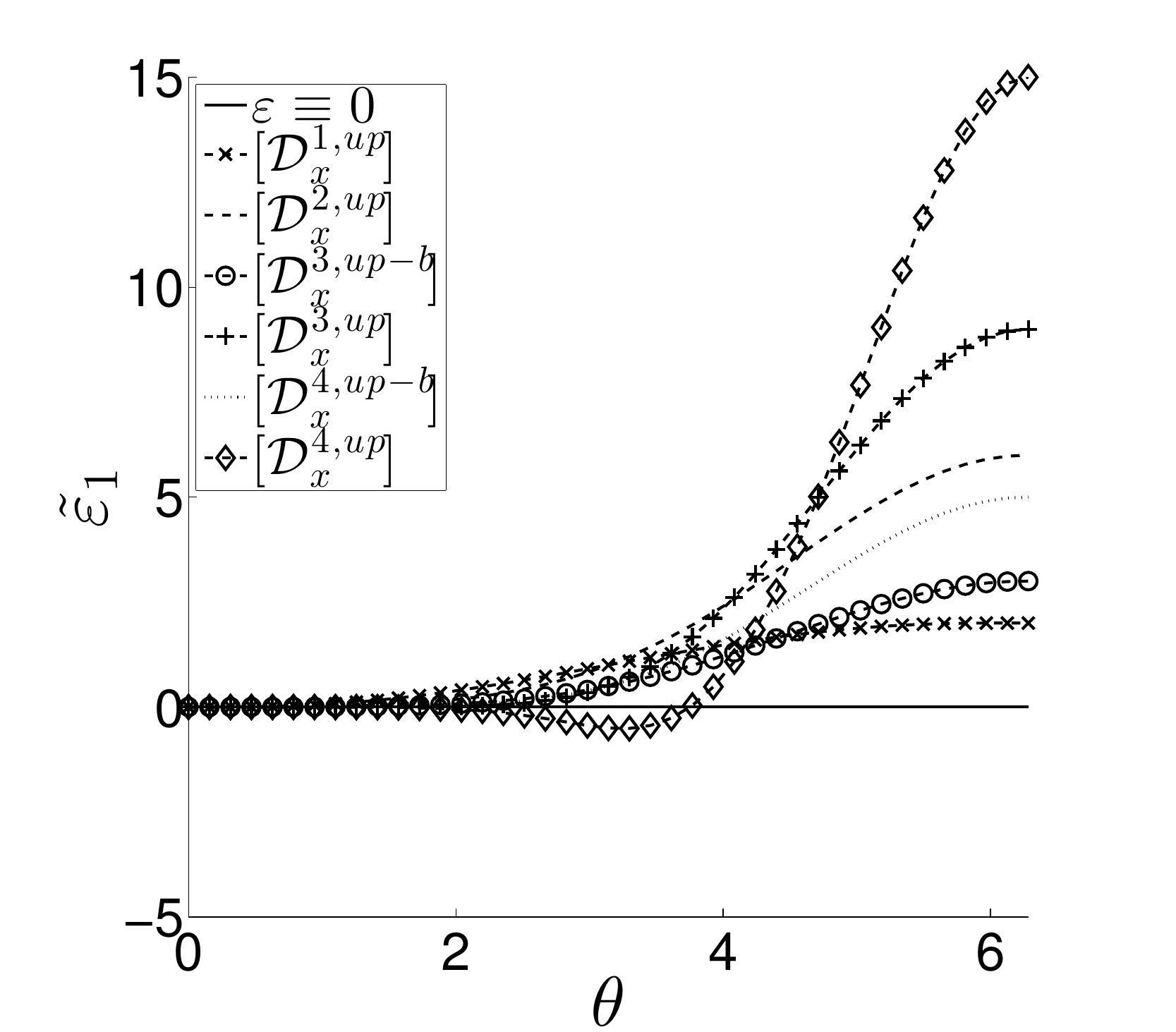}
    \caption{$\tilde{\varepsilon}_1$ vs. $\theta$.}
    \label{fg:appr_stab_dissip_pos}
\end{subfigure}
\begin{subfigure}[b]{0.45\textwidth}\centering
    \includegraphics[width=\textwidth]{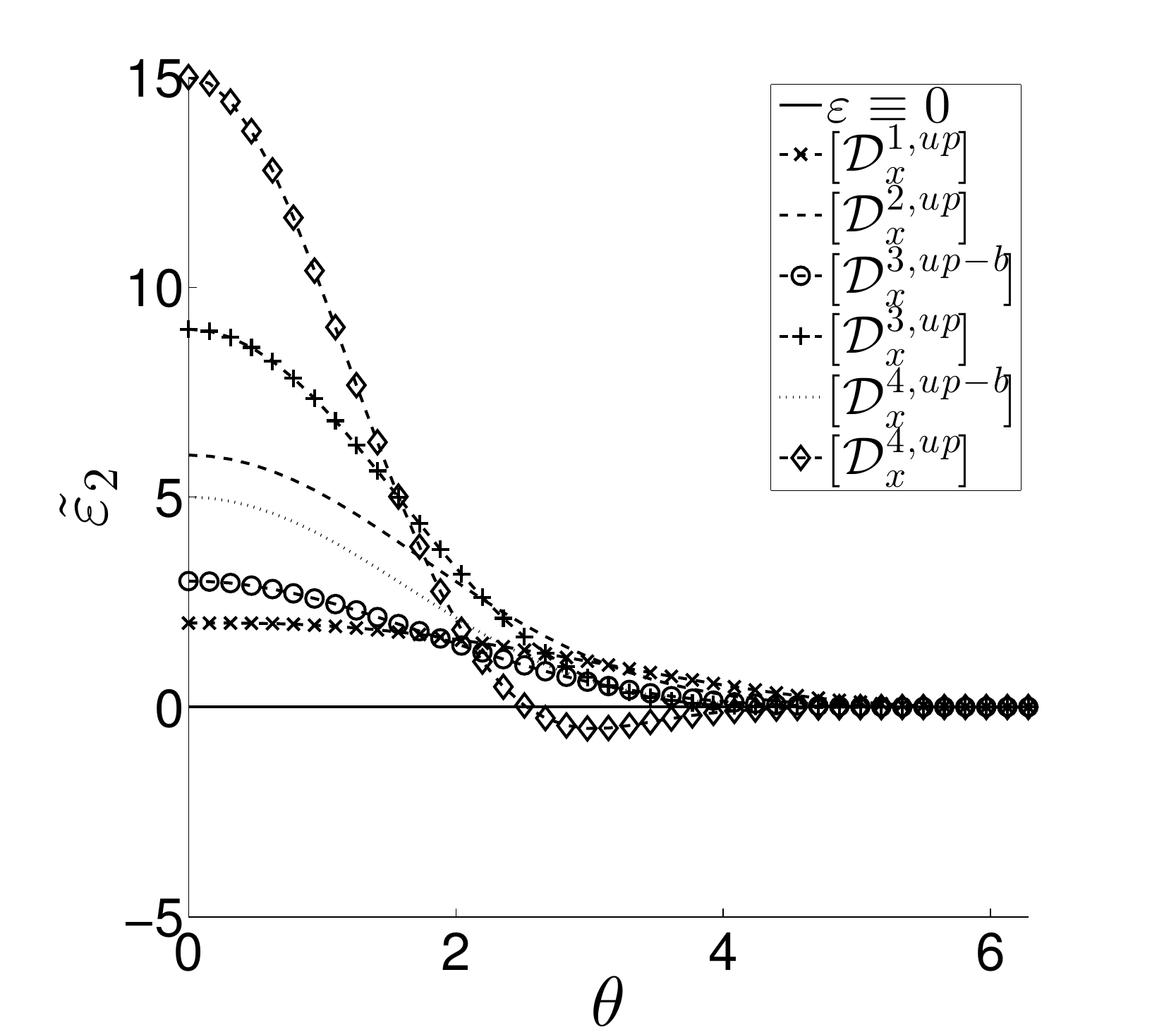}
    \caption{$\tilde{\varepsilon}_2$ vs. $\theta$.}
    \label{fg:appr_stab_dissip_neg}
\end{subfigure}
\caption{Dissipation property of various FD-FV operators $[\mathcal{D}_x]$: (\ref{fg:appr_stab_dissip_pos}) $\theta\Re(\lambda_1)$ vs. $\theta$, (\ref{fg:appr_stab_dissip_neg}) $\theta\Re(\lambda_2)$ vs. $\theta$.}
\label{fg:appr_stab_dissip}
\end{figure}

In Figure~\ref{fg:appr_stab_disper} we plot the dispersion relations in the range $\theta\in[0,2\pi]$, in which the right end of the interval corresponds to resolving the wave length by one cell.
Clearly, we see that the wave speed associated with $\lambda_1$ is closer to the true wave speed as the order of the operator increases; 
whereas the wave associated with $\lambda_2$ travels in the opposite direction as the physical wave.
Note, however, the second wave has a magnitude diminishing very fast as $h\to0$, see the accuracy analysis in Appendix~\ref{app:order}.
Similarly, the dissipation properties of the FD-FV operators are plotted in the same range in Figure~\ref{fg:appr_stab_dissip}.
From the figures, the dissipation rate of $[\mathcal{D}_x^{4,up}]$ is positive around $\theta\approx\pi$, or two cells per wave length.
Thus the method is linearly unstable in this regime on the semi-discretized level.
This observation is consistent with a class of leapfrog methods~\cite{AIserles:1986a}, which claims that the number of upwind points in the stencil cannot exceed the number of downwind points by more than two.
Furthermore, the dispersion and dissipation relations associated with $\lambda_{1,2}$ show certain symmetry properties about $\theta=\pi$; 
this is further discussed in Remark~\ref{rk:order_sym} of Appendix~\ref{app:order}.
\begin{remark}\label{rk:appr_stab_alt}
An alternative way to investigate the linear stability of the FD-FV methods is to write the semi-discretized system in matrix form, and plot the eigenvalue distribution of the resulting linear ODE system.
This approach not only applies to periodic boundary conditions but also works well with other types of boundary conditions, such as the examples given in Section~\ref{sec:appr_bc}.
\end{remark}

Finally, we consider the fully-discretized FD-FV methods by pairing each $[\mathcal{D}_x]$ operator with an appropriate time-integrator.
For simplicity, only explicit Runge-Kutta type methods are considered in this paper.
In the view of Theorem~\ref{thm:appr_acc_fdfv}, for each $[\mathcal{D}_x]$ that is a $p^\powth\textrm{\sc-}$order accuracy operator we choose a $(p+1)^\powth\textrm{\sc-}$order accurate scheme to integrate the resulting ODE system in time.
In particular, we consider the $2^\pownd\textrm{\sc-}$order and $3^\powrd\textrm{\sc-}$order total variational diminishing Runge-Kutta schemes by S. Gottlieb et al.~\cite{SGottlieb:1998a}, ``the'' $4^\powth\textrm{\sc-}$order method as described by E. Hairer et al.~\cite{EHairer:1993a}, and the $5^\powth\textrm{\sc-}$order method by H. A. Luther~\cite{HALuther:1966a}.
These methods are denoted by $\oname{RK2}, \oname{RK3}, \oname{RK4}$, and $\oname{RK5}$, respectively;
and their corresponding FD-FV operators $[\mathcal{D}_x]$ as well as the maximum allowed Courant number (denoted by $\nu_{\max}$) are summarized in Table~\ref{tb:appr_stab_cfl}.
\begin{table}
\caption{Pair of FD-FV operator and Runge-Kutta time-integrator, and the corresponding maximum allowed Courant numbers.}
\label{tb:appr_stab_cfl}
\begin{tabular}{|c|c|c|c|}
\hline
Operator & Time-integrator & Formal order of the FD-FV method & $\nu_{\max}$ \\ \hline
$[\mathcal{D}_x^{1,up}]$ & $\oname{RK2}$ & $2$ & $1.0$ \\ \hline
$[\mathcal{D}_x^{2,up}]$ & $\oname{RK3}$ & $3$ & $0.409$ \\ \hline
$[\mathcal{D}_x^{3,up-biased}]$ & $\oname{RK4}$ & $4$ & $0.808$ \\ \hline
$[\mathcal{D}_x^{3,up}]$ & $\oname{RK4}$ & $4$ & $0.309$ \\ \hline
$[\mathcal{D}_x^{4,up-biased}]$ & $\oname{RK5}$ & $5$ & $0.494$ \\ \hline
\end{tabular}
\end{table}

\subsection{Extension to systems of equations}
\label{sec:appr_sys}
In this section, we extend the methods presented before to a system of hyperbolic conservation laws.
Considering the governing equation~(\ref{eq:appr_1dhcl}), we define the nodal variable $\bs{w}_{j+1/2}$ and the cell average $\overline{\bs{w}}_j$ similarly as (\ref{eq:appr_gen_n}--\ref{eq:appr_gen_c}).
Then in analogy to (\ref{eq:appr_gen_c_eq}), the semi-discretized formula for $\overline{\bs{w}}_j$ is given by:
\begin{equation}\label{eq:appr_sys_c_eq}
\frac{d\overline{\bs{w}}_j}{dt} + \frac{1}{h}\left[\bs{f}\left(\bs{w}_{j+1/2}\right)-\bs{f}\left(\bs{w}_{j-1/2}\right)\right] = 0.
\end{equation}
To construct the spatial discretization for nodal values, we use the hyperbolic assumption and compute the characteristic decomposition of the flux Jacobian $\bs{J} = \bs{R}\bs{\Lambda}\bs{L}$, where $\bs{J}\eqdef\partial\bs{f}/\partial\bs{w}$.
We denote the diagonal entries of the diagonal matrix $\bs{\Lambda}$ by $\lambda_1\le\lambda_2\le\cdots\le\lambda_d$, and let $\bs{r}_k$ and $\bs{l}_k$ be the right and left eigenvalues associated with $\lambda_k, 1\le k\le d$, respectively, such that $\bs{R} = [\bs{r}_1,\ \bs{r}_2,\ \cdots,\ \bs{r}_d]$ and $\bs{R}^{-1} = \bs{L} = [\bs{l}_1,\ \bs{l}_2,\ \cdots,\ \bs{l}_d]^T$.
Using these notations, equation~(\ref{eq:appr_1dhcl}) is rewritten locally as decoupled scalar problems as:
\begin{displaymath}
\bs{L}\bs{w}_t + \bs{\Lambda}\bs{L}\bs{w}_x = 0.
\end{displaymath}
To this end, the semi-discretized formula for $\bs{w}_{j+1/2}$ is given by:
\begin{displaymath}
\frac{d(\bs{L}\bs{w}_{j+1/2})}{dt} + \bs{\Lambda}[\mathcal{D}_x(\bs{L}\bs{w})]_{j+1/2} = 0\ \textrm{ or }\ 
\frac{d(\omega_{k,j+1/2})}{dt} + \lambda_k[\mathcal{D}_x\omega]_{j+1/2} = 0,\ k=1,\cdots,q,
\end{displaymath}
where $\bs{L}$ and $\bs{\Lambda}$ are computed using flux Jacobian $\bs{J}(\bs{w}_{j+1/2})$, and $\omega_{k,j+1/2}\eqdef\bs{l}_k^T\bs{w}_{j+1/2}$.
Thus by using a particular FD-FV operator, such as $[\mathcal{D}_x^{1,up}]$, to solve the hyperbolic system~(\ref{eq:appr_1dhcl}), we mean to apply the scalar operator to each characteristics.
\begin{example}\label{ex:appr_sys}
Consider the hyperbolic equation:
\begin{displaymath}
\left[\begin{array}{c}u\\v\end{array}\right]_t + \left[\begin{array}{c}v\\u\end{array}\right]_x = 0\quad\Longrightarrow\quad
\left[\begin{array}{c}u+v\\u-v\end{array}\right]_t + \left[\begin{array}{cc}1&0\\0&-1\end{array}\right]\left[\begin{array}{c}u+v\\u-v\end{array}\right]_x = 0,
\end{displaymath}
the spatial discretization for $(u_{j+1/2},v_{j+1/2})$ using $[\mathcal{D}_x^{1,up}]$ is given by:
\begin{displaymath}
\left\{\begin{array}{>{\displaystyle}l}
\frac{d(u_{j+1/2}+v_{j+1/2})}{dt} + \frac{2}{h}\left[(u_{j+1/2}+v_{j+1/2})-(\overline{u}_j+\overline{v}_j)\right] = 0 \\
\frac{d(u_{j+1/2}-v_{j+1/2})}{dt} - \frac{2}{h}\left[(\overline{u}_{j+1}-\overline{v}_{j+1})-(u_{j+1/2}-v_{j+1/2})\right] = 0 
\end{array}\right..
\end{displaymath}
\end{example}

\subsection{Boundary conditions}
\label{sec:appr_bc}
To conclude this section, we describe the numerical treatments of the Dirichlet boundary condition for incoming characteristics, and briefly study the linear stability property by applying the methods to solve simple model problems.
For simplicity, we consider the advection problem:
\begin{displaymath}
w_t + cw_x = 0,\quad (x,t)\in[0,x_0]\times[0,T],\quad\textrm{ where }\quad c>0,\quad x_0\gg0,
\end{displaymath}
with the Dirichlet boundary condition:
\begin{displaymath}
w(0,t) = g(t)
\end{displaymath}
for some known data $g(t),\ t\in[0,T]$.
Then we follow the next two-step procedure to enforce this boundary condition in the spatial discretization of the FD-FV method given by $[\mathcal{D}_x^p]$, where $p$ denotes the order of this operator.
Supposing $x_{1/2}=0$, we first set $w_{1/2} = g$.
Next, for all small indices $j$ such that calculating $[\mathcal{D}_x^pw]_{j+1/2}$ requires data collocated to the left of $x=0$, we use an operator $[\mathcal{D}_x^qw]$ of order $q$ and with a biased stencil at this node, such that the corresponding stencil stays within $[0,x_0]$.
In order not to degrade the overall order of accuracy of the method, $q$ should be at least $p\!-\!1$~\cite{BGustafsson:2008a}.
\begin{example}\label{ex:appr_bc_d}
Consider the FD-FV method given by $[\mathcal{D}_x^{3,up}]$ (i.e., the $3^\powrd\textrm{\sc-}$order backward operator of Table~\ref{tb:appr_ddo}), then applying this operator at $x_{1/2}$ and $x_{3/2}$ requires a stencil that penetrates into outside the computational domain.
To enforce the boundary condition at $x=0$, we set $w_{1/2}=g$ and apply $[\mathcal{D}_x^{2,up}]$ or $[\mathcal{D}_x^{3,up-biased}]$ at $x_{3/2}$ to discretize the equation for the nodal variable in space.
\end{example}

Although the hyperbolic problems are in general not well-posed with Neumann boundary conditions, in practice we may use the numerical boundary condition of Neumann type for nonlinear hyperbolic systems.
For example, we consider the 1D Euler equations:
\begin{equation}\label{eq:appr_bc_1deuler}
\left[\begin{array}{c}
\rho \\ \rho u \\ E
\end{array}\right]_t + 
\left[\begin{array}{c}
\rho u \\ \rho u^2 + p \\ (E+p)u
\end{array}\right]_x = 0,
\end{equation}
which is defined on $(x,t)\in[0,x_0]\times[0,T]$.
Here $\rho$, $u$, and $p$ are density, velocity, and pressure, respectively; and $E$ is the total energy per unit volume that is determined by the ideal gas law: $E = p/(\gamma-1)+\rho u^2/2$ with $\gamma\equiv1.4$.
Then the wall-boundary condition at $x=0$:
\begin{displaymath}
u(0,t)\equiv 0
\end{displaymath}
leads to the Neumann boundary condition for the total energy at this location: $E_x(0,t)\equiv0$.
To see this, we first have $u(0,t)=0$ and $u_t(0,t)=0$, then from the momentum and energy equations we obtain:
\begin{displaymath}
0 = (\rho u)_t + (\rho u^2+p)_x = p_x\quad\Longrightarrow\quad
E_x = \frac{p_x}{\gamma-1}+\frac{1}{2}\rho_xu^2+\rho uu_x = 0\quad\textrm{ at }\quad x=0.
\end{displaymath}
Hence we describe here the numerical Neumann boundary condition, such as $w_x(0,t) = g(t)$ at the left boundary $x_{1/2}=0$.
This procedure is also composed of two steps: First we compute $w_{1/2}$ such that $[\mathcal{D}_x^{q,forward}]_{1/2}=g$, then for all small indices $j$ such that calculating $[\mathcal{D}_x^pw]_{j+1/2}$ requires data collocated to the left of $x_{1/2}$ we use the same strategy as in the Dirichlet case.
Similarly as before, $q$ should be at least $p\!-\!1$ in order not to degrade the formal order of accuracy.
\begin{example}\label{ex:appr_bc_n}
Consider again the FD-FV method given by $[\mathcal{D}_x^{3,up}]$ and homogeneous boundary condition $w_x=0$, then we first compute $w_{1/2}$ by setting $[\mathcal{D}_x^{3,forward}w]_{1/2}=0$ to have:
\begin{displaymath}
w_{1/2} = \frac{1}{10}(\overline{w}_2^n-8w_{3/2}^n+17\overline{w}_1).
\end{displaymath}
Next we apply $[\mathcal{D}_x^{2,up}]$ or $[\mathcal{D}_x^{3,up-biased}]$ at $x_{3/2}$ as in Example~\ref{ex:appr_bc_d}.
\end{example}

\begin{figure}\centering
\begin{subfigure}[b]{0.48\textwidth}\centering
    \includegraphics[width=\textwidth]{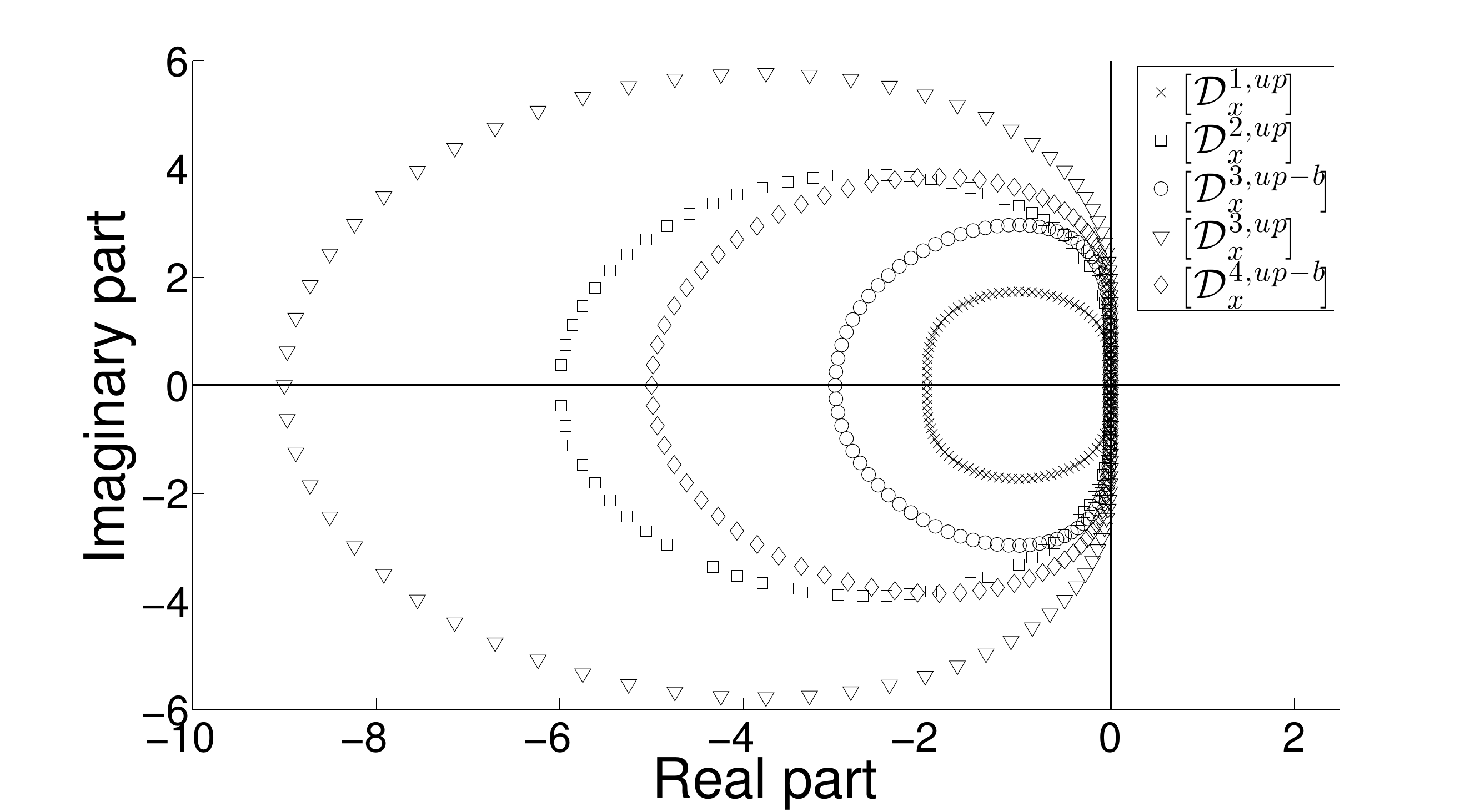}
    \caption{Periodic boundary condition.}
    \label{fg:appr_bc_spec_period}
\end{subfigure}
\begin{subfigure}[b]{0.48\textwidth}\centering
    \includegraphics[width=\textwidth]{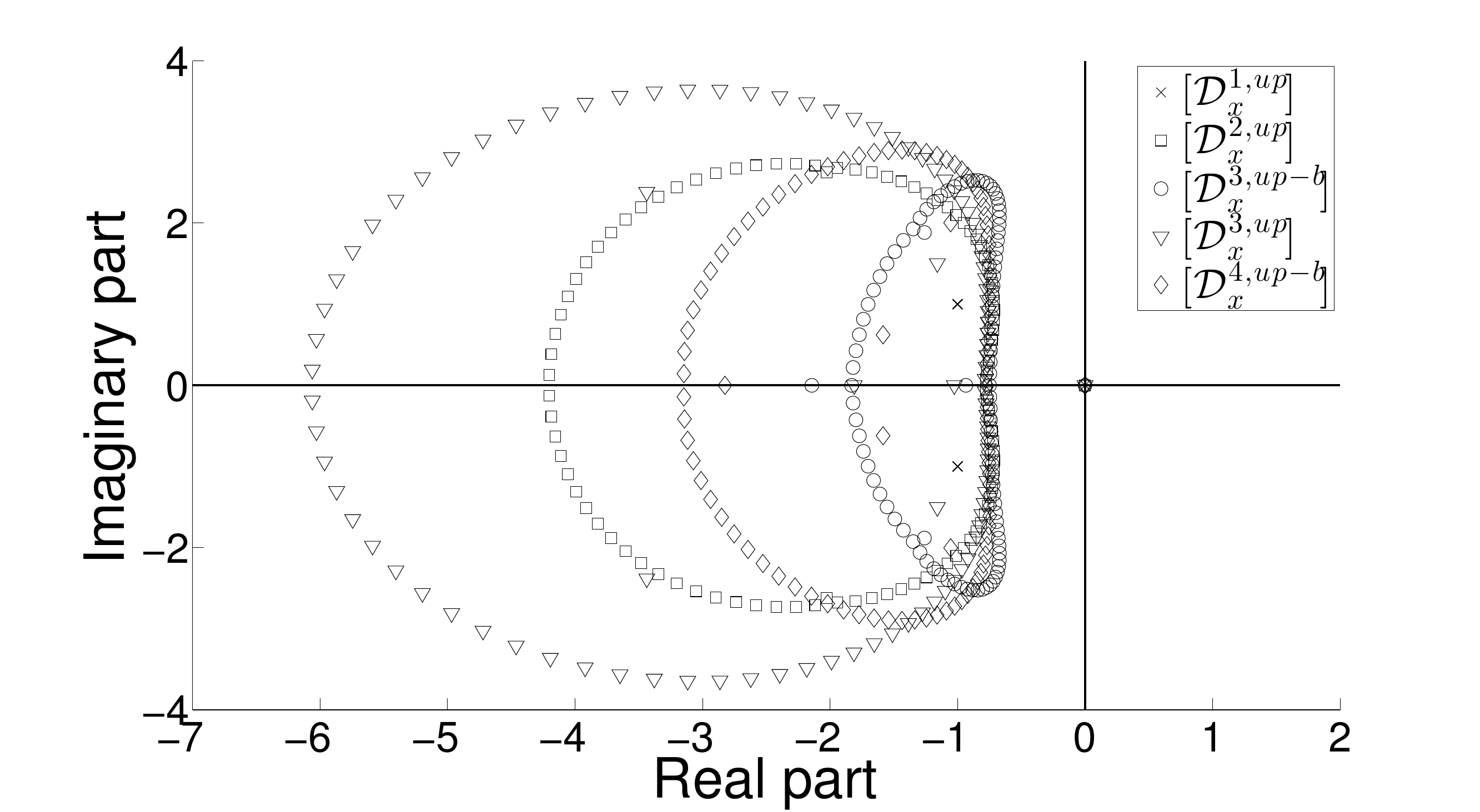}
    \caption{Dirichlet boundary condition.}
    \label{fg:appr_bc_spec_dirich}
\end{subfigure}
\caption{The spectra of various FD-FV schemes to solve $w_t+w_x=0$ on $x\in[0,1]$ using $N=50$ cells: (\ref{fg:appr_bc_spec_period}) periodic boundary condition, and (\ref{fg:appr_bc_spec_dirich}) Dirichlet boundary condition at $x=0$.}
\label{fg:appr_bc_spec}
\end{figure}

Finally, we study the spectral property of the resulting ODE systems obtained by applying the FD-FV methods to solve a model problem with various boundary conditions.
In particular, we solve the advection problem $w_t+w_x=0$ on the domain $x\in[0,1]$ by meshing the latter into $N=50$ uniform cells.
In the first scenario, we consider the periodic boundary condition $w(0,t)=w(1,t)$, and denote by $W$ the array of unknowns:
\begin{displaymath}
\bs{W} = [\overline{w}_1,\ \overline{w}_2,\ \cdots,\ \overline{w}_N,\ w_{1/2},\ w_{3/2},\ \cdots,\ w_{N-1/2}]^T.
\end{displaymath}
Then the spatial discretization leads to a ODE system for $\bs{W}$:
\begin{equation}\label{eq:appr_bc_ode}
\bs{W}_t = \frac{1}{h}\bs{D}\bs{W},\quad h = \frac{1}{N}.
\end{equation}
Here $\bs{D}$ is a $2N\times 2N$ matrix that is determined by the operator $[\mathcal{D}_x]$.
For example, the matrix $\bs{D}$ corresponding to $[\mathcal{D}_x^{1,up}]$ is given by:
\begin{displaymath}
\bs{D} = \left[\begin{array}{cc}
\bs{0} & \bs{I}-\bs{L}_1 \\
2\bs{L}_2 & -2\bs{I}
\end{array}\right],\quad
\setlength{\arraycolsep}{3pt} 
\bs{L}_1 = \left[\begin{array}{ccccc}
0 & 1 & \cdots & 0 & 0 \\
0 & 0 & \cdots & 0 & 0 \\
\vdots & \vdots & \ddots & \vdots & \vdots \\
0 & 0 & \cdots & 0 & 1 \\
1 & 0 & \cdots & 0 & 0 
\end{array}\right]\quad\textrm{ and }\quad
\bs{L}_2 = \left[\begin{array}{ccccc}
0 & 0 & \cdots & 0 & 1 \\
1 & 0 & \cdots & 0 & 0 \\
\vdots & \vdots & \ddots & \vdots & \vdots \\
0 & 0 & \cdots & 0 & 0 \\
0 & 0 & \cdots & 1 & 0 
\end{array}\right],
\setlength{\arraycolsep}{1pt} 
\end{displaymath}
with $\bs{0}$ and $\bs{I}$ being the $N\times N$ zero matrix and identity matrix, respectively.
The semi-discretization of the FD-FV method is linearly stable if and only if all the eigenvalues of the matrix $\bs{D}$ of (\ref{eq:appr_bc_ode}) stay within the left complex domain (i.e., the real parts are non-positive).
To this end, we plot the eigenvalue distributions of the linearly stable operators in Section~\ref{sec:appr_stab}, as shown in Figure~\ref{fg:appr_bc_spec_period}.
Clearly, all the eigenvalues are in the left half complex domain, indicating that the methods are linearly stable.

In the second scenario, we consider the Dirichlet boundary condition at $x_{1/2}=0$, namely $w(0,t)=g(t)$; and we denote by $\bs{W}$ the array of $2N-1$ unknowns:
\begin{displaymath}
\bs{W} = [\overline{w}_1,\ \overline{w}_2,\ \cdots,\ \overline{w}_N,\ w_{3/2},\ \cdots,\ w_{N-1/2}]^T.
\end{displaymath}
Then applying the various FD-FV methods to solve this problem yields the ODE system:
\begin{displaymath}
\bs{W}_t = \frac{1}{h}\bs{D}\bs{W} + \bs{g},
\end{displaymath}
where $\bs{g}$ designates the contribution from boundary conditions.
The eigenvalue distributions computed using several FD-FV operators are presented in Figure~\ref{fg:appr_bc_spec_dirich}, which again indicates linear stability for all operators.
Note that we use $[\mathcal{D}_x^{3,up-biased}]$ at $x_{3/2}$ and use $[\mathcal{D}_x^{3,up}]$ at the last node $x_{N+1/2}$ in the cases of FD-FV methods that are defined by $[\mathcal{D}_x^{3,up}]$ and $[\mathcal{D}_x^{4,up-biased}]$.


\section{Extension to two space dimensions}
\label{sec:ext}
The extension of the preceding methodology to two space dimensions is complicated by the fact that for a two-dimensional cell, the face-averaged flux cannot be evaluated exactly by using nodal values as in one-dimensional case.
In particular, consider the general first-order conservation law written in Cartesian coordinates:
\begin{equation}\label{eq:ext_2dhcl}
\bs{w}_t + \bs{f}(\bs{w})_x + \bs{g}(\bs{w})_y = 0,
\end{equation}
where $\bs{w}(x,y,t)$ takes value in $\mathbb{R}^d$, and $\bs{f}, \bs{g}: \mathbb{R}^d\to\mathbb{R}^d$ are both continuously differentiable.
Let the rectangular computational domain $\Omega$ be divided into cells with non-overlapping interiors: $\Omega = \cup_{j_1,j_2}\mathcal{C}_{j_1,j_2}$, where $\mathcal{C}_{j_1,j_2} = [x_{\mhf{j_1}},x_{\phf{j_1}}]\times[y_{\mhf{j_2}},y_{\phf{j_2}}]$; 
and similar to the 1D case, we consider uniform grids for simplicity and let $x_{\phf{j_1}}-x_{\mhf{j_1}}\equiv h_1$, $y_{\phf{j_2}}-y_{\mhf{j_2}}\equiv h_2$, and define $x_{j_1} = (x_{\mhf{j_1}}+x_{\phf{j_1}})/2$, $y_{j_2} = (y_{\mhf{j_2}}+y_{\phf{j_2}})/2$.
Clearly, the cell-averaged value is defined such that:
\begin{equation}\label{eq:ext_c}
\overline{\bs{w}}_{j_1,j_2}^n \approx \frac{1}{h_1h_2}\int_{\mathcal{C}_{j_1,j_2}}\bs{w}(x,y,t^n)dxdy.
\end{equation}
However, integrating~(\ref{eq:ext_2dhcl}) over $\mathcal{C}_{j_1,j_2}$ yields:
\begin{displaymath}
\begin{array}{>{\displaystyle}r>{\displaystyle}l}
 &\frac{d}{dt}\left(\frac{1}{h_1h_2}\int_{C_{j_1,j_2}}\bs{w}(x,y,t)dxdy\right) \\
+&\frac{1}{h_1}\left[\frac{1}{h_2}\int_{y_{j_2-1/2}}^{y_{j_2+1/2}}\bs{f}\left(\bs{w}(x_{j_1+1/2},y,t)\right)dy-
                     \frac{1}{h_2}\int_{y_{j_2-1/2}}^{y_{j_2+1/2}}\bs{f}\left(\bs{w}(x_{j_1-1/2},y,t)\right)dy\right] \\
+&\frac{1}{h_2}\left[\frac{1}{h_1}\int_{x_{j_1-1/2}}^{x_{j_1+1/2}}\bs{g}\left(\bs{w}(x,y_{j_2+1/2},t)\right)dx-
                     \frac{1}{h_1}\int_{x_{j_1-1/2}}^{x_{j_1+1/2}}\bs{g}\left(\bs{w}(x,y_{j_2-1/2},t)\right)dx\right] = 0,
\end{array}
\end{displaymath}
which contains the face-averaged fluxes that are not available directly.
In order to obtain a high-order scheme, these fluxes need to be evaluated to a certain order of accuracy.
One way to do this is to define the face-averaged flux itself as a variable, however it requires additional governing equations to describe the dynamics of this variable, which increase the complexity of the methodology significantly.
Instead, we choose a more direct extension of the 1D methodology: we select several nodes on a particular cell face and approximate the face-averaged flux by a weighted sum of the fluxes evaluated at these nodes.
That is, we first choose a quadrature rule for integration on the interval $[-1,1]$ with quadrature points $-1\le q_1,<\cdots<q_m\le 1$ and corresponding weights $\omega_1,\cdots,\omega_m>0$ such that $\omega_1+\cdots+\omega_m=2$.
Then the nodal variables on the edge $[x_{\mhf{j_1}},x_{\phf{j_1}}]\times\{y_{\phf{j_2}}\}$ are collocated at $(x_{j_1+q_kh_1/2},y_{\phf{j_2}}), k=1,\cdots,m$; and the face-averaged flux is:
\begin{displaymath}
\frac{1}{h_1}\int_{x_{\mhf{j_1}}}^{x_{\phf{j_1}}}\bs{g}(\bs{w}(x,y_{\phf{j_2}},t)dx\approx\frac{1}{2}\sum_{k=1}^m\omega_k\bs{g}(\bs{w}(x_{j_1+q_kh_1/2},y_{\phf{j_2}},t)).
\end{displaymath}
Similarly, the nodal variables on the edge $\{x_{\phf{j_1}}\}\times[y_{\mhf{j_2}},y_{\phf{j_2}}]$ are collocated at the points $(x_{\phf{j_1}},y_{j_2+q_kh_2/2}), k=1,\cdots,m$; and the face-averaged flux is: 
\begin{displaymath}
\frac{1}{h_2}\int_{y_{\mhf{j_2}}}^{y_{\phf{j_2}}}\bs{f}(\bs{w}(x_{\phf{j_1}},y,t)dy\approx\frac{1}{2}\sum_{k=1}^m\omega_k\bs{f}(\bs{w}(x_{\phf{j_1}},y_{\phf{j_2}+q_kh_2/2},t)).
\end{displaymath}

By doing so, we lose the exactness of the cell-boundary fluxes as in the 1D case, hence the order of this approximation must be chosen carefully without degrading the designed order of the FD-FV scheme.
In particular, we use at least $p^\powth\textrm{\sc-}$order quadrature rule to design a FD-FV scheme with the same order.

For simplicity, we only consider two 1D FD-FV operators, namely $[\mathcal{D}_x^{1,up}]$ and $[\mathcal{D}_x^{2,up}]$, and extend them to 2D problems. 
The extensions for the other operators, however, follow the same procedure.
The common component of these FD-FV methods is the semi-discretized formula for the cell-averages, namely:
\begin{equation}\label{eq:ext_semi_c}
\frac{d\overline{\bs{w}}_{j_1,j_2}}{dt} + 
\frac{\sum_{k=1}^m\omega_k\left(\bs{f}_{\phf{j_1},j_2;k}-\bs{f}_{\mhf{j_1},j_2;k}\right)}{2h_1} + 
\frac{\sum_{k=1}^m\omega_k\left(\bs{g}_{j_1,\phf{j_2};k}-\bs{g}_{j_1,\mhf{j_2};k}\right)}{2h_2} = 0, 
\end{equation}
where we define $\bs{x}_{\phf{j_1},j_2;k} \eqdef (x_{\phf{j_1}},y_{j_2}+q_kh_2/2)$, $\bs{x}_{j_1,\phf{j_2};k} \eqdef (x_{j_1}+q_kh_1/2,y_{\phf{j_2}})$,
and for any such subscript combination $\bs{x}_\ast$ we denote the corresponding nodal value and nodal flux by $\bs{w}_\ast = \bs{w}(\bs{x}_\ast,t)$ and $\bs{f}(\bs{w}_\ast)$ or $\bs{g}(\bs{w}_\ast)$, respectively.
Note that there is no ambiguity on which index should the last subscript $k$ corresponds to, because it always acts on the index that is not a half integer.
We shall also use the notations $x_{j_1;k} = x_{j_1}+q_kh_1/2$ and $y_{j_2;k} = y_{j_2}+q_kh_2/2$ in subsequent sections.

With (\ref{eq:ext_semi_c}) established, we focus on the semi-discretization formula for nodal values in the following sections.
For simplicity, we restrict ourself first to the case of simple scalar advection equation:
\begin{equation}\label{eq:ext_2dadv}
w_t + c_1w_x + c_2w_y = 0
\end{equation}
for some constants $c_1>0$ and $c_2>0$.
Extension to general conservation systems is similar to that in Section~\ref{sec:appr_sys}, and the details will be provided at the end of this section.

\subsection{A second-order scheme based on $[\mathcal{D}_x^{1,up}]$}
\label{sec:ext_2nd}
To construct a second-order FD-FV scheme, we need at least second-order accurate quadrature rule.
In particular, we select $m=1$, $q_1=0$, and $\omega_1 = 2$.
In this case, we have exactly one node on each cell face; hence instead of using the notations $\bs{x}_{\phf{j_1},j_2;k}$ and $\bs{x}_{j_1,\phf{j_2};k}$ introduced before, we will write $\bs{x}_{\phf{j_1},j_2}$ and $\bs{x}_{j_1,\phf{j_2}}$ for simplicity.
The same rule holds for the other quantities associated with the nodal values.

The semi-discretized formula for a second-order FD-FV scheme using $[\mathcal{D}_x^{1,up}]$ is given by the following:
\begin{equation}\label{eq:ext_2nd_semi}
\begin{array}{>{\displaystyle}r>{\displaystyle}c>{\displaystyle}l}
\frac{d\overline{w}_{j_1,j_2}}{dt} + \frac{c_1w_{\phf{j_1},j_2}-c_1w_{\mhf{j_1},j_2}}{h_1} + \frac{c_2w_{j_1,\phf{j_2}}-c_1w_{j_1,\mhf{j_2}}}{h_2} & = & 0, \\
\frac{dw_{\phf{j_1},j_2}}{dt} + \frac{c_1(2w_{\phf{j_1},j_2}-2\overline{w}_{j_1,j_2})}{h_1} + \frac{c_2(w_{\phf{j_1},j_2}-w_{\phf{j_1},j_2-1})}{h_2}& = & 0, \\
\frac{dw_{j_1,\phf{j_2}}}{dt} + \frac{c_1(w_{j_1,\phf{j_2}}-w_{j_1-1,\phf{j_2}})}{h_1} + \frac{c_2(2w_{j_1,\phf{j_2}}-2\overline{w}_{j_1,j_2})}{h_2}& = & 0.
\end{array}
\end{equation}
In these equations, the one-dimensional FD-FV operator $[\mathcal{D}_x^{1,up}]$ is used to discretize the spatial derivative $\partial_xw_{\phf{j_1},j_2}$ in~(\ref{eq:ext_2nd_semi})$_2$; and a similar operator $[\mathcal{D}_y^{1,up}]$ is used to discretize $\partial_yw_{j_1,\phf{j_2}}$ in~(\ref{eq:ext_2nd_semi})$_3$.
The other two spatial derivatives in the last two equations of~(\ref{eq:ext_2nd_semi}) are discretized by a conventional finite-difference method, i.e., the first-order upwind order that is denoted by $\mathcal{D}_x^{1,up}$ or $\mathcal{D}_y^{1,up}$.
These notations differ from the FD-FV operators in the absence of the bracket.
These lead to the simplified semi-discretized formula for nodal values:
\begin{equation}\label{eq:ext_2nd_semi_n}
\begin{array}{>{\displaystyle}r>{\displaystyle}c>{\displaystyle}l}
\frac{dw_{\phf{j_1},j_2}}{dt} + c_1[\mathcal{D}_x^{1,up}w]_{\phf{j_1},j_2} + c_2\mathcal{D}_y^{1,up}w_{\phf{j_1},j_2} & = & 0, \\
\frac{dw_{j_1,\phf{j_2}}}{dt} + c_1\mathcal{D}_x^{1,up}w_{j_1,\phf{j_2}} + c_2[\mathcal{D}_y^{1,up}w]_{j_1,\phf{j_2}} & = & 0,
\end{array}
\end{equation}
where the exact formula for $\mathcal{D}_x^{1,up}$ and $\mathcal{D}_y^{1,up}$ are:
\begin{displaymath}
\mathcal{D}_y^{1,up}w_{\phf{j_1},j_2} = \frac{w_{\phf{j_1},j_2}-w_{\phf{j_1},j_2-1}}{h_1},\quad
\mathcal{D}_x^{1,up}w_{j_1,\phf{j_2}} = \frac{w_{j_1,\phf{j_2}}-w_{j_1-1,\phf{j_2}}}{h_2}.
\end{displaymath}

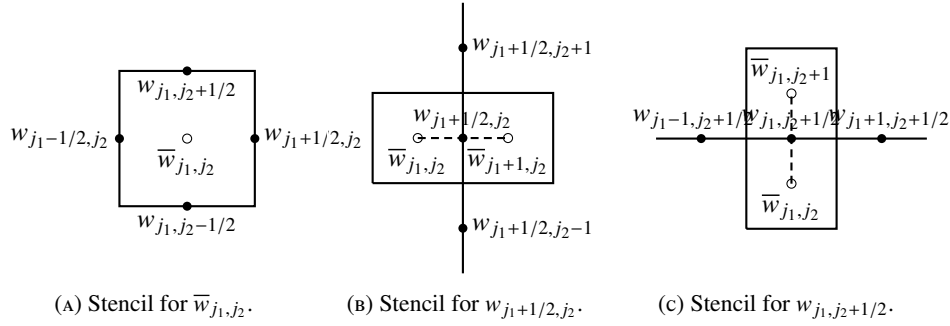
\begin{figure}\centering
\begin{subfigure}[b]{0.32\textwidth}\centering
    \begin{tikzpicture}
        \draw [color=white] (-1.8,-1.8) rectangle (1.8,1.8);
        \draw [line width=0.8] (-0.9,-0.9) rectangle (0.9,0.9);
        \draw (   0,   0) circle (0.06cm);
        \fill (   0, 0.9) circle (0.06cm);
        \fill (   0,-0.9) circle (0.06cm);
        \fill ( 0.9,   0) circle (0.06cm);
        \fill (-0.9,   0) circle (0.06cm);
        \draw (   0,   0) node [below] {$\overline{w}_{j_1,j_2}$};
        \draw (   0, 0.9) node [below] {$w_{j_1,\phf{j_2}}$};
        \draw (   0,-0.9) node [below] {$w_{j_1,\mhf{j_2}}$};
        \draw ( 0.9,   0) node [right] {$w_{\phf{j_1},j_2}$};
        \draw (-0.9,   0) node [left] {$w_{\mhf{j_1},j_2}$};
    \end{tikzpicture}
    \caption{Stencil for $\overline{w}_{j_1,j_2}$.}
    \label{fg:ext_2nd_semi_c}
\end{subfigure}
\begin{subfigure}[b]{0.32\textwidth}\centering
    \begin{tikzpicture}
        \draw [color=white] (-1.8,-1.8) rectangle (1.8,1.8);
        \draw [line width=0.8] (-1.2,-0.6) -- (1.2,-0.6) -- (1.2,0.6) -- (-1.2,0.6) -- (-1.2,-0.6);
        \draw [line width=0.8] (0,-1.8) -- (0,1.8);
        \draw [densely dashed, line width=0.8] (-0.6,0) -- (0.6,0);
        \draw (-0.6, 0) circle (0.06cm);
        \draw ( 0.6, 0) circle (0.06cm);
        \fill ( 0, 0) circle (0.06cm);
        \fill ( 0, 1.2) circle (0.06cm);
        \fill ( 0,-1.2) circle (0.06cm);
        \draw (-0.6, 0) node [below] {$\overline{w}_{j_1,j_2}$};
        \draw ( 0.6, 0) node [below] {$\overline{w}_{j_1+1,j_2}$};
        \draw ( 0, 0) node [above] {$w_{\phf{j_1},j_2}$};
        \draw ( 0, 1.2) node [right] {$w_{\phf{j_1},j_2+1}$};
        \draw ( 0,-1.2) node [right] {$w_{\phf{j_1},j_2-1}$};
    \end{tikzpicture}
    \caption{Stencil for $w_{\phf{j_1},j_2}$.}
    \label{fg:ext_2nd_semi_x}
\end{subfigure}
\begin{subfigure}[b]{0.32\textwidth}\centering
    \begin{tikzpicture}
        \draw [color=white] (-1.8,-1.8) rectangle (1.8,1.8);
        \draw [line width=0.8] (-1.8,0) -- (1.8,0);
        \draw [line width=0.8] (-0.6,-1.2) -- (-0.6,1.2) -- (0.6,1.2) -- (0.6,-1.2) -- (-0.6,-1.2);
        \draw [densely dashed, line width=0.8] (0,-0.6) -- (0,0.6);
        \draw ( 0,-0.6) circle (0.06cm);
        \draw ( 0, 0.6) circle (0.06cm);
        \fill ( 0, 0) circle (0.06cm);
        \fill ( 1.2, 0) circle (0.06cm);
        \fill (-1.2, 0) circle (0.06cm);
        \draw ( 0,-0.6) node [below] {$\overline{w}_{j_1,j_2}$};
        \draw ( 0, 0.6) node [above] {$\overline{w}_{j_1,j_2+1}$};
        \draw ( 0, 0) node [above] {$w_{j_1,j_2+1/2}$};
        \draw ( 1.3, 0) node [above] {$w_{j_1+1,j_2+1/2}$};
        \draw (-1.3, 0) node [above] {$w_{j_1-1,j_2+1/2}$};
    \end{tikzpicture}
    \caption{Stencil for $w_{j_1,\phf{j_2}}$.}
    \label{fg:ext_2nd_semi_y}
\end{subfigure}
\caption{Stencil of the semi-discretization formula~(\ref{eq:ext_2nd_semi}) for the (\ref{fg:ext_2nd_semi_c}) cell-averaged value, and (\ref{fg:ext_2nd_semi_x}--\ref{fg:ext_2nd_semi_y}) nodal values.
$\circ$ -- location of cell averages; $\bullet$ -- location of nodal values.}
\label{fg:ext_2nd_semi}
\end{figure}

To this end, the stencil of the semi-discretization formula for $\overline{w}_{j_1,j_2}$, $w_{\phf{j_1},j_2}$, and $w_{j_1,\phf{j_2}}$ given in (\ref{eq:ext_2nd_semi}) are illustrated in Figure~\ref{fg:ext_2nd_semi}.
Note that we also include the stencils for situations $c_1>0$ and/or $c_2>0$.
Similar to the one-dimensional case, although we use discrete differential operators of first-order spatial accuracy, the method~(\ref{eq:ext_2nd_semi}) leads to a second-order accurate method.
To show this, we consider the Cauchy problem of solving (\ref{eq:ext_2dadv}) on $\mathbb{R}^2\times[0,T]$ using a uniform mesh with cell size $h_1\times h_2$.
In particular, let the initial condition be a simple planar wave of the form:
\begin{equation}\label{eq:ext_2nd_wave}
w(x,y,0) = e^{i(k_1x+k_2y)},
\end{equation}
Then we formalize the statement that (\ref{eq:ext_2nd_semi}) is second-order accurate in space by the following statement:
\begin{conjecture}\label{con:ext_acc_2nd}
Fixing any $k_1$ and $k_2$, then applying the method~(\ref{eq:ext_2nd_semi}) to solve the preceding Cauchy problem with exact initial data~(\ref{eq:ext_2nd_wave}) leads to semi-discretized solutions $\overline{w}_{j_1,j_2}(T)$, $w_{\phf{j_1},j_2}(T)$ and $w_{j_1,\phf{j_2}}(T)$, which differ from the exact values by $O(h^2)$.
\end{conjecture}
Here we assume that $O(h_1)=O(h_2)=O(h)$ always holds, which is particularly true if we keep the aspect ratio of each cell the same when refining the mesh.
This conjecture is based on a similar argument to the proof of Theorem~\ref{thm:appr_acc_fdfv}, and more details concerning why this conjecture could be true is provided in Appendix~\ref{app:2nd}.
Finally, we combine the semi-discretization~(\ref{eq:ext_2nd_semi}) and the second-order total variational diminishing Runge-Kutta scheme in time to obtain a globally second-order accurate method.
We omit the linear stability analysis of the resulting fully-discrete scheme, and mention that the stability limit given by the Courant-Friedrichs-Lewy condition is:
\begin{equation}\label{eq:ext_2nd_cfl}
\Delta t < \frac{\nu_{\max}}{c_1/h_1 + c_2/h_2},\quad \nu_{\max} = 1.0.
\end{equation}

\subsection{A third-order scheme based on $[\mathcal{D}_x^{2,up}]$}
\label{sec:ext_3rd}
\begin{figure}\centering
\begin{subfigure}[b]{0.32\textwidth}\centering
    \begin{tikzpicture}
        \draw [color=white] (-1.8,-1.8) rectangle (1.8,1.8);
        \draw [line width=0.8] (-0.9,-0.9) rectangle (0.9,0.9);
        \draw (0,0) circle (0.06cm);
        \fill ( 0.9, 0.52) circle (0.06cm);
        \fill ( 0.9,-0.52) circle (0.06cm);
        \fill (-0.9, 0.52) circle (0.06cm);
        \fill (-0.9,-0.52) circle (0.06cm);
        \fill ( 0.52, 0.9) circle (0.06cm);
        \fill (-0.52, 0.9) circle (0.06cm);
        \fill ( 0.52,-0.9) circle (0.06cm);
        \fill (-0.52,-0.9) circle (0.06cm);
        \draw (0,0.05) node [below] {$\overline{w}_{j_1,j_2}$};
        \draw ( 0.9,-0.52) node [right] {$w_{\phf{j_1},j_2;1}$};
        \draw ( 0.9, 0.52) node [right] {$w_{\phf{j_1},j_2;2}$};
        \draw (-0.67,-0.9) node [below] {$w_{j_1,\mhf{j_2};1}$};
        \draw ( 0.67,-0.9) node [below] {$w_{j_1,\mhf{j_2};2}$};
        
        \draw [line width=0.8, densely dashed] (-0.9,0.52) -- (0.9,0.52);
        \draw (-0.05,0.47) rectangle (0.05,0.57);
        \draw ( 0.0,1.2) node [above] {$\overline{w}_{j_1,j_2;2}^x$};
        \draw [line width=0.6,->] (0.0,1.2) -- (0.0, 0.6);
    \end{tikzpicture} 
    \caption{Stencil for $\overline{w}_{j_1,j_2}$.}
    \label{fg:ext_3rd_semi_c}
\end{subfigure}
\begin{subfigure}[b]{0.32\textwidth}\centering
    \begin{tikzpicture}
        \draw [color=white] (-1.8,-1.8) rectangle (1.8,1.8);
        \draw [line width=0.8] (-1.2,-0.6) -- (1.2,-0.6) -- (1.2,0.6) -- (-1.2,0.6) -- (-1.2,-0.6);
        \draw [line width=0.8] (0,-1.8) -- (0,1.8);
        \draw [densely dashed, line width=0.8] (-1.2, 0.346) -- ( 1.2, 0.346);
        \draw (-0.65, 0.296) rectangle (-0.55, 0.396);
        \draw ( 0.55, 0.296) rectangle ( 0.65, 0.396);
        \draw (-0.6, 0) circle (0.06cm);
        \draw ( 0.6, 0) circle (0.06cm);
        \fill ( 0, 0.346) circle (0.06cm);
        \fill ( 0,-0.346) circle (0.06cm);
        \fill (-1.2, 0.346) circle (0.06cm);
        \fill ( 1.2, 0.346) circle (0.06cm);
        \fill ( 0,-0.854) circle (0.06cm);
        \fill ( 0, 0.854) circle (0.06cm);
        \fill ( 0, 1.546) circle (0.06cm);
        \fill (-0.946, 0.6) circle (0.06cm);
        \fill (-0.346, 0.6) circle (0.06cm);
        \fill (-0.946,-0.6) circle (0.06cm);
        \fill (-0.346,-0.6) circle (0.06cm);
        \fill ( 0.946, 0.6) circle (0.06cm);
        \fill ( 0.346, 0.6) circle (0.06cm);
        \fill ( 0.946,-0.6) circle (0.06cm);
        \fill ( 0.346,-0.6) circle (0.06cm);
        \draw (-0.6, 0) node [below] {$\overline{w}_{j_1,j_2}$};
        \draw ( 0.6, 0) node [below] {$\overline{w}_{j_1+1,j_2}$};
        \draw (-0.8, 0.8) node [above] {$w_{\phf{j_1},j_2;2}$};
        \draw [->, line width=0.6] (-0.35,0.8) -- (-0.05,0.4);
    \end{tikzpicture}
    \caption{Stencil for $w_{\phf{j_1},j_2;2}$.}
    \label{fg:ext_3rd_semi_x}
\end{subfigure}
\begin{subfigure}[b]{0.32\textwidth}\centering
    \begin{tikzpicture}
        \draw [color=white] (-1.8,-1.8) rectangle (1.8,1.8);
        \draw [line width=0.8] (-1.8,0) -- (1.8,0);
        \draw [line width=0.8] (-0.6,-1.2) -- (-0.6,1.2) -- (0.6,1.2) -- (0.6,-1.2) -- (-0.6,-1.2);
        \draw [densely dashed, line width=0.8] (-0.346,-1.2) -- (-0.346,1.2);
        \draw (-0.396, 0.55) rectangle (-0.296, 0.65);
        \draw (-0.396,-0.65) rectangle (-0.296,-0.55);
        \draw ( 0,-0.6) circle (0.06cm);
        \draw ( 0, 0.6) circle (0.06cm);
        \fill (-0.346, 0) circle (0.06cm);
        \fill ( 0.346, 0) circle (0.06cm);
        \fill (-0.346, 1.2) circle (0.06cm);
        \fill (-0.346,-1.2) circle (0.06cm);
        \fill (-1.546, 0) circle (0.06cm);
        \fill (-0.854, 0) circle (0.06cm);
        \fill ( 0.854, 0) circle (0.06cm);
        \fill (-0.6, 0.254) circle (0.06cm);
        \fill (-0.6, 0.946) circle (0.06cm);
        \fill (-0.6,-0.254) circle (0.06cm);
        \fill (-0.6,-0.946) circle (0.06cm);
        \fill ( 0.6, 0.254) circle (0.06cm);
        \fill ( 0.6, 0.946) circle (0.06cm);
        \fill ( 0.6,-0.254) circle (0.06cm);
        \fill ( 0.6,-0.946) circle (0.06cm);
        \draw ( 0.1,-0.6) node [below] {$\overline{w}_{j_1,j_2}$};
        \draw ( 0.15, 0.6) node [below] {$\overline{w}_{j_1,j_2+1}$};
        \draw (-1.3, 0.2) node [above] {$w_{j_1,j_2+1/2;1}$};
        \draw [->, line width=0.6] (-0.9, 0.25) -- (-0.4,0.05);
    \end{tikzpicture}
    \caption{Stencil for $w_{j_1,\phf{j_2}};1$.}
    \label{fg:ext_3rd_semi_y}
\end{subfigure}
\caption{Stencil of the semi-discretization for: (\ref{fg:ext_3rd_semi_c}), the cell-averaged value; and (\ref{fg:ext_3rd_semi_x}--\ref{fg:ext_3rd_semi_y}), the nodal values.
$\circ$ -- location of cell averages; $\bullet$ -- location of nodal values; $\square$ -- reconstructed line-averages.}
\label{fg:ext_3rd_semi}
\end{figure}
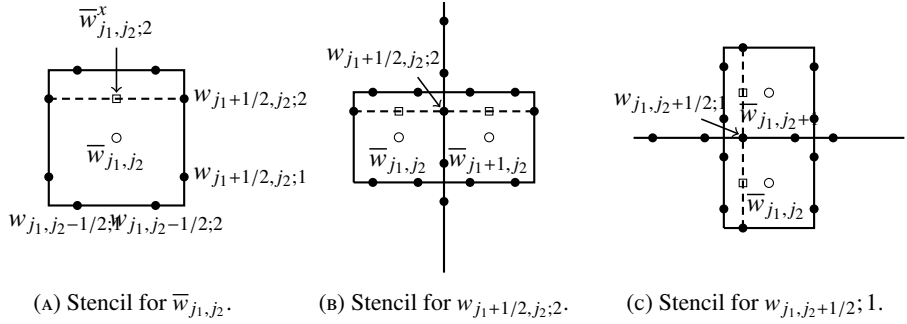

Next we consider constructing a third-order FD-FV method by utilizing the operator $[\mathcal{D}_x^{2,up}]$.
In this case, the quadrature rule to evaluate the edge-averaged flux needs to be at least third-order accurate.
To this end, we choose the Gauss-Legendre quadrature points with $m=2$, $q_{1,2} = \pm1/\sqrt{3}$, and $\omega_{1,2}=1$.
Then the stencil for the semi-discretization of a cell-averaged value $\overline{w}_{j_1,j_2}$ is arranged as Figure~\ref{fg:ext_3rd_semi_c}, and the formula is:
\begin{equation}\label{eq:ext_3rd_semi_c}
\begin{array}{>{\displaystyle}r>{\displaystyle}c>{\displaystyle}l}
\frac{d\overline{w}_{j_1,j_2}}{dt} + 
\frac{c_1(w_{\phf{j_1},j_2;1}+w_{\phf{j_1},j_2;2})-c_1(w_{\mhf{j_1},j_2;1}+w_{\mhf{j_1},j_2;2})}{2h_1} \\
+ \frac{c_2(w_{j_1,\phf{j_2};1}+w_{j_1,\phf{j_2};2})-c_2(w_{j_1,\mhf{j_2};1}+w_{j_1,\mhf{j_2};2})}{2h_2} 
& = & 0, \\
\end{array}
\end{equation}
The semi-discretization formula for the nodal values, however, needs more elaboration.
For example, considering apply $[\mathcal{D}_x^{2,up}]$ to the nodal value at $\bs{x}_{\phf{j_1},j_2;2}$ (Figure~\ref{fg:ext_3rd_semi_c}), we need the line-averaged data for the dashed line in the figure.
In particular, we need an approximation of the quantity:
\begin{displaymath}
\overline{w}^x_{j_1,j_2;2}(t) \approx \frac{1}{h_1}\int_{x_{\mhf{j_1}}}^{x_{\phf{j_1}}}w(x,y_{j_2;2},t)dx
\end{displaymath}
to at least third-order accuracy in space.
For this purpose, we omit the variable $t$ and define the variable $\overline{w}_{j_1}^x(y)$ as:
\begin{displaymath}
\overline{w}_{j_1}^x(y) \eqdef \frac{1}{h_1}\int_{x_{\mhf{j_1}}}^{x_{\phf{j_1}}}w(x,y)dx,
\end{displaymath}
and our target is to find the approximations $\overline{w}^x_{j_1,j_2;k}\approx\overline{w}_{j_1}^x(y_{j_2;k}),\ k=1,2$.
Due to our choice of the quadrature rule to locate the nodal values, $\overline{w}^x_{j_1}(y_{j_1\pm1/2})$ can be evaluated to the third-order accuracy by:
\begin{displaymath}
\overline{w}^x_{j_1}(y_{j_2\pm1/2}) = \frac{1}{2}w_{j_1,j_2\pm1/2;1} + \frac{1}{2}w_{j_1,j_2\pm1/2;2} + O(h_1^3).
\end{displaymath}
In addition, seeing that:
\begin{displaymath}
\frac{1}{h_1h_2}\int_{C_{j_1,j_2}}w(x,y)dxdy = \frac{1}{h_2}\int_{y_{\mhf{j_2}}}^{y_{\phf{j_1}}}\overline{w}^x_{j_1}(y)dy,
\end{displaymath}
in which the left hand side is represented by $\overline{w}_{j_1,j_2}$, we obtain a third-order accurate formula for $\overline{w}^x_{j_1,j_2;k}$ for $k=1,2$, namely:
\begin{displaymath}
\begin{array}{>{\displaystyle}r>{\displaystyle}c>{\displaystyle}l}
\overline{w}^x_{j_1,j_2;k} 
&=& \left(\frac{3q_k^2}{4}+\frac{q_k}{2}-\frac{1}{4}\right)\frac{w_{j_1,j_2+1/2;1}+w_{j_1,j_2+1/2;2}}{2} + \\
& & \left(\frac{3q_k^2}{4}-\frac{q_k}{2}-\frac{1}{4}\right)\frac{w_{j_1,j_2-1/2;1}+w_{j_1,j_2-1/2;2}}{2} + \left(\frac{3}{2}-\frac{3q_k^2}{2}\right)\overline{w}_{j_1,j_2}.
\end{array}
\end{displaymath}
This enables us to construct the discrete differential formula $[\mathcal{D}_x^{2,up}w]_{\phf{j_1},j_2;k}$ by:
\begin{displaymath}
[\mathcal{D}_x^{2,up}w]_{\phf{j_1},j_2;k} = \frac{2}{h_1}\left(
        2w_{\phf{j_1},j_2;2}-3\overline{w}^x_{j_1,j_2;k}+w_{\mhf{j_1},j_2;k}
        \right),
\end{displaymath}
and similarly the formula for $[\mathcal{D}_y^{2,up}w]_{j_1,\phf{j_2};k},\ k=1,2$.
Finally, for the other two derivatives $\partial_yw_{\phf{j_1},j_2;k}$ and $\partial_xw_{j_1,\phf{j_2};k}$, we use the similar strategy as in Section~\ref{sec:ext_2nd} and use the second-order upwind finite-difference operators along the corresponding grid lines.
We denote these discrete operators by $\mathcal{D}_y^{2,up}w_{\phf{j_1},j_2;k}$ and $\mathcal{D}_x^{2,up}w_{j_1,\phf{j_2};k}$, respectively, and write the semi-discretized equations for all the nodal values as:
\begin{equation}\label{eq:ext_3rd_semi_n}
\begin{array}{>{\displaystyle}r>{\displaystyle}c>{\displaystyle}l}
\frac{dw_{\phf{j_1},j_2;k}}{dt} + c_1[\mathcal{D}_x^{2,up}w]_{\phf{j_1},j_2;k} + c_2\mathcal{D}_y^{2,up}w_{\phf{j_1},j_2;k} & = & 0, \\
\frac{dw_{j_1,\phf{j_2};k}}{dt} + c_1\mathcal{D}_x^{2,up}w_{j_1,\phf{j_2};k} + c_2[\mathcal{D}_y^{2,up}w]_{j_1,\phf{j_2};k} & = & 0.
\end{array}
\end{equation}
The equations~(\ref{eq:ext_3rd_semi_c}) and (\ref{eq:ext_3rd_semi_n}) constitute the whole set of semi-discretization formula for our third-order FD-FV scheme to solve~(\ref{eq:ext_2dadv});
and the stencil for each cell-averaged value and nodal value (also taking into account of the possibility $c_1<0$ and/or $c_2<0$) are illustrated in Figure~\ref{fg:ext_3rd_semi}, in which the square designates the intermediate line-averaged quantities.

In principal, the method presented in this section could be extended to construct FD-FV scheme of arbitrary order of accuracy; 
however, how to proof that the formal order of accuracy of the resulting FD-FV scheme is one-order higher than the order of the discrete differential operators $[\mathcal{D}_{x,y}]$ and $\mathcal{D}_{x,y}$ that are used to construct the method is not clear.
Nevertheless, we make the conjecture that this is true, and verify it by numerical examples in next section.

Finally, we combine the spatial discretization scheme presented with the third-order total variational diminishing Runge-Kutta scheme in time to obtain a globally third-order accurate method, and use the Courant-Friedrichs-Lewy condition:
\begin{equation}\label{eq:ext_3rd_cfl}
\Delta t < \frac{\nu_{\max}}{c_1/h_1 + c_2/h_2},\quad \nu_{\max} = 0.306
\end{equation}
in the numerical examples.

\subsection{Extension to system of equations}
\label{sec:ext_sys}
The methods in previous sections can be written in a more general form for scalar problem:
\begin{displaymath}
w_t + f(w)_x + g(w)_y = 0
\end{displaymath}
as follows:
\begin{equation}\label{eq:ext_sys_semi}
\begin{array}{>{\displaystyle}r>{\displaystyle}c>{\displaystyle}l}
\frac{d\overline{w}_{j_1,j_2}}{dt} + 
\frac{\sum_k\omega_k(f_{\phf{j_1},j_2;k}-f_{\mhf{j_1},j_2;k})}{h_1} +
\frac{\sum_k\omega_k(g_{j_1,\phf{j_2};k}-g_{j_1,\mhf{j_2};k})}{h_2} & = & 0, \\
\frac{dw_{\phf{j_1},j_2;k}}{dt} + 
f'(w_{\phf{j_1},j_2;k})[\mathcal{D}_xw]_{\phf{j_1},j_2;k} + 
g'(w_{\phf{j_1},j_2;k})\mathcal{D}_yw_{\phf{j_1},j_2;k} & = & 0, \\ 
\frac{dw_{j_1,\phf{j_2};k}}{dt} + 
f'(w_{j_1,\phf{j_2};k})\mathcal{D}_xw_{j_1,\phf{j_2};k} + 
g'(w_{j_1,\phf{j_2};k})[\mathcal{D}_yw]_{j_1,\phf{j_2};k} & = & 0,
\end{array}
\end{equation}
where $[\mathcal{D}_x]$ and $[\mathcal{D}_y]$ denote one-dimensional FD-FV operators, and $\mathcal{D}_x$ and $\mathcal{D}_y$ denote conventional finite-difference operators.

The method~(\ref{eq:ext_sys_semi}) could be extended to the system~(\ref{eq:ext_2dhcl}) by using the characteristic decomposition of the Jacobian matrices:
\begin{displaymath}
\pp{\bs{f}}{\bs{w}} = \bs{J}_f = \bs{R}_f\bs{\Lambda}_f\bs{L}_f,\quad
\pp{\bs{g}}{\bs{w}} = \bs{J}_g = \bs{R}_g\bs{\Lambda}_g\bs{L}_g.
\end{displaymath}
And we have:
\begin{displaymath}
\begin{array}{>{\displaystyle}r>{\displaystyle}c>{\displaystyle}l}
\frac{d\overline{\bs{w}}_{j_1,j_2}}{dt} + 
\frac{\sum_k\omega_k(\bs{f}_{\phf{j_1},j_2;k}-\bs{f}_{\mhf{j_1},j_2;k})}{h_1} +
\frac{\sum_k\omega_k(\bs{f}_{j_1,\phf{j_2};k}-\bs{g}_{j_1,\mhf{j_2};k})}{h_2} & = & 0, \\
\frac{d\bs{w}_{\phf{j_1},j_2;k}}{dt} + 
\bs{R}_f\bs{\Lambda}_f[\mathcal{D}_x(\bs{L}_fw)]_{\phf{j_1},j_2;k} + 
\bs{R}_g\bs{\Lambda}_g\mathcal{D}_y(\bs{L}_gw)_{\phf{j_1},j_2;k} & = & 0, \\ 
\frac{d\bs{w}_{j_1,\phf{j_2};k}}{dt} + 
\bs{R}_f\bs{\Lambda}_f\mathcal{D}_x(\bs{L}_fw)_{j_1,\phf{j_2};k} + 
\bs{R}_g\bs{\Lambda}_g[\mathcal{D}_y(\bs{L}_gw)]_{j_1,\phf{j_2};k} & = & 0,
\end{array}
\end{displaymath}
where for simplicity we omit the subscripts on all the matrices associated with the Jacobian matrices.

\section{Numerical assessment}
\label{sec:num}
We apply the FD-FV methods developed in preceding sections to solve various problems in one and two space dimensions.
We focus on problems with solutions that are sufficiently smooth in the paper, but show one simple example with non-convex flux function and discontinuous initial condition to demonstrate the potential of the FD-FV methods to capture correct solution structure for general fluxes.
The latter issue, however, is not fully investigated and is left for future work.

\subsection{One-dimensional problems}
\label{sec:num_1d}
In this section we focus on the smooth problems, such as the advection equation and Euler equations that admit smooth solutions.
We also solve the Buckley-Leverett equation, an example of non-convex flux, with discontinuous initial data.
Note that although our methods are constructed based on Taylor series expansions, they perform well in solving this problem and capture the nonlinear waves accurately.

\subsubsection{One-dimensional advection equations}
\label{sec:num_1dadv}
\begin{figure}\centering
    \includegraphics[width=\textwidth]{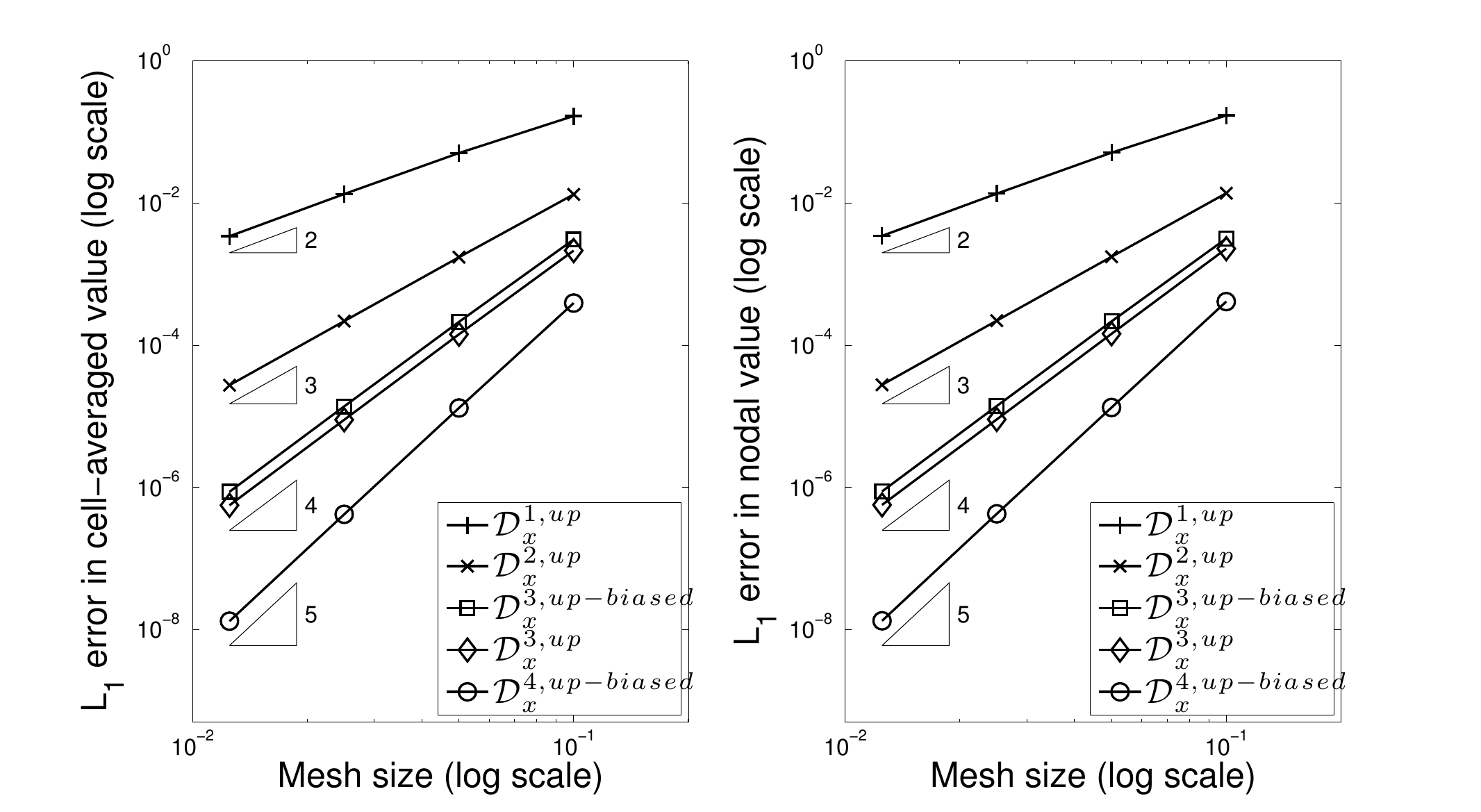}
    \caption{$L_1$-errors computed by applying the FD-FV methods in Table~\ref{tb:appr_stab_cfl} to solve~(\ref{eq:num_1dadv_per}): cell averages (left), and nodal values (right).}
    \label{fg:num_1dadv_per}
\end{figure}

First we consider an advection problem with periodic boundary conditions:
\begin{equation}\label{eq:num_1dadv_per}
\left\{\begin{array}{>{\displaystyle}r>{\displaystyle}c>{\displaystyle}l>{\displaystyle}c>{\displaystyle}l}
u_t + 2u_x & = & 0, & \quad\quad & (x,t)\in[-1,1]\times[0,1], \\
u(x,0) & = & 1 + \frac{1}{2}\sin(\pi x), & \quad\quad & x\in[-1,1], \\
u(-1,t) & = & u(1,t) & \quad\quad & t\in[0,1].
\end{array}\right.
\end{equation}
The $L_1$-norms of the errors for both cell averages and nodal values are computed using the exact solution at $T=1$, and these errors are plotted in logarithmic scales in Figure~\ref{fg:num_1dadv_per}.
All the five methods in Table~\ref{tb:appr_stab_cfl} are tested, and the Courant numbers are selected to be $90\%$ the values of $\nu_{\max}$ in the last column.
For each of the method, we compute the numerical solutions using four uniform meshes with numbers of cells ranging from $20$ to $160$.

\begin{figure}\centering
    \includegraphics[width=\textwidth]{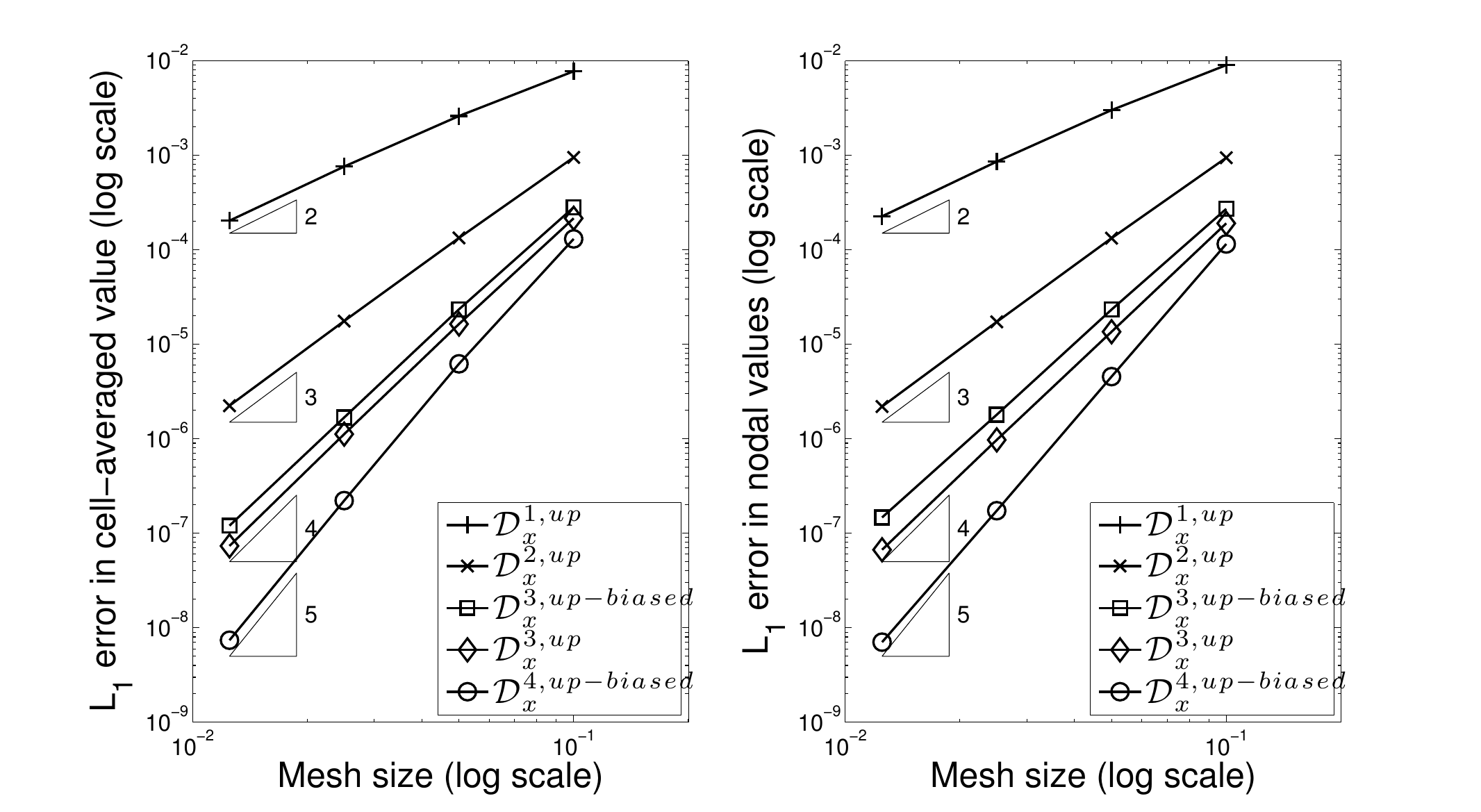}
    \caption{$L_1$-errors computed by applying the FD-FV methods in Table~\ref{tb:appr_stab_cfl} to solve~(\ref{eq:num_1dadv_dir}): cell averages (left), and nodal values (right).}
    \label{fg:num_1dadv_dir}
\end{figure}

Next, we consider an advection problem with Dirichlet boundary condition:
\begin{equation}\label{eq:num_1dadv_dir}
\left\{\begin{array}{>{\displaystyle}r>{\displaystyle}c>{\displaystyle}l>{\displaystyle}c>{\displaystyle}l}
u_t + u_x & = & 0, & \quad\quad & (x,t)\in[-0.5,0.5]\times[0,0.5], \\
u(x,0) & = & 1 + \frac{1}{2}x^3\sin(2\pi x)\bb{1}_{\{x\le0\}}, & \quad\quad & x\in[-0.5,0.5], \\
u(-0.5,t) & = & 1 + \frac{1}{2}(t+\frac{1}{2})^3\sin(2\pi(t+\frac{1}{2})), & \quad\quad & t\in[0,0.5],
\end{array}\right.
\end{equation}
where $\bb{1}_{\{x\le0\}}$ takes the value $1$ if $x\le 0$ or $0$ if $x>0$.
Using the same numerical methods and the uniform grids with the same number of cells as before, the $L_1$-norms of the numerical errors are computed and plotted in logarithmic scales in Figure~\ref{fg:num_1dadv_dir}.
Note that in this problem, no boundary condition is specified at the right end because the latter only has outgoing characteristics.
From a numerical point of view, however, applying $[\mathcal{D}_x^{3,up-biased}]$ and $[\mathcal{D}_x^{4,up-biased}]$ at the last node requires data to the right of the boundary, which is not available.
For this reason, we locally use a purely upwind operator, namely $[\mathcal{D}_x^{3,up}]$, at the last node in these two situations.
Because this change of operator is local in space, it does not affect the overall order of accuracy, especially for the fifth-order scheme constructed by $[\mathcal{D}_x^{4,up-biased}]$.

Figures~\ref{fg:num_1dadv_per}--\ref{fg:num_1dadv_dir} show that the observed orders of accuracy of the FD-FV methods are consistent with that predicted by Theorem~\ref{thm:appr_acc_fdfv}, at least for linear problems.
In particular, we achieve the desired order of accuracy for both initial value problem and initial boundary value problem.

\begin{remark}\label{rk:num_1dadv_dir}
We do not use any special treatment to enforce the Dirichlet boundary condition at intermediate Runge-Kutta stages to solve (\ref{eq:num_1dadv_dir}).
As pointed out by Mark H. Carpenter et al.~\cite{MHCarpenter:1995a}, this strategy may not be appropriate for high-order Runge-Kutta schemes. 
However, there seems to be no fix of this issue in literature for general explicit Runge-Kutta methods. 
Thus we do not apply the fix for the fourth-order Runge-Kutta in their work in our computations; 
nevertheless our numerical solutions show no degrading of the spatial orders of accuracy.
\end{remark}

\subsubsection{One-dimensional Euler equations}
\label{sec:num_1deuler}
Our next example concerns the nonlinear Euler equations in 1D.
In particular, we consider solving the problem~(\ref{eq:appr_bc_1deuler}) on $(x,t)\in[-1,1]\times[0,0.3]$ with the following smooth initial conditions and periodic boundary conditions:
\begin{equation}\label{eq:num_1deuler_per}
\left\{\begin{array}{>{\displaystyle}l>{\displaystyle}c>{\displaystyle}l}
\rho(x,0) = 1+\frac{1}{2}\sin(\pi x), \ 
u(x,0) = 2+\frac{1}{2}\sin(\pi x), \ 
p(x,0) = 1+\frac{1}{2}\sin(\pi x), &\quad& x\in[-1,1];\\ \\
\rho(-1,t) = \rho(1,t),\ 
u(-1,t) = u(1,t),\ 
p(-1,t) = p(1,t), &\quad& t\in[0,0.3].
\end{array}\right.
\end{equation}
The solutions to this problem at $T=0.3$ is sufficiently smooth for our purpose to measure the order of accuracy by applying the FD-FV methods to solve a nonlinear problem.
In particular, we compute the reference solution using $[\mathcal{D}_x^{4,up-biased}]$ and the fifth-order Runge-Kutta method on a very fine grid with $2,\!560$ uniform cells.

\begin{figure}\centering
    \includegraphics[width=\textwidth]{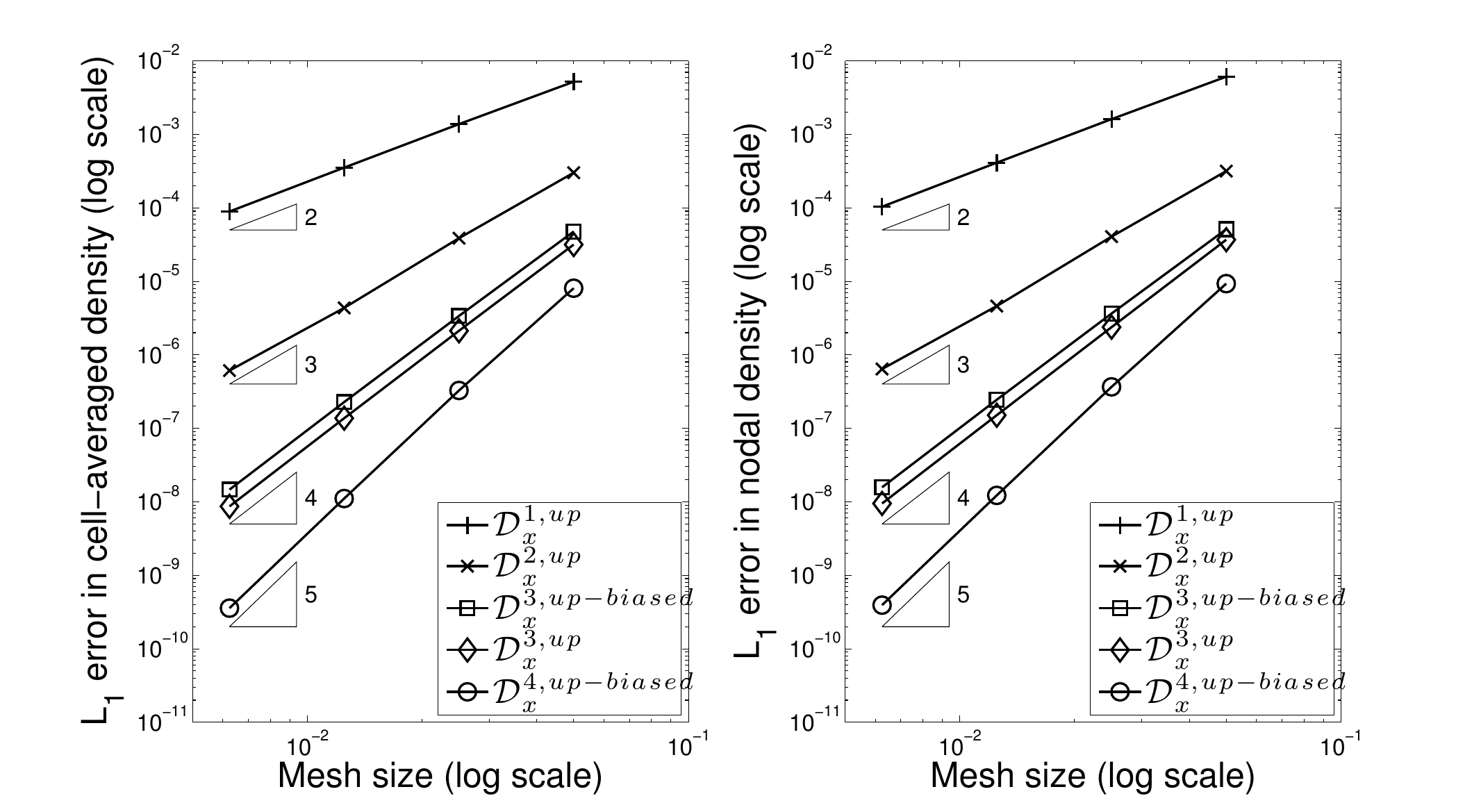}
    \caption{$L_1$-errors in density computed by applying the FD-FV methods in Table~\ref{tb:appr_stab_cfl} to solve~(\ref{eq:appr_bc_1deuler}) and (\ref{eq:num_1deuler_per}): cell-averaged density (left), and nodal density (right).}
    \label{fg:num_1deuler_per}
\end{figure}

We test the same set of numerical methods to solve this problem, using four successively refined meshes with the number of uniform cells ranging from $40$ to $320$.
The computed $L_1$-norms in density are plotted in logarithmic scales in Figure~\ref{fg:num_1deuler_per}; and the convergence curves clearly show that the desired spatial orders are achieved for the nonlinear Euler equations.

For completeness, we also present the $L_1$-norms of the errors in velocity and pressure in Table~\ref{tb:num_1deuler_per}, and they also demonstrate the target orders of convergence.
\begin{table}\centering
\begin{tabular}{|c|c|cc|cc|cc|cc|}
\hline
\multicolumn{2}{|c|}{\multirow{2}{*}{}} & \multicolumn{2}{c|}{$u$} & \multicolumn{2}{c|}{$u(\overline{\bs{w}})$}
& \multicolumn{2}{c|}{$p$} & \multicolumn{2}{c|}{$p(\overline{\bs{w}})$} \\ \hline
Method & Cells & \multicolumn{1}{c|}{$L_1$ errors} & Rates & \multicolumn{1}{c|}{$L_1$ errors} & Rates & \multicolumn{1}{c|}{$L_1$ errors} & Rates & \multicolumn{1}{c|}{$L_1$ errors} & Rates \\ \hline
\multirow{4}{*}{$[\mathcal{D}_x^{1,up}]$}& $40$  & $6.61\textrm{\sc{e}-}3$ &        & $4.72\textrm{\sc{e}-}3$ &		      & $7.54\textrm{\sc{e}-}3$ &        & $5.98\textrm{\sc{e}-}3$ &		\\
											& $80$  & $1.80\textrm{\sc{e}-}3$ & $1.87$ & $1.33\textrm{\sc{e}-}3$ & $1.83$ & $2.01\textrm{\sc{e}-}3$ & $1.91$ & $1.61\textrm{\sc{e}-}3$ & $1.90$ \\
											& $160$ & $4.70\textrm{\sc{e}-}4$ & $1.94$ & $3.50\textrm{\sc{e}-}4$ & $1.92$ & $5.17\textrm{\sc{e}-}4$ & $1.96$ & $4.14\textrm{\sc{e}-}4$ & $1.96$ \\
											& $320$ & $1.20\textrm{\sc{e}-}4$ & $1.97$ & $8.98\textrm{\sc{e}-}5$ & $1.96$ & $1.31\textrm{\sc{e}-}4$ & $1.98$ & $1.05\textrm{\sc{e}-}4$ & $1.98$ \\ \hline
\multirow{4}{*}{$[\mathcal{D}_x^{2,up}]$}& $40$  & $4.05\textrm{\sc{e}-}4$ &        & $3.54\textrm{\sc{e}-}4$ &		      & $4.62\textrm{\sc{e}-}4$ &        & $4.10\textrm{\sc{e}-}4$ &		\\
											& $80$  & $5.18\textrm{\sc{e}-}5$ & $2.97$ & $4.54\textrm{\sc{e}-}5$ & $2.96$ & $5.92\textrm{\sc{e}-}5$ & $2.96$ & $5.23\textrm{\sc{e}-}5$ & $2.97$ \\
											& $160$ & $5.88\textrm{\sc{e}-}6$ & $3.14$ & $5.11\textrm{\sc{e}-}6$ & $3.15$ & $6.68\textrm{\sc{e}-}6$ & $3.15$ & $5.83\textrm{\sc{e}-}6$ & $3.17$ \\
											& $320$ & $8.21\textrm{\sc{e}-}7$ & $2.84$ & $7.22\textrm{\sc{e}-}7$ & $2.82$ & $9.33\textrm{\sc{e}-}7$ & $2.84$ & $8.26\textrm{\sc{e}-}7$ & $2.82$ \\ \hline
\multirow{4}{*}{$[\mathcal{D}_x^{3,up-biased}]$} 
											& $40$  & $7.66\textrm{\sc{e}-}5$ &        & $6.61\textrm{\sc{e}-}5$ &		  & $7.85\textrm{\sc{e}-}5$ &        & $6.49\textrm{\sc{e}-}5$ &		\\
											& $80$  & $5.20\textrm{\sc{e}-}6$ & $3.88$ & $4.46\textrm{\sc{e}-}6$ & $3.89$ & $5.29\textrm{\sc{e}-}6$ & $3.89$ & $4.26\textrm{\sc{e}-}6$ & $3.93$ \\
											& $160$ & $3.36\textrm{\sc{e}-}7$ & $3.95$ & $2.87\textrm{\sc{e}-}7$ & $3.96$ & $3.39\textrm{\sc{e}-}7$ & $3.96$ & $2.72\textrm{\sc{e}-}7$ & $3.97$ \\
											& $320$ & $2.13\textrm{\sc{e}-}8$ & $3.98$ & $1.82\textrm{\sc{e}-}8$ & $3.98$ & $2.13\textrm{\sc{e}-}8$ & $3.99$ & $1.71\textrm{\sc{e}-}8$ & $3.99$ \\ \hline
\multirow{4}{*}{$[\mathcal{D}_x^{3,up}]$}& $40$  & $4.39\textrm{\sc{e}-}5$ &        & $3.63\textrm{\sc{e}-}5$ &		      & $4.32\textrm{\sc{e}-}5$ &        & $3.27\textrm{\sc{e}-}5$ &		\\
											& $80$  & $2.82\textrm{\sc{e}-}6$ & $3.96$ & $2.43\textrm{\sc{e}-}6$ & $3.90$ & $2.82\textrm{\sc{e}-}6$ & $3.94$ & $2.21\textrm{\sc{e}-}6$ & $3.88$ \\
											& $160$ & $1.79\textrm{\sc{e}-}7$ & $3.98$ & $1.54\textrm{\sc{e}-}7$ & $3.98$ & $1.78\textrm{\sc{e}-}7$ & $3.98$ & $1.41\textrm{\sc{e}-}7$ & $3.97$ \\
											& $320$ & $1.12\textrm{\sc{e}-}8$ & $3.99$ & $9.70\textrm{\sc{e}-}8$ & $3.99$ & $1.12\textrm{\sc{e}-}8$ & $3.99$ & $8.89\textrm{\sc{e}-}8$ & $3.99$ \\ \hline
\multirow{4}{*}{$[\mathcal{D}_x^{4,up-biased}]$} 
											& $40$  & $8.12\textrm{\sc{e}-}6$  &        & $6.63\textrm{\sc{e}-}6$  &		& $8.40\textrm{\sc{e}-}6$  &        & $6.51\textrm{\sc{e}-}6$  &		\\
											& $80$  & $2.94\textrm{\sc{e}-}7$  & $4.79$ & $2.69\textrm{\sc{e}-}7$  & $4.62$ & $3.09\textrm{\sc{e}-}7$  & $4.76$ & $2.64\textrm{\sc{e}-}7$  & $4.62$ \\
											& $160$ & $9.78\textrm{\sc{e}-}9$  & $4.91$ & $9.09\textrm{\sc{e}-}9$  & $4.89$ & $1.02\textrm{\sc{e}-}8$  & $4.92$ & $9.13\textrm{\sc{e}-}9$  & $4.86$ \\
											& $320$ & $3.15\textrm{\sc{e}-}10$ & $4.96$ & $2.93\textrm{\sc{e}-}10$ & $4.95$ & $3.27\textrm{\sc{e}-}10$ & $4.97$ & $2.94\textrm{\sc{e}-}10$ & $4.96$ \\ \hline
\end{tabular}
\caption{$L_1$ errors in velocity and pressure computed by applying the FD-FV methods in Table~\ref{tb:appr_stab_cfl} to solve (\ref{eq:appr_bc_1deuler}) and (\ref{eq:num_1deuler_per}).}
\label{tb:num_1deuler_per}
\end{table}
\begin{remark}\label{rk:num_1deuler_per}
In order to compute the numerical errors in the prime variables $u$ and $p$, especially for the cell-averaged ones, we need to compute these quantities as functions of the conservative variables.
In particular, denoting the cell-averaged values for the conservative variables by $\overline{w} = [\overline{\rho}\ \overline{\rho u}\ \overline{E}]^T$, the ``cell-averaged'' prime variables $u(\overline{\bs{w}})$ and $p(\overline{\bs{w}})$ are computed using:
\begin{displaymath}
u(\overline{\bs{w}}) = \frac{\overline{\rho u}}{\overline{\rho}},\quad
p(\overline{\bs{w}}) = (\gamma-1)\left(\overline{E}-\frac{1}{2}\frac{(\overline{\rho u})^2}{\overline{\rho}}\right).
\end{displaymath}
\end{remark}

\subsubsection{One-dimensional Buckley-Leverett equation}
\label{sec:num_1dbl}
The last 1D example is a scalar conservation law~(\ref{eq:appr_gen_1dhcl}) defined by the Buckley-Leverett flux function $f(w) = \frac{4w^2}{4w^2+(1-w)^2}$~\cite{SAKarabasov:2009a}.
We solve this problem on the domain $(x,t)\in[-0.5,2.5]\times[0,0.5]$, with the initial condition:
\begin{equation}\label{eq:num_1dbl_ic}
w(x,0) = \left\{\begin{array}{>{\displaystyle}l>{\displaystyle}c>{\displaystyle}l}
0.1 & \quad\quad & -0.5<x<0.0, \\
0.5 & \quad\quad & 0.0<x<1.0, \\
0.2 & \quad\quad & 1.0<x<2.0.
\end{array}\right.
\end{equation}
The solution to this problem has more complex behavior: There is an inflection point inside each of the upstream wave and downstream wave that is given by the solution $0<w_0<1, \frac{d^2f(w_0)}{dx^2} = 0$.
This leads to the solution structure containing both shock and rarefaction in each of the two waves.

\begin{figure}\centering
    \includegraphics[trim=1.2in 0.1in 1.2in 0.3in, clip, width=.48\textwidth]{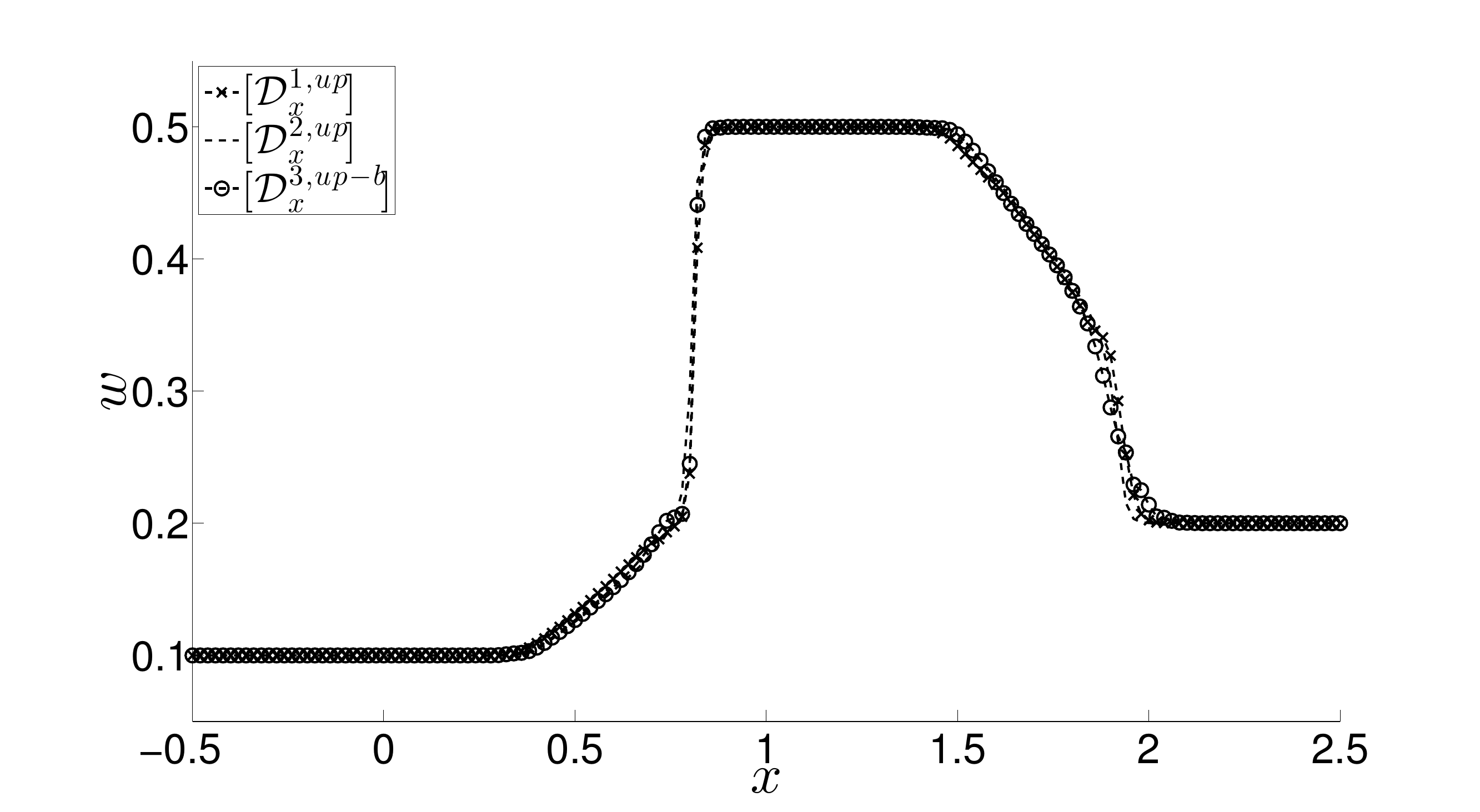}
    \includegraphics[trim=1.2in 0.1in 1.2in 0.3in, clip, width=.48\textwidth]{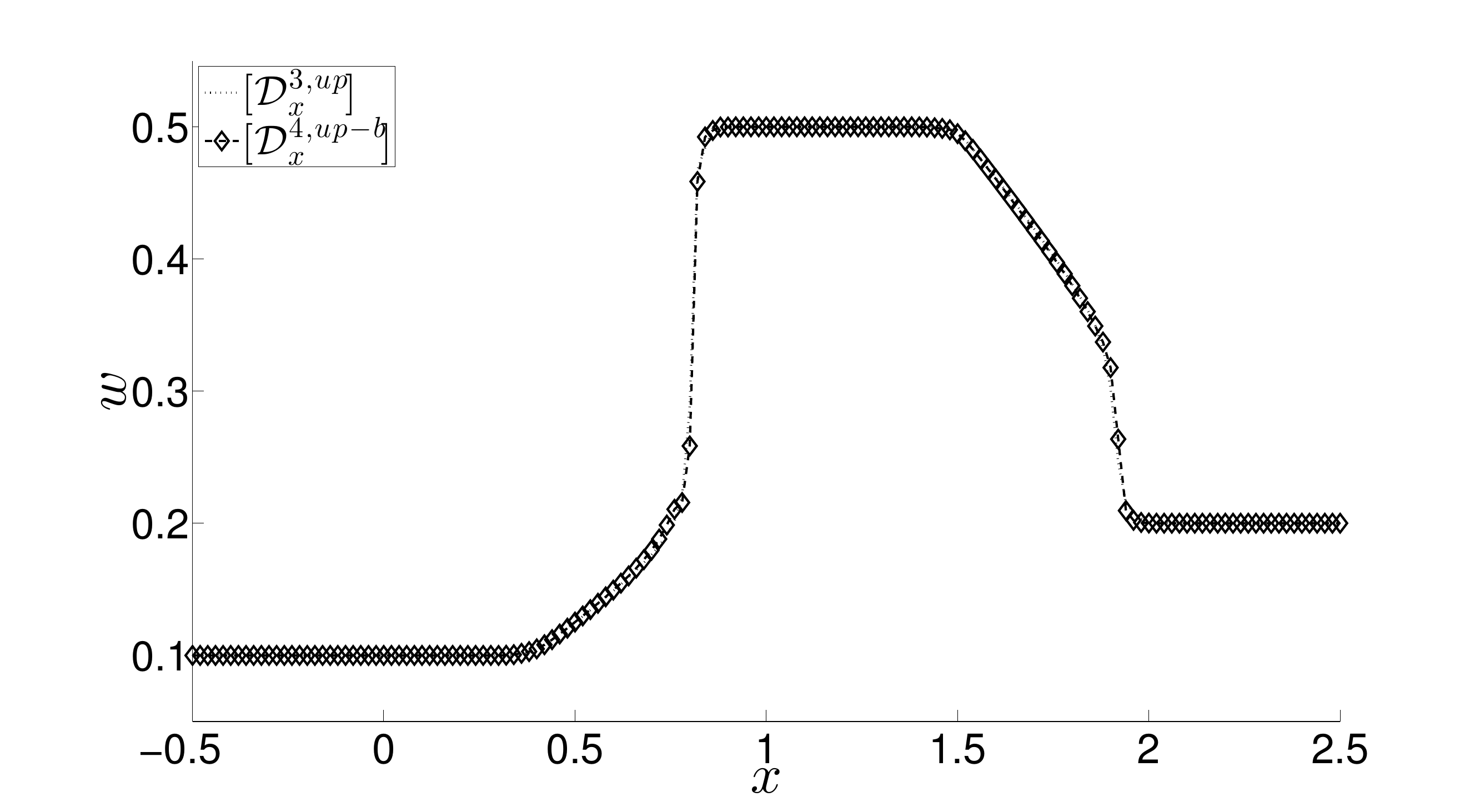}
\caption{Solutions to the Buckley-Leverett problem with initial condition~(\ref{eq:num_1dbl_ic}): (left) -- $[\mathcal{D}_x^{1,up}]$, $[\mathcal{D}_x^{2,up}]$, and $[\mathcal{D}_x^{3,up-biased}]$; and (right) -- $[\mathcal{D}_x^{3,up}]$ and $[\mathcal{D}_x^{4,up-biased}]$.}
\label{fg:num_1dbl}
\end{figure}

The numerical solutions computed using the preceding numerical methods on a uniform mesh using $150$ cells are plotted in Figures~\ref{fg:num_1dbl}.
In these figures, only the nodal values are used to plot the solutions, and the cell-averaged values are almost on top of the corresponding solution curves.

\begin{remark}\label{rk:num_1dbl}
Although the methods are constructed using the Taylor series expansions, they happen to lead to correct solution structures for this particular problem.
However, we emphasize that without further modification to enhance the nonlinear stability, the methods are in general not suitable to compute discontinuous solutions.
For example, applying the methods in Table~\ref{tb:appr_stab_cfl} to solve the classical Sod's shock tube problem~\cite{GASod:1978a} for Euler equations leads to spurious oscillations known as Gibbs phenomenon.
Nevertheless, we present the simple problem with Buckley-Leverett flux function in this paper to demonstrate its potential to capture the correct solution structure for non-convex fluxes.
A full study of the nonlinear stability and convergence to the entropy solution of the FD-FV methods, however, is beyond the scope of this paper, and its analysis and enhancement will be addressed in future work.
\end{remark}

\subsection{Two-dimensional problems}
\label{sec:num_2d}

\subsubsection{Two-dimensional Euler equations}
\label{sec:num_2deuler}
The first 2D problem is the isentropic vortex advection problem~\cite{CWShu:2003a}.
In particular, we solve the 2D Euler equations:
\begin{equation}\label{eq:num_2deuler}
\left[\begin{array}{c}
\rho \\ \rho u \\ \rho v \\ E
\end{array}\right]_t + 
\left[\begin{array}{c}
\rho u \\ \rho u^2 + p \\ \rho uv \\ (E+p)u
\end{array}\right]_x + 
\left[\begin{array}{c}
\rho v \\ \rho uv \\ \rho v^2+p \\ (E+p)v
\end{array}\right]_y
= 0
\end{equation}
on the domain $(x,y,t)\in[-5,5]^2\times[0,10]$.
Here $u$ and $v$ are velocity components, and the total energy is determined by $E=p/(\gamma-1)+\rho(u^2+v^2)/2$ with $\gamma\equiv1.4$.
We specify periodic boundary condition on all edges, and set the initial condition to be an isentropic vortex given by:
\begin{equation}\label{eq:num_2deuler_ic}
\left\{\begin{array}{>{\displaystyle}r>{\displaystyle}c>{\displaystyle}l}
\rho(x,y,0) & = & \left(1-\frac{(\gamma-1)\epsilon^2}{8\gamma \pi^2}\exp\left(1-x^2-y^2\right)\right)^{\frac{1}{\gamma-1}}, \\
u(x,y,0)    & = & 1-\frac{\epsilon y}{2\pi}\exp\left(\frac{1}{2}(1-x^2-y^2)\right), \\
v(x,y,0)    & = & 1+\frac{\epsilon x}{2\pi}\exp\left(\frac{1}{2}(1-x^2-y^2)\right), \\
p(x,y,0)    & = & \left(1-\frac{(\gamma-1)\epsilon^2}{8\gamma \pi^2}\exp\left(1-x^2-y^2\right)\right)^{\frac{\gamma}{\gamma-1}},
\end{array}\right.
\end{equation}
where the constant is chosen to be $\epsilon=5$.
The exact solution at one period $T=10$ is almost identical to the initial data, which is used to compute the numerical errors.

We test the two numerical methods in Section~\ref{sec:ext} using four uniform grids with number of cells ranging from $20^2$ to $160^2$, and summarize the $L_1$-norms of the numerical errors in Table~\ref{tb:num_2deuler_per}.
\begin{table}\centering
\begin{tabular}{|c||cc|cc||cc|cc|}
\hline
\multirow{3}{*}{Grid sizes} 
& \multicolumn{4}{c||}{Method of Section~\ref{sec:ext_2nd}} & \multicolumn{4}{c|}{Method of Section~\ref{sec:ext_3rd}} \\ \cline{2-9}
 & \multicolumn{1}{c|}{$L_1$ errors} & Rates & \multicolumn{1}{c|}{$L_1$ errors} & Rates & \multicolumn{1}{c|}{$L_1$ errors} & Rates & \multicolumn{1}{c|}{$L_1$ errors} & Rates \\ \cline{2-9}
 & \multicolumn{2}{c|}{$\rho$} & \multicolumn{2}{c||}{$\overline{\rho}$} & \multicolumn{2}{c|}{$\rho$} & \multicolumn{2}{c|}{$\overline{\rho}$} \\ \hline
 $20\times20$   & $1.94\textrm{\sc{e}+}0$ &        & $1.87\textrm{\sc{e}+}0$ &        & $2.50\textrm{\sc{e}-}1$ &        & $2.38\textrm{\sc{e}-}1$ &        \\ 
 $40\times40$   & $8.48\textrm{\sc{e}-}1$ & $1.19$ & $8.12\textrm{\sc{e}-}1$ & $1.20$ & $4.52\textrm{\sc{e}-}2$ & $2.47$ & $4.34\textrm{\sc{e}-}2$ & $2.46$ \\ 
 $80\times80$   & $2.56\textrm{\sc{e}-}1$ & $1.73$ & $2.41\textrm{\sc{e}-}1$ & $1.75$ & $6.33\textrm{\sc{e}-}3$ & $2.83$ & $6.10\textrm{\sc{e}-}3$ & $2.83$ \\ 
 $160\times160$ & $6.98\textrm{\sc{e}-}2$ & $1.88$ & $6.26\textrm{\sc{e}-}2$ & $1.95$ & $8.24\textrm{\sc{e}-}4$ & $2.94$ & $7.83\textrm{\sc{e}-}4$ & $2.96$ \\ \hline
 & \multicolumn{2}{c|}{$u$} & \multicolumn{2}{c||}{$u(\overline{\bs{w}})$} & \multicolumn{2}{c|}{$u$} & \multicolumn{2}{c|}{$u(\overline{\bs{w}})$} \\ \hline
 $20\times20$   & $6.12\textrm{\sc{e}+}0$ &        & $6.36\textrm{\sc{e}+}0$ &        & $7.29\textrm{\sc{e}-}1$ &        & $6.93\textrm{\sc{e}-}1$ &        \\ 
 $40\times40$   & $2.74\textrm{\sc{e}+}0$ & $1.16$ & $2.79\textrm{\sc{e}+}0$ & $1.19$ & $1.30\textrm{\sc{e}-}1$ & $2.49$ & $1.29\textrm{\sc{e}-}1$ & $2.43$ \\ 
 $80\times80$   & $8.25\textrm{\sc{e}-}1$ & $1.73$ & $8.19\textrm{\sc{e}-}1$ & $1.77$ & $1.81\textrm{\sc{e}-}2$ & $2.84$ & $1.87\textrm{\sc{e}-}2$ & $2.78$ \\ 
 $160\times160$ & $2.18\textrm{\sc{e}-}1$ & $1.92$ & $2.14\textrm{\sc{e}-}1$ & $1.94$ & $2.31\textrm{\sc{e}-}3$ & $2.97$ & $2.50\textrm{\sc{e}-}3$ & $2.90$ \\ \hline
 & \multicolumn{2}{c|}{$v$} & \multicolumn{2}{c||}{$v(\overline{\bs{w}})$} & \multicolumn{2}{c|}{$v$} & \multicolumn{2}{c|}{$v(\overline{\bs{w}})$} \\ \hline
 $20\times20$   & $5.46\textrm{\sc{e}+}0$ &        & $5.76\textrm{\sc{e}+}0$ &        & $6.06\textrm{\sc{e}-}1$ &        & $5.70\textrm{\sc{e}-}1$ &        \\ 
 $40\times40$   & $2.06\textrm{\sc{e}+}0$ & $1.40$ & $2.08\textrm{\sc{e}+}0$ & $1.47$ & $1.11\textrm{\sc{e}-}1$ & $2.45$ & $1.10\textrm{\sc{e}-}1$ & $2.37$ \\ 
 $80\times80$   & $5.92\textrm{\sc{e}-}1$ & $1.80$ & $5.83\textrm{\sc{e}-}1$ & $1.84$ & $1.56\textrm{\sc{e}-}2$ & $2.83$ & $1.60\textrm{\sc{e}-}2$ & $2.78$ \\ 
 $160\times160$ & $1.55\textrm{\sc{e}-}1$ & $1.93$ & $1.50\textrm{\sc{e}-}1$ & $1.97$ & $2.01\textrm{\sc{e}-}3$ & $2.96$ & $2.16\textrm{\sc{e}-}3$ & $2.89$ \\ \hline
 & \multicolumn{2}{c|}{$p$} & \multicolumn{2}{c||}{$p(\overline{\bs{w}})$} & \multicolumn{2}{c|}{$p$} & \multicolumn{2}{c|}{$p(\overline{\bs{w}})$} \\ \hline
 $20\times20$   & $2.46\textrm{\sc{e}+}0$ &        & $2.41\textrm{\sc{e}+}0$ &        & $3.65\textrm{\sc{e}-}1$ &        & $3.53\textrm{\sc{e}-}1$ &        \\ 
 $40\times40$   & $1.09\textrm{\sc{e}+}0$ & $1.17$ & $1.03\textrm{\sc{e}+}0$ & $1.23$ & $5.73\textrm{\sc{e}-}2$ & $2.67$ & $5.71\textrm{\sc{e}-}2$ & $2.63$ \\ 
 $80\times80$   & $3.26\textrm{\sc{e}-}1$ & $1.74$ & $2.91\textrm{\sc{e}-}1$ & $1.82$ & $7.88\textrm{\sc{e}-}3$ & $2.86$ & $8.08\textrm{\sc{e}-}3$ & $2.82$ \\ 
 $160\times160$ & $8.38\textrm{\sc{e}-}2$ & $1.89$ & $7.56\textrm{\sc{e}-}2$ & $1.95$ & $1.02\textrm{\sc{e}-}3$ & $2.94$ & $1.21\textrm{\sc{e}-}3$ & $2.74$ \\ \hline
\end{tabular}
\caption{$L_1$ errors of solving the isentropic vortex problem by the two methods of Section~\ref{sec:ext}.}
\label{tb:num_2deuler_per}
\end{table}
From the results we see that the two conjectures on the order of the methods are confirmed numerically, at least when the mesh is sufficiently refined.
Note that the discussion of Remark~\ref{rk:num_1deuler_per} holds here, which explains the computation of $u(\overline{w})$, $v(\overline{w})$, and $p(\overline{w})$ in Table~\ref{tb:num_2deuler_per} from the computed cell-averaged conservative variables.

\subsubsection{Two-dimensional linear elasticity problem with orthotropic materials}
\label{sec:num_2dle}
In the last example, we compute the wave propagation in an anisotropic medium that is modeled by linear elastic orthotropic material.
This is a linear problem that describes the solid dynamics of linear elastic material with infinitesimal strain assumption and the plane strain assumption.
The governing equation is given by:
\begin{equation}\label{eq:num_2dle_pde}
\left[\begin{array}{c}
\sigma_{11} \\ \sigma_{22} \\ \sigma_{12} \\ v_1 \\ v_2
\end{array}\right]_t - 
\left[\begin{array}{c}
\termc{C}_{11}v_1 \\ \termc{C}_{21}v_1 \\ \termc{C}_{33}v_2 \\ \rho_0^{-1}\sigma_{11} \\ \rho_0^{-1}\sigma_{12} 
\end{array}\right]_x -
\left[\begin{array}{c}
\termc{C}_{12}v_2 \\ \termc{C}_{22}v_2 \\ \termc{C}_{33}v_1 \\ \rho_0^{-1}\sigma_{12} \\ \rho_0^{-1}\sigma_{22} 
\end{array}\right]_y
= 0,
\end{equation}
where $\rho_0$ is a constant with the usual meaning of the density of the material, $v_1$ and $v_2$ denotes the velocity in $x\textrm{\sc-}$ and $y\textrm{\sc-}$directions, respectively; and $\sigma_{ij}$ and $\termc{C}_{ij}$ are components of the stress tensor and the stiffness matrix (in Voigt notation), respectively.

In particular, to derive (\ref{eq:num_2dle_pde}), we start with the following classical equation for linear elasticity with plane strain:
\begin{equation}\label{eq:num_2dle_orig}
\rho_0\bs{u}_{tt} = \nabla\cdot\bs{\sigma},
\end{equation}
where $\bs{u} = [u_1, u_2]^T$ is the displacement vector, and $\bs{\sigma}$ is the symmetric Cauchy stress tensor that is given by:
\begin{displaymath}
\bs{\sigma} = \left[\begin{array}{cc}
\sigma_{11} & \sigma_{12} \\ \sigma_{21} & \sigma_{22}
\end{array}\right],\quad
\sigma_{12} = \sigma_{21}.
\end{displaymath}
The stress tensor $\bs{\sigma}$ is connected to the strain tensor $\nabla\bs{u}$ by the constitutive equation, which is given in Voigt notation by:
\begin{equation}\label{eq:num_2dle_const}
\left[\begin{array}{c}
\sigma_{11} \\ \sigma_{22} \\ \sigma_{12}
\end{array}\right] = 
\left[\begin{array}{ccc}
\termc{C}_{11} & \termc{C}_{12} & 0 \\
\termc{C}_{21} & \termc{C}_{22} & 0 \\
0 & 0 & \termc{C}_{33}
\end{array}\right]
\left[\begin{array}{c}
u_{1,1} \\ u_{2,2} \\ u_{1,2}+u_{2,1}
\end{array}\right],
\end{equation}
where we use $_{,1}$ and $_{,2}$ to denote $_{,x}$ and $_{,y}$, respectively.

Then the hyperbolic equations~(\ref{eq:num_2dle_pde}) are obtained from (\ref{eq:num_2dle_orig}) and (\ref{eq:num_2dle_const}) by choosing $\sigma_{11}, \sigma_{22}, \sigma_{12}$, and $v_1\equiv du_1/dt$, $v_2\equiv du_2/dt$ as the variables.
Note that (\ref{eq:num_2dle_pde}) is a linear system, hence solving (\ref{eq:num_2dle_pde}) is equivalent to solving the conservation law obtained from the classical momentum equation.

In this test, we choose $\rho_0=3.0$ and the following stiffness matrix:
\begin{displaymath}
\left[\begin{array}{ccc}
\termc{C}_{11} & \termc{C}_{12} & 0 \\
\termc{C}_{21} & \termc{C}_{22} & 0 \\
0 & 0 & \termc{C}_{33}
\end{array}\right] =
\left[\begin{array}{ccc}
2.0 & 0.99 & 0 \\
0.99 & 0.5 & 0 \\
0 & 0 & 2.0
\end{array}\right],
\end{displaymath}
which is only slightly positively definite and creates a strong anisotropic effect.
We solve the problem on the rectangular domain $(x,y,t)\in[-1,1]^2\times[0,0.7]$ with all initial conditions set to zero.

\begin{figure}\centering
\begin{subfigure}[b]{.32\textwidth}\centering
    \includegraphics[width=\textwidth]{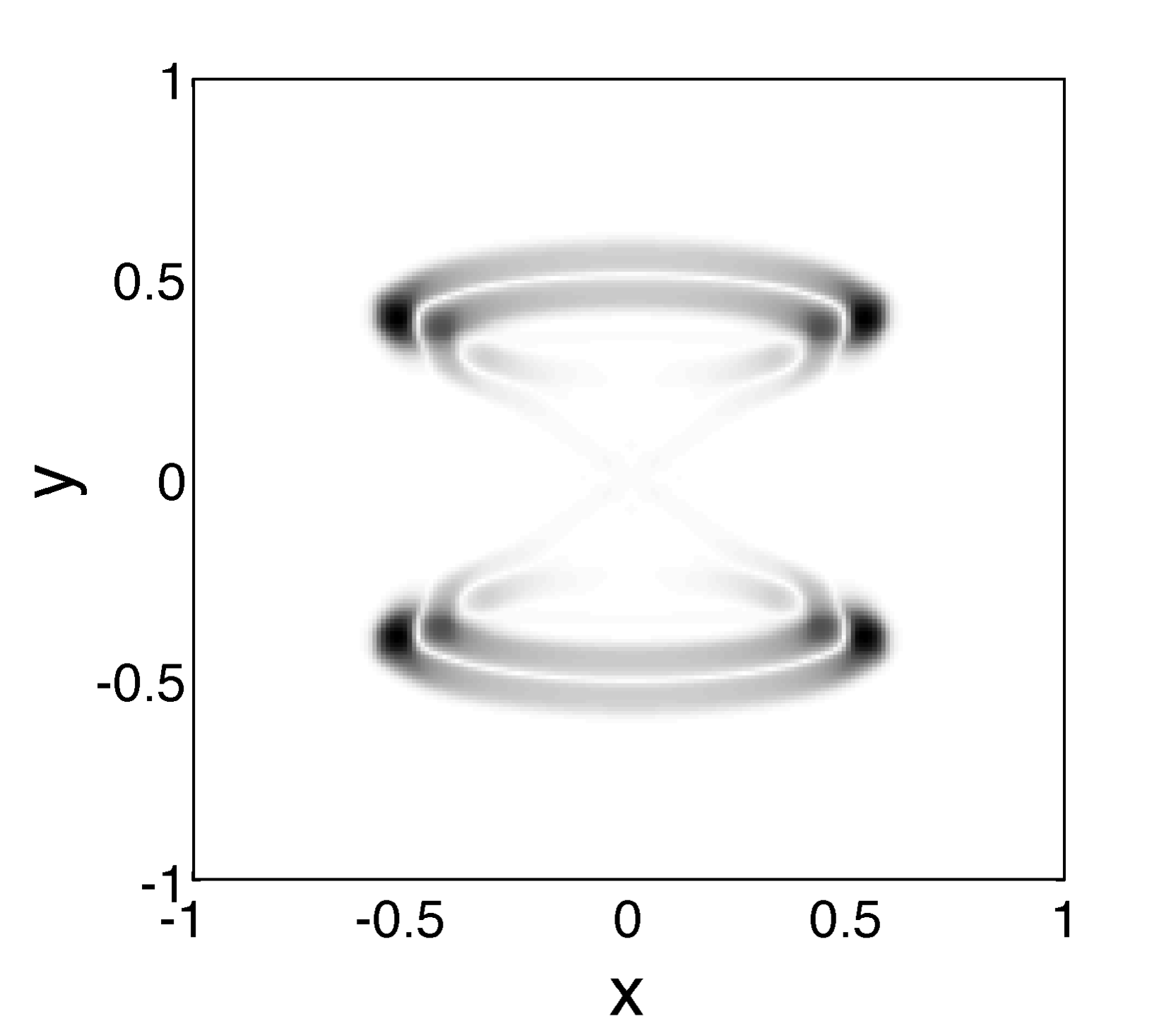}
    \caption{Finite volume method.}
    \label{fg:num_2dle_fvm}
\end{subfigure}
\begin{subfigure}[b]{.32\textwidth}\centering
    \includegraphics[width=\textwidth]{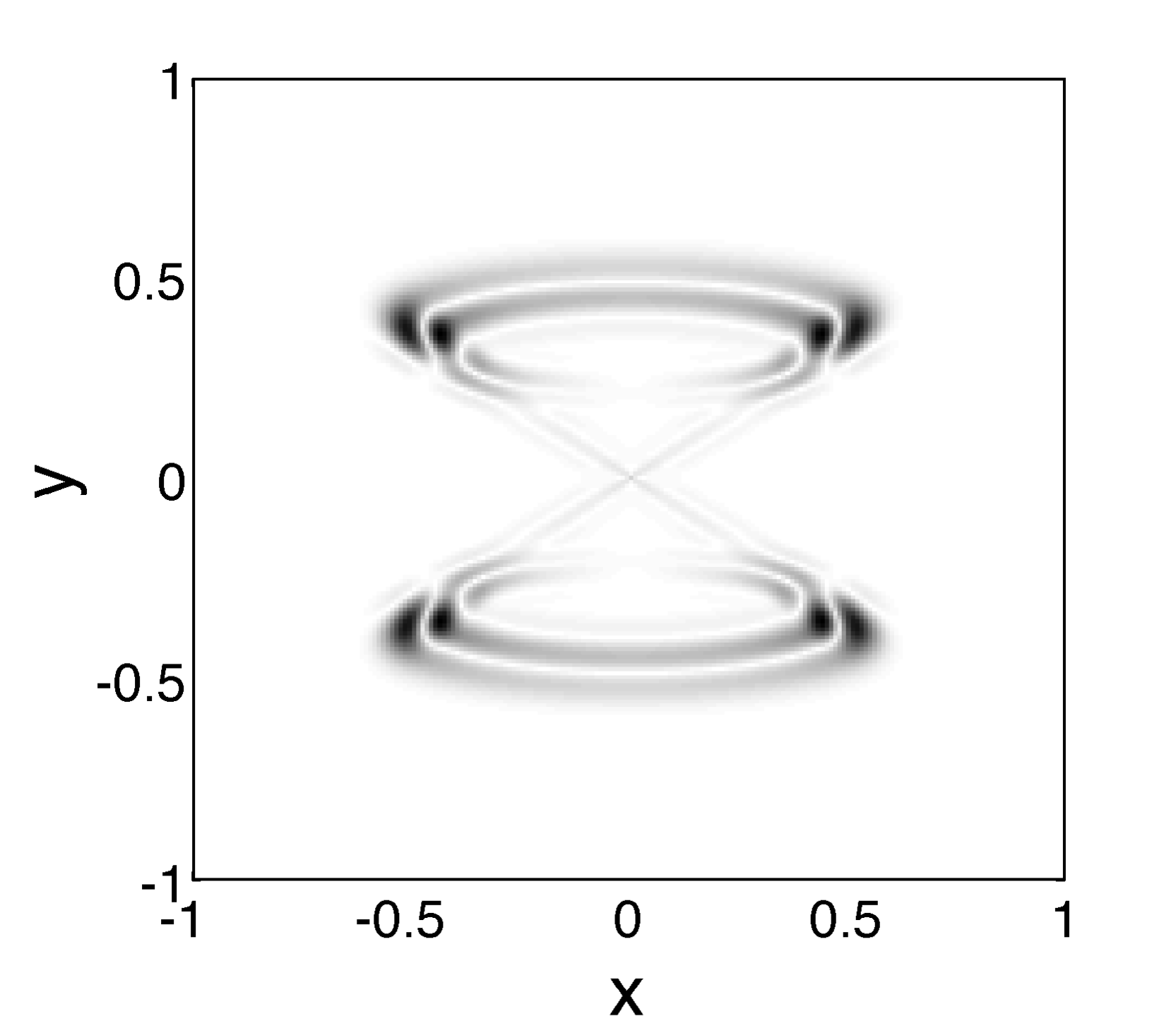}
    \caption{$2^{\pownd}\textrm{\sc-}$order FD-FV method.}
    \label{fg:num_2dle_up1}
\end{subfigure}
\begin{subfigure}[b]{.32\textwidth}\centering
    \includegraphics[width=\textwidth]{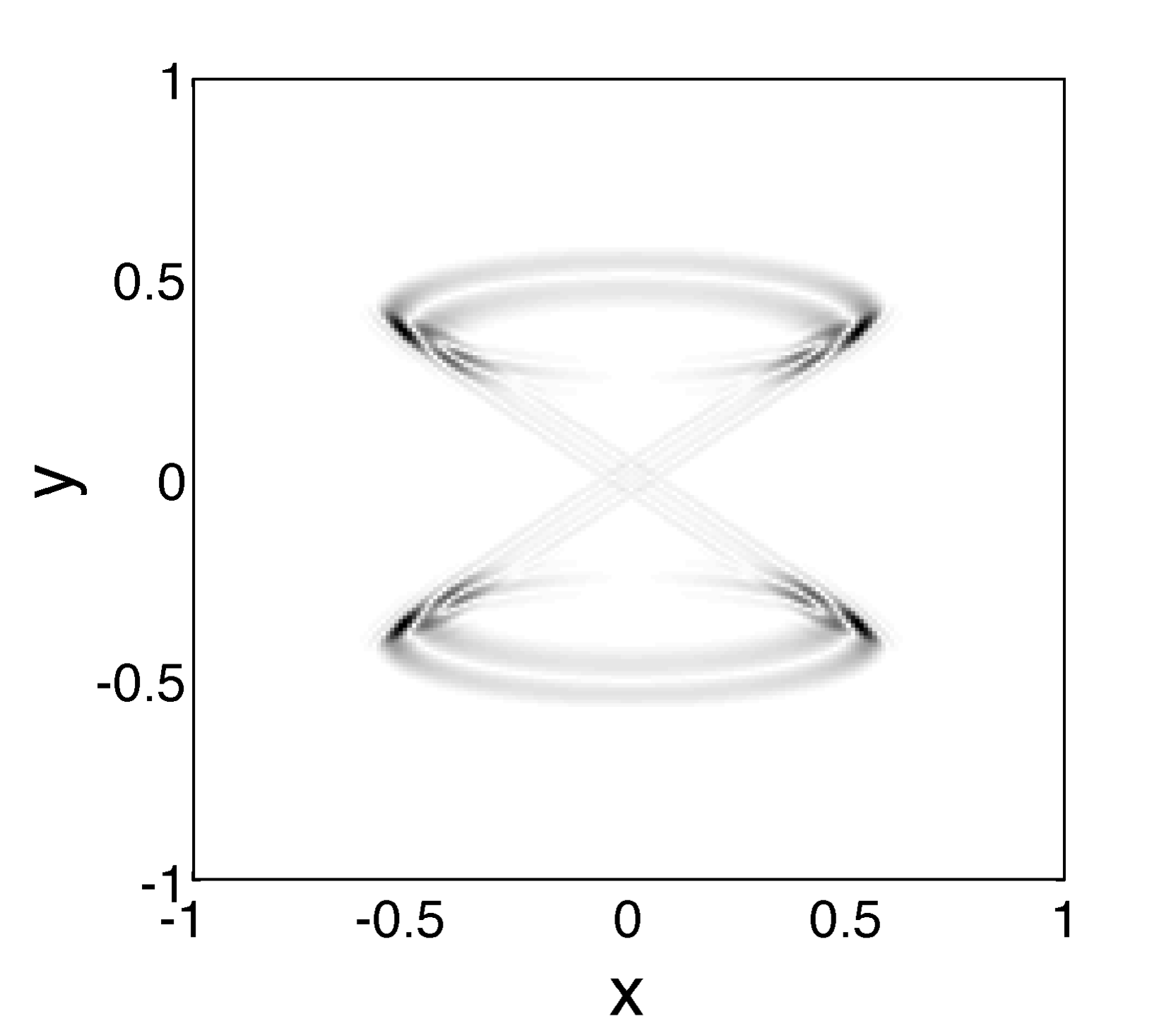}
    \caption{$3^{\powrd}\textrm{\sc-}$order FD-FV method.}
    \label{fg:num_2dle_up2}
\end{subfigure}
\caption{The velocity magnitude $\nrm{\bs{v}}=\sqrt{v_1^2+v_2^2}$ computed on a uniform $201^2$ mesh.}
\label{fg:num_2dle}
\end{figure}

To create the wave in the medium, we perturb the latter using a point pulse in the $x\textrm{\sc-}$momentum equation by adding a Ricker wavelet $s(t)$ to the right hand side of the second last equation of (\ref{eq:num_2dle_pde}):
\begin{displaymath}
s(t) = (1-2\pi^2f^2\hat{t}^2)e^{-\pi^2f^2\hat{t}^2},
\end{displaymath}
where $f=4.0$ is the peak frequency and $\hat{t} = t-(2\pi^2f^2)^{-1/2}$ so that the pulse has zero strength at $t=0.0$.
At $T=0.7$, the fastest wave has not reached the domain boundary, which allows the homogeneous Dirichlet numerical boundary condition for all variables in our computation.
This can be verified by the reference solution computed using a second-order finite volume method~\cite{BvanLeer:1979a} on a uniform $201^2$ mesh.
In particular, we use the Rusanov flux~\cite{VVRusanov:1962a}, the van Albada limiter~\cite{GDvanAlbada:1982a}, and the second-order total variational diminishing Runge-Kutta method to compute the reference solution, which is plotted in Figure~\ref{fg:num_2dle_fvm}.
The numerical solutions computed by the two FD-FV methods are plotted in Figures~\ref{fg:num_2dle_up1}--\ref{fg:num_2dle_up2}, in which the cell-averaged data are used to make the figures.
Clearly, the results confirm that the FD-FV methods compute solutions with better resolution comparing to that by the finite volume method on the same mesh.
This observation is in consistent with our prediction in that: (1) on the same mesh the FD-FV methods use more variables than the finite volume method, and (2) the FD-FV method in Section~\ref{sec:ext_3rd} is formally one-order higher than the second-order finite volume method.

\section{Conclusions}
\label{sec:concl}
We formulate a computational framework for the discretization of general first-order hyperbolic conservation law in space using both cell-averages and nodal values.
In the one-dimensional case, we prove that the formal order of spatial accuracy of the proposed FD-FV method is in general one-order higher than the order of the differential operator involved in the construction of the method.
This is unlike the conventional finite difference or finite volume methods, in which case the order of the scheme is typically the same as that of the discrete differential operator.
We extend the methodology to Cartesian grids in two space dimensions, and provide strong evidence that this property of spatial order accuracy carries to 2D smooth problems.
All the 1D and 2D methods are equipped with explicit time-integrators with appropriate order, and linear stability analysis shows that these methods are linearly stable when a Courant-Friedrich-Lewy condition is satisfied.

The numerical performance of the proposed methods is assessed by solving various smooth problems in one and two space dimensions, and the theoretical orders of accuracy are confirmed by the numerical tests.
These problems include both linear problems, such as scalar advection equation and solid dynamic equation for linear elastic orthotropic materials, and nonlinear problems, such as the 1D and 2D Euler equations.
We also show that the method is able to capture the correct solution structure for a problem with non-convex Buckley-Leverett flux function, at least in the specific problem that is considered.
A full investigation of enhancing the nonlinear stability and the convergence to the entropy solution of the proposed FD-FV methods is work in progress and will be left for future publication; 
nevertheless this Buckley-Leverett example shows the potential of the FD-FV methods for computing the correct solutions for hyperbolic problems with general flux functions.


\bibliographystyle{plain}      

\appendix

\section{The connection of the FD-FV operator and the Hermite interpolation polynomials}
\label{app:hermite}
We sketch here a proof of the existence of the constants $\alpha_l$ and $\beta_l$ that satisfy the condition~(\ref{eq:appr_acc_ddo_iff}) providing $p\le 4q$, i.e., the stencil of the operator $[\mathcal{D}_x]$ is large enough.
Furthermore, if $p<4q$ we can simply add more equations to (\ref{eq:appr_acc_ddo_iff}), so that we obtain a system of $4q+1$ equations for $4q+1$ coefficients of the form:
\begin{equation}\label{eq:hermite_sys}
\sum_{l=-q+1}^q\frac{l^{m+1}-(l-1)^{m+1}}{(m+1)!}\alpha_l + \sum_{l=-q}^q\frac{l^m}{m!}\beta_l = \gamma_mh^{-m},\quad
m=0,\cdots4q,
\end{equation}
where $\gamma_m$ are arbitrary constants.
Instead of showing the matrix on the left hand side of (\ref{eq:hermite_sys}) is nonsingular, we link it to Hermite interpolation polynomials and show that the solutions to this system exist for any right hand sides.

To this end, we consider the $2q+1$ cell faces in the stencil: $x_{j+1/2-q}<x_{j+1/2-(q-1)}<\cdots<x_{j+1/2}<x_{j+1/2+1}<\cdots<x_{j+1/2+q}$, and construct the $4q+2$ ``basis'' polynomials of degree $4q+1$ using the Hermite interpolation theory~\cite{EIsaacson:1994a}.
We denote these polynomials by $\Phi_{-q},\Phi_{-q+1},\cdots,\Phi_{q-1},\Phi_{q}$ and $\Psi_{-q},\Psi_{-q+1},\cdots,\Psi_{q-1},\Psi_q$, and they are determined by the conditions:
\begin{displaymath}
\Phi_l(x_{j+1/2+k})=\delta_{lk},\quad
\Phi_l'(x_{j+1/2+k})=0;\qquad
\Psi_l(x_{j+1/2+k})=0,\quad
\Psi_l'(x_{j+1/2+k})=\delta_{lk}
\end{displaymath}
for all $-q\le k,l\le q$, where $\delta_{lk}$ is the Kronecker delta.
Next, we define $4q+1$ polynomials of degree $4q$ by:
\begin{displaymath}
\phi_l(x) = \frac{1}{h}\sum_{k=l}^q\Phi'_k(x),\quad-q+1\le l\le q;\qquad
\psi_l(x) = \Psi'_l(x),\quad-q\le l\le q.
\end{displaymath}
Then it can be easily verified that:
\begin{equation}\label{eq:hermite_c_prop}
\frac{1}{h}\int_{\mathcal{C}_{j+k}}\phi_l(x)dx = \delta_{lk},\quad-q+1\le l,k\le q,\qquad
\phi_l(x_{j+1/2+k}) = 0,\quad -q+1\le l \le q,\ -q\le k\le q;
\end{equation}
and
\begin{equation}\label{eq:hermite_n_prop}
\frac{1}{h}\int_{\mathcal{C}_{j+k}}\psi_l(x)dx = 0,\quad-q\le l\le q,\ -q+1\le k\le q,\qquad
\psi_l(x_{j+1/2+k}) = \delta_{lk},\quad -q\le l,k \le q.
\end{equation}
Then we will verify that:
\begin{displaymath}
\alpha_l = \sum_{m=0}^{4q}\gamma_m\partial_x^m\phi_l(x),\quad -q+1\le l\le q,\quad\textrm{ and }\quad
\beta_l = \sum_{m=0}^{4q}\gamma_m\partial_x^m\psi_l(x),\quad-q\le l\le q
\end{displaymath}
satisfy the system~(\ref{eq:hermite_sys}).
Indeed, consider any polynomial $p(x)$ of degree no larger than $4q$, it may be uniquely written as:
\begin{displaymath}
p(x) = \sum_{l=-q+1}^qc_l\phi_l(x)+\sum_{l=-q}^qd_l\psi_l(x).
\end{displaymath}
Using (\ref{eq:hermite_c_prop}--\ref{eq:hermite_n_prop}), we obtain: 
\begin{displaymath}
\overline{p}_{j+l} = \frac{1}{h}\int_{\mathcal{C}_{j+l}}p(x)dx = c_l,\quad
-q+1\le l\le q;\qquad
p_{j+1/2+l} = p(x_{j+1/2+l}) = d_l,\quad
-q\le l\le q.
\end{displaymath}
Thus:
\begin{displaymath}
\begin{array}{>{\displaystyle}r>{\displaystyle}c>{\displaystyle}l}
\sum_{l=-q+1}^q\alpha_l\overline{p}_{j+l}+\sum_{l=-q}^q\beta_lp_{j+1/2+l} &=& 
\sum_{l=-q+1}^qc_l\sum_{m=0}^{4q}\gamma_m\partial_x^m\phi_l(x) +
\sum_{l=-q}^qd_l\sum_{m=0}^{4q}\gamma_m\partial_x^m\psi_l(x) \\
    &=&
\sum_{m=0}^{4q}\gamma_m\partial_x^m\left(\sum_{l=-q+1}^qc_l\phi_l(x)+\sum_{l=-q}^qd_l\psi_l(x)\right) = \sum_{m=0}^{4q}\gamma_m\partial_x^mp(x).
\end{array}
\end{displaymath}
Because this equality holds for any polynomial $p(x)$ of degree equal to or less than $4q$, by applying the Taylor series expansions to the terms on the left hand side, we see that (\ref{eq:hermite_sys}) holds.
\begin{remark}\label{rk:hermite_gen}
The preceding procedure may be modified slightly for the case of nonuniform grids, which then establish the general existence result for the FD-FV operator $[\mathcal{D}_x]$ of any order of accuracy.
\end{remark}

\section{Proof of Theorem~\ref{thm:appr_acc_fdfv}}
\label{app:order}
Rewriting the system~(\ref{eq:appr_acc_amp_ode}) in matrix form:
\begin{displaymath}
\frac{d}{dt}\left[\begin{array}{c}
A \\ N
\end{array}\right] + ik\bb{C}\left[\begin{array}{c}
A \\ N
\end{array}\right] = 0,\quad
\bb{C} = \bb{C}(\theta)\eqdef\left[\begin{array}{cc}
0 & i \\ a(\theta) & b(\theta)
\end{array}\right],
\end{displaymath}
and using the initial condition $A(0)=N(0)=1$, we may solve the ODE exactly to obtain:
\begin{equation}\label{eq:order_sol_mat}
\left[\begin{array}{c}
A(t)\\N(t)\end{array}\right] = e^{-ck\bb{C}(\theta)t}\left[\begin{array}{c}1\\1\end{array}\right].
\end{equation}
Let the two eigenvalues of $\bb{C}(\theta)$ be:
\begin{equation}\label{eq:order_eig}
\lambda_1(\theta) = \frac{1}{2}\left(b(\theta)+\sqrt{b(\theta)^2+4ia(\theta)}\right)\quad\textrm{ and }\quad
\lambda_2(\theta) = \frac{1}{2}\left(b(\theta)-\sqrt{b(\theta)^2+4ia(\theta)}\right),
\end{equation}
respectively, we compute from~(\ref{eq:order_sol_mat}):
\begin{equation}\label{eq:order_sol}
\begin{array}{>{\displaystyle}r>{\displaystyle}c>{\displaystyle}l}
A(t) &=& \frac{\lambda_2-i}{\lambda_2-\lambda_1}e^{-ck\lambda_1t} - \frac{\lambda_1-i}{\lambda_2-\lambda_1}e^{-ck\lambda_2t}, \\
N(t) &=& \frac{\lambda_1(\lambda_2-i)}{i(\lambda_2-\lambda_1)}e^{-ck\lambda_1t} - \frac{\lambda_2(\lambda_1-i)}{i(\lambda_2-\lambda_1)}e^{-ck\lambda_2t}.
\end{array}
\end{equation}

Next we compute the Laurent series of $\lambda_{1,2}$ about $\theta=0$ (as we will see soon, if $b_0=\sum_{l=-q}^q\beta_l=0$, then the eigenvalues have no pole at $0$ and we have the Taylor series expansion as in conventional finite difference analysis).
For this purpose, we define $a_m$ and $b_m$ by:
\begin{displaymath}
a_m = \frac{1}{(m+1)!}\sum_{l=-q+1}^q(l^{m+1}-(l-1)^{m+1})\alpha_l,\quad
b_m = \frac{1}{m!}\sum_{l=-q}^ql^{m}\beta_l,
\end{displaymath}
which simplifies the expression for $a(\cdot)$ and $b(\cdot)$ to:
\begin{displaymath}
a(\theta) = \frac{a_0}{\theta}+\sum_{m=0}^\infty i^{m+1}a_{m+1}\theta^m\quad\textrm{ and }\quad
b(\theta) = \frac{b_0}{\theta}+\sum_{m=0}^\infty i^{m+1}b_{m+1}\theta^m.
\end{displaymath}
By the assumption that $[\mathcal{D}_x]$ is $p^\powth\textrm{\sc-}$order accurate and Theorem~\ref{thm:appr_acc_ddo}, we have:
\begin{displaymath}
a(\theta) + b(\theta) = i + r(\theta)\theta^p,\quad\textrm{ where }\quad r(\theta)\eqdef\sum_{m=0}^\infty i^{m+p+1}(a_{m+p+1}+b_{m+p+1})\theta^m.
\end{displaymath}
Then we are ready to compute the square root in (\ref{eq:order_eig}).
Because $2i-b=-b_0/\theta+O(1)\ne0$ for small $\theta\ne0$, we have:
\begin{displaymath}
b^2 + 4ia = b^2+4i(i-b+r\theta^p) = (2i-b)^2+4ir\theta^p = \left(2i-b+\frac{2ir\theta^p}{2i-b}\right)^2+\frac{4r^2\theta^{2p}}{(2i-b)^2},
\end{displaymath}
in which there is the estimate:
\begin{displaymath}
\frac{2ir\theta^p}{2i-b} = 
\frac{2i^{p+2}c_p\theta^p+O(\theta^{p+1})}{2i-b_0/\theta+O(1)} =
\frac{2i^{p+2}c_p\theta^{p+1}+O(\theta^{p+2})}{-b_0+O(\theta)} = 
O(\theta^{p+1}).
\end{displaymath}
It follows that:
\begin{displaymath}
\sqrt{b^2+4ia} = 2i-b+\frac{2ir\theta^p}{2i-b} + O(\theta^{2p+3}) = 2i-b-\frac{2i^{p+2}c_p}{b_0}\theta^{p+1}+O(\theta^{p+2}).
\end{displaymath}
Note that we make the particular choice between the two possible roots so that using (\ref{eq:order_eig}), $\lambda_1\to i$ as $\theta\to0$.
Then the two eigenvalues are given by:
\begin{displaymath}
\lambda_1 = i-\frac{i^{p+2}c_p}{b_0}\theta^{p+1}+O(\theta^{p+2})\quad\textrm{ and }\quad
\lambda_2 = \frac{b_0}{\theta} + (b_1-1)i + O(\theta).
\end{displaymath}
Plugging these expressions into the first terms in the right hand sides of (\ref{eq:order_sol}) leads to:
\begin{equation}\label{eq:order_sol_est}
\begin{array}{>{\displaystyle}r>{\displaystyle}c>{\displaystyle}l}
\frac{\lambda_2-i}{\lambda_2-\lambda_1}e^{-ck\lambda_1t} &=& e^{-ickt}\left[1+\frac{i^{p+2}cc_pkt}{b_0}\theta^{p+1} + O(\theta^{p+2})\right], \\
\frac{\lambda_1(\lambda_2-i)}{i(\lambda_2-\lambda_1)}e^{-ck\lambda_1t} &=& e^{-ickt}\left[1+\frac{i^{p+1}c_p(ickt-1)}{b_0}\theta^{p+1} + O(\theta^{p+2})\right];
\end{array}
\end{equation}
whereas plugging them into the second terms shows that:
\begin{equation}\label{eq:order_sol_bdd}
e^{-ck\lambda_2 t} = e^{-\frac{cb_ot}{h}}e^{-ick(b_1-1)t}(1+O(\theta)),
\end{equation}
which decays to zero exponentially as $t\to\infty$ or $h\to0$ provided that $cb_0>0$ (we suppose $t,h>0$), otherwise this term grows without bound as the mesh is refined.
Hence we see that if $cb_0>0$, as stated in the theorem, the following equations hold:
\begin{displaymath}
\begin{array}{>{\displaystyle}r>{\displaystyle}c>{\displaystyle}l}
A(t) &=& e^{-ickt}\left[1+\frac{i^{p+2}cc_pkt}{b_0}\theta^{p+1} + O(\theta^{p+2})\right], \\
N(t) &=& e^{-ickt}\left[1+\frac{i^{p+1}c_p(ickt-1)}{b_0}\theta^{p+1} + O(\theta^{p+2})\right];
\end{array}
\end{displaymath}
and by Definition~\ref{def:appr_acc_fdfv}, the semi-discretized FD-FV method is at least $(p+1)^\powth\textrm{\sc-}$order accurate, which completes the proof.
\begin{remark}\label{rk:order_sym}
By defining $z_{1,2} = \theta\lambda_{1,2}(\theta)$, it is not difficult to see that both $z_1$ and $z_2$ are symmetric functions of $\theta$ about the value $\theta=\pi$.
Thus with our preceding choice of $\lambda_{1,2}$, $\lambda_1$ has significant effect as the mesh is refined, i.e., $\theta\to0$; whereas the other eigenvalue $\lambda_2$ has more significant effect when the mesh is approaching the under-resolved limit, i.e., $\theta\to2\pi$ or one cell per wave length.
\end{remark}

\section{Evidences for Conjecture~\ref{con:ext_acc_2nd}}
\label{app:2nd}
Similar to the analysis in one space dimension, we use the superscript $^\ast$ to denote the exact solutions. 
Suppose $x_{\phf{j_1}}=(\phf{j_1})h_1$ and $y_{\phf{j_2}}=(\phf{j_2})h_2$, we have the exact solutions at $t=T$:
\begin{displaymath}
\begin{array}{>{\displaystyle}l}
\overline{w}_{j_1,j_2}^\ast = \frac{1}{i^2\theta_1\theta_2}e^{-i(c_1k_1+c_2k_2)T}e^{i(j_1\theta_1+j_2\theta_2)}\left(e^{i\theta_1/2}-e^{-i\theta_1/2}\right)\left(e^{i\theta_2/2}-e^{-i\theta_2/2}\right),\\ \\
w_{\phf{j_1},j_2}^\ast = e^{-i(c_1k_1+c_2k_2)T}e^{i((\phf{j_1})\theta_1+j_2\theta_2)},\quad
w_{j_1,\phf{j_2}}^\ast = e^{-i(c_1k_1+c_2k_2)T}e^{i(j_1\theta_1+(\phf{j_2})\theta_2)},
\end{array}
\end{displaymath}
where $\theta_1 = k_1h_1$ and $\theta_2 = k_2h_2$.

It is not difficult to see that providing exact initial data at $t=0$, the solutions to the semi-discretized problem are
\begin{displaymath}
\begin{array}{>{\displaystyle}l}
\overline{w}_{j_1,j_2}(t) = \frac{1}{i^2\theta_1\theta_2}A(t)e^{i(j_1\theta_1+j_2\theta_2)}\left(e^{i\theta_1/2}-e^{-i\theta_1/2}\right)\left(e^{i\theta_2/2}-e^{-i\theta_2/2}\right),\\ \\
w_{\phf{j_1},j_2}(t) = N_1(t)e^{i((\phf{j_1})\theta_1+j_2\theta_2)},\quad
w_{j_1,\phf{j_2}}(t) = N_2(t)e^{i(j_1\theta_1+(\phf{j_2})\theta_2)},
\end{array}
\end{displaymath}
where $A(t)$, $N_1(t)$, and $N_2(t)$ solve the ODE system:
\begin{equation}\label{eq:2nd_ode}
\frac{d}{dt}\left[\begin{array}{>{\displaystyle}c}
A \\ N_1 \\ N_2
\end{array}\right]
= 
\left[
\begin{array}{>{\displaystyle}c>{\displaystyle}c>{\displaystyle}c}
0 & -ic_1k_1\tau_2^{-1} & -ic_2k_2\tau_1^{-1} \\
2c_1k_1\tau_1\tau_2\theta_1^{-1}e^{-i\theta_1/2} & -2c_1k_1\theta_1^{-1}-ic_2k_2\tau_2e^{-i\theta_2/2} & 0 \\
2c_2k_2\tau_1\tau_2\theta_2^{-1}e^{-i\theta_2/2} & 0 & -2c_2k_2\theta_2^{-1}-ic_1k_1\tau_1e^{-i\theta_1/2} 
\end{array}
\right]
\left[\begin{array}{>{\displaystyle}c}
A \\ N_1 \\ N_2
\end{array}\right],
\end{equation}
with $\tau_j = (e^{i\theta_j/2}-e^{-i\theta_j/2})/(i\theta_j) = 1 + O(\theta_j^2),\ j=1,2$.
One way to compute the three eigenvalues associated with (\ref{eq:2nd_ode}) is directly finding the roots of the cubic characteristic polynomial, just as what we have done in Appendix~\ref{app:order}.
However, this method does not extend to more general cases when the ODE system is larger; 
thus we do not perform the calculation here but simply make the following conjecture and will verify it by numerical examples afterwards:
\begin{equation}\label{eq:2nd_eigs}
\lambda_0 = -i(c_1k_1+c_2k_2)+O(\theta^2),\quad\textrm{and}\quad
\lambda_j = -\frac{2c_jk_j}{\theta_j} + O(1),\ j=1,2.
\end{equation}
In particular, we consider three sets of parameters:
\begin{quote}
\begin{enumerate}[(a)]
\item $c_1k_1=1.0$, $c_2k_2=1.0$, $\theta_1/\theta_2=5.0$;
\item $c_1k_1=5.0$, $c_2k_2=1.0$, $\theta_1/\theta_2=1.0$;
\item $c_1k_1=1.0$, $c_2k_2=1.2$, $\theta_1/\theta_2=1.0$.
\end{enumerate}
\end{quote}
For each set of the parameters, we plot the following curves to demonstrate~(\ref{eq:2nd_eigs}): 
(1) $(\Im\lambda_0+c_1k_1+c_2k_2)/\theta_2^2$ vs. $\theta_2$, 
(2) $(\Re\lambda_0)/\theta_2^2$ vs. $\theta_2$,
(3) $(\Re\lambda_j)\theta_j/(c_jk_j)$ vs. $\theta_2$, $j=1,2$,
and (4) $(\Im\lambda_j)/(c_jk_j)$ vs. $\theta_2$, $j=1,2$.
For the conjecture to be valid, all curves should converges to a finite value as $\theta_2\to0$; 
and in the case of (3), this constant value should be $-2$ for $j=1,2$.
\begin{figure}\centering
\begin{subfigure}[b]{0.45\textwidth}\centering
    \includegraphics[trim=3.6in 0.0in 3.6in 0.0in, clip, width=\textwidth]{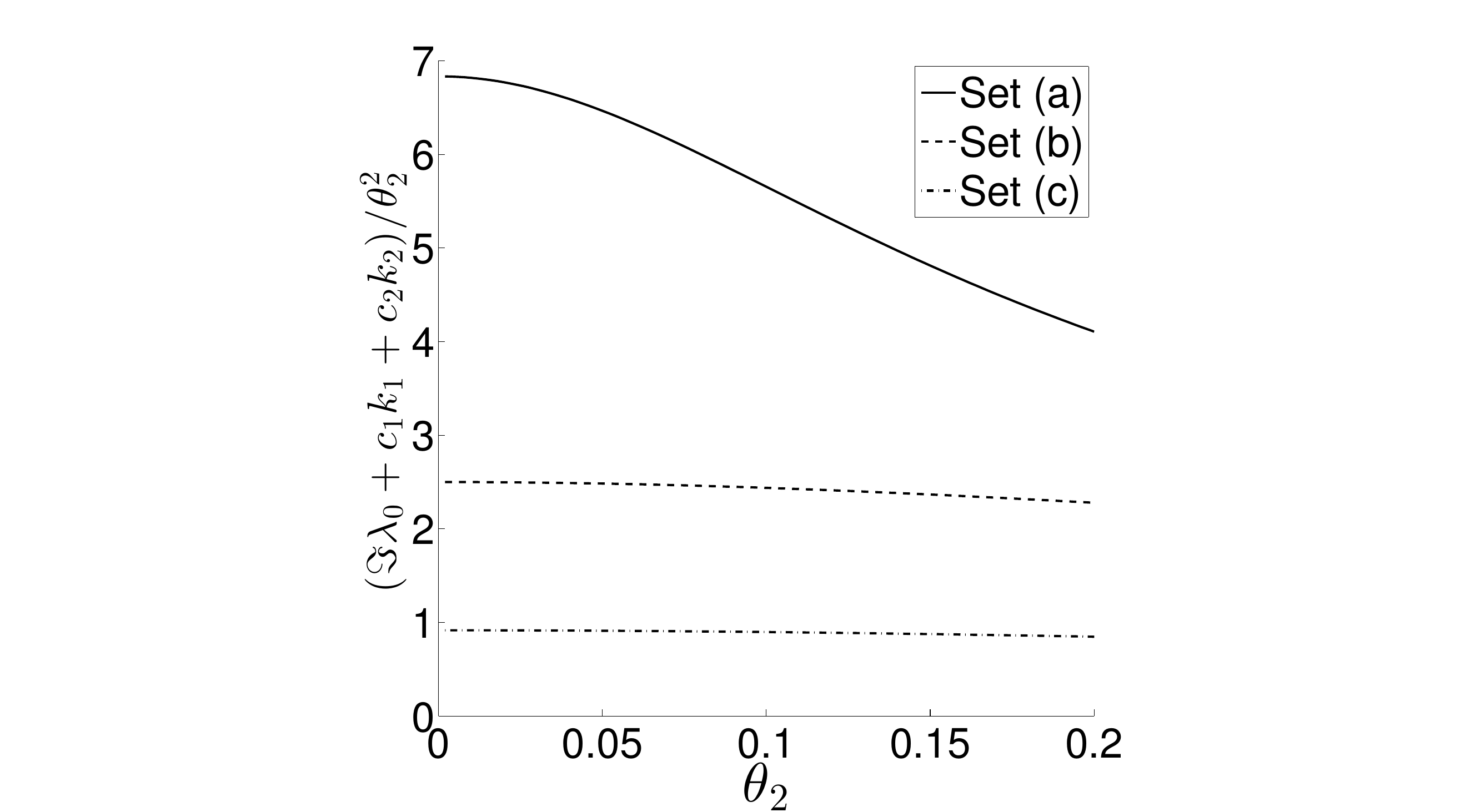}
    \caption{$\frac{\Im\lambda_0+c_1k_1+c_2k_2}{\theta_2^2}$ vs. $\theta_2$.}
    \label{fg:app_2nd_imag_l0}
\end{subfigure}
\begin{subfigure}[b]{0.45\textwidth}\centering
    \includegraphics[trim=3.6in 0.0in 3.6in 0.0in, clip, width=\textwidth]{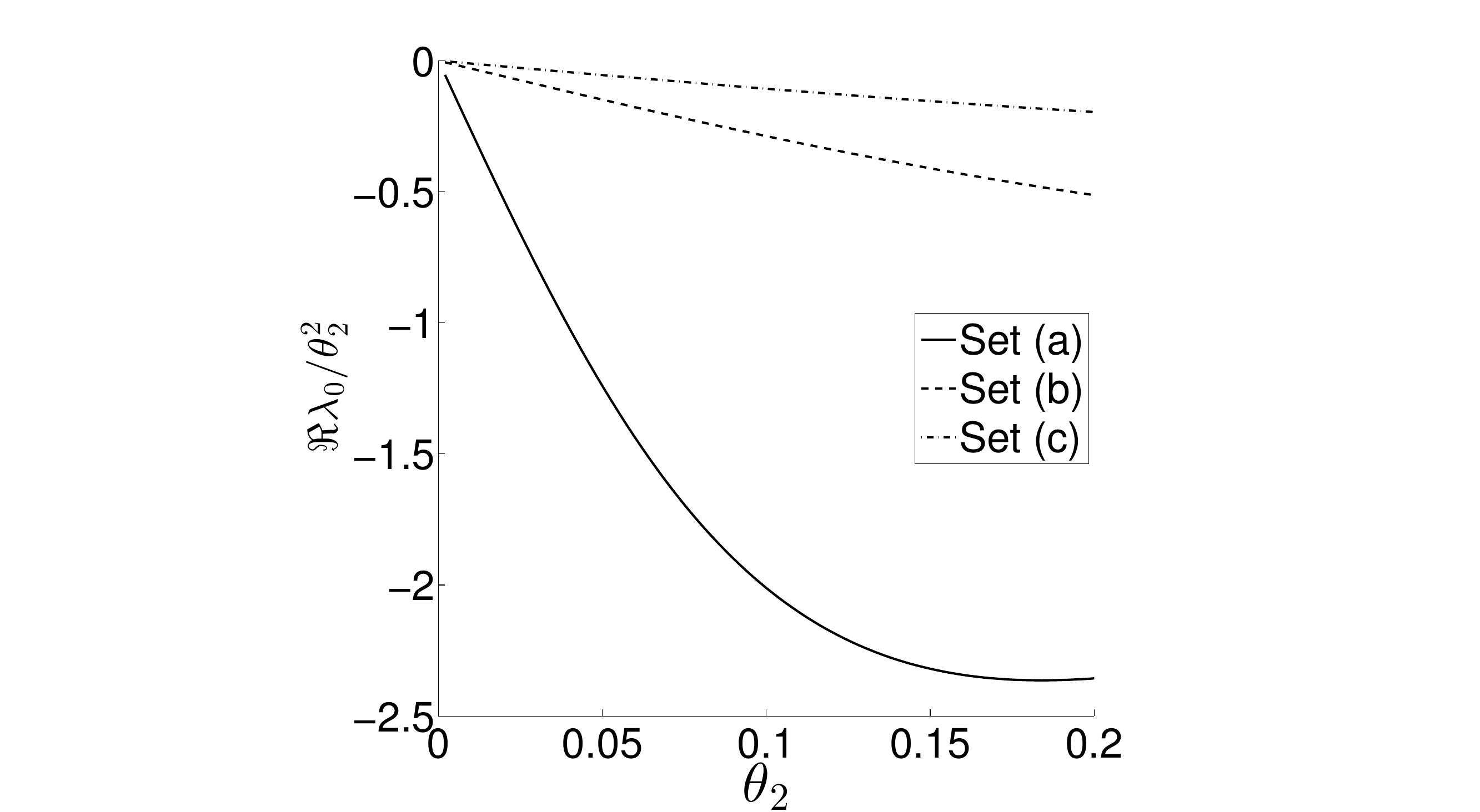}
    \caption{$\frac{\Re\lambda_0}{\theta_2^2}$ vs. $\theta_2$.}
    \label{fg:app_2nd_real_l0}
\end{subfigure}
\caption{The (\ref{fg:app_2nd_imag_l0}) imaginary and (\ref{fg:app_2nd_real_l0}) real parts of the first eigenvalue $\lambda_0$ as $\theta_2\to0$.}
\label{fg:app_2nd_l0}
\end{figure}
\begin{figure}\centering
\begin{subfigure}[b]{0.45\textwidth}\centering
    \includegraphics[trim=3.6in 0.0in 3.6in 0.0in, clip, width=\textwidth]{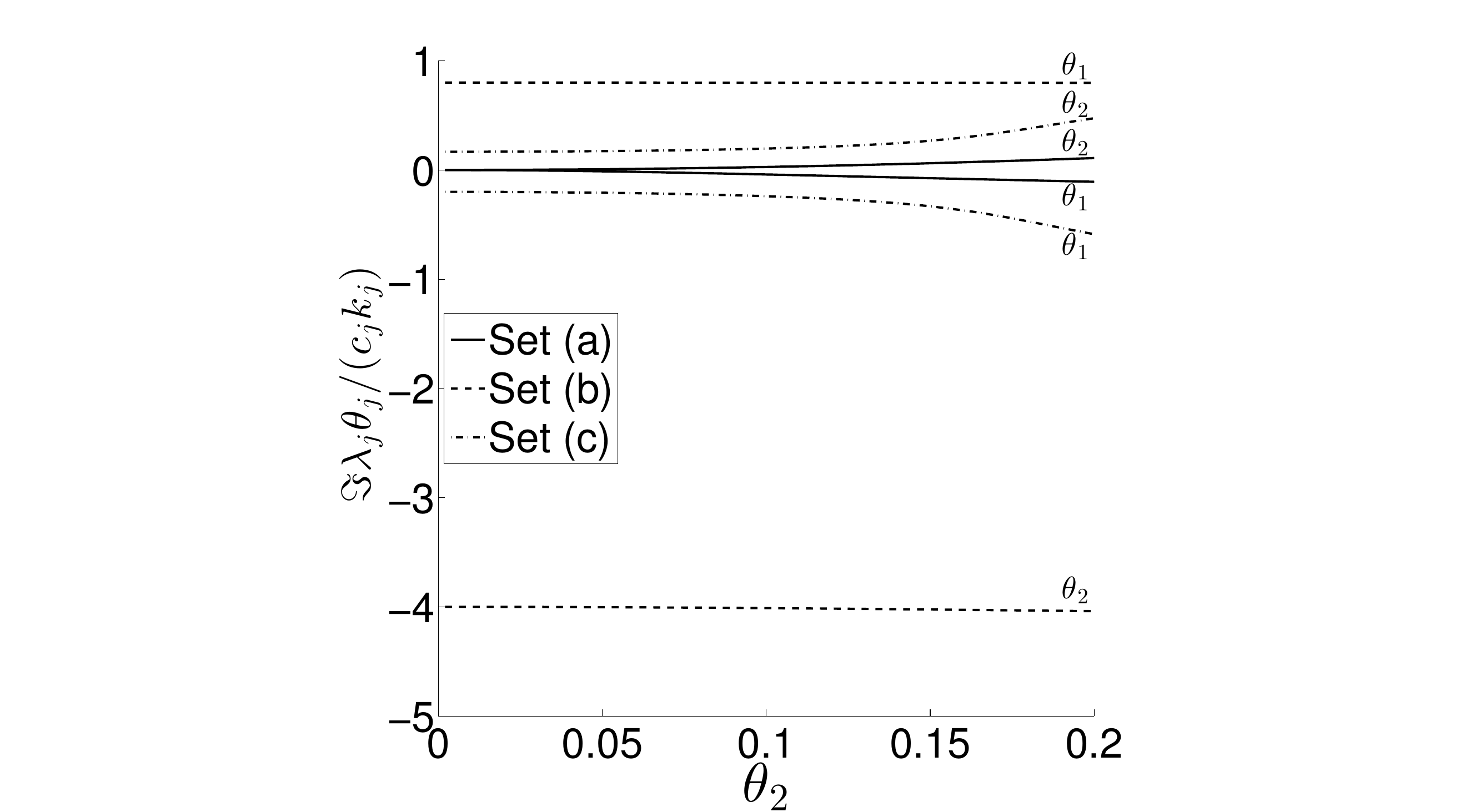}
    \caption{$\frac{\Im\lambda_j}{c_jk_j}$ vs. $\theta_2$, $j=1,2$.}
    \label{fg:app_2nd_imag_lj}
\end{subfigure}
\begin{subfigure}[b]{0.45\textwidth}\centering
    \includegraphics[trim=3.6in 0.0in 3.6in 0.0in, clip, width=\textwidth]{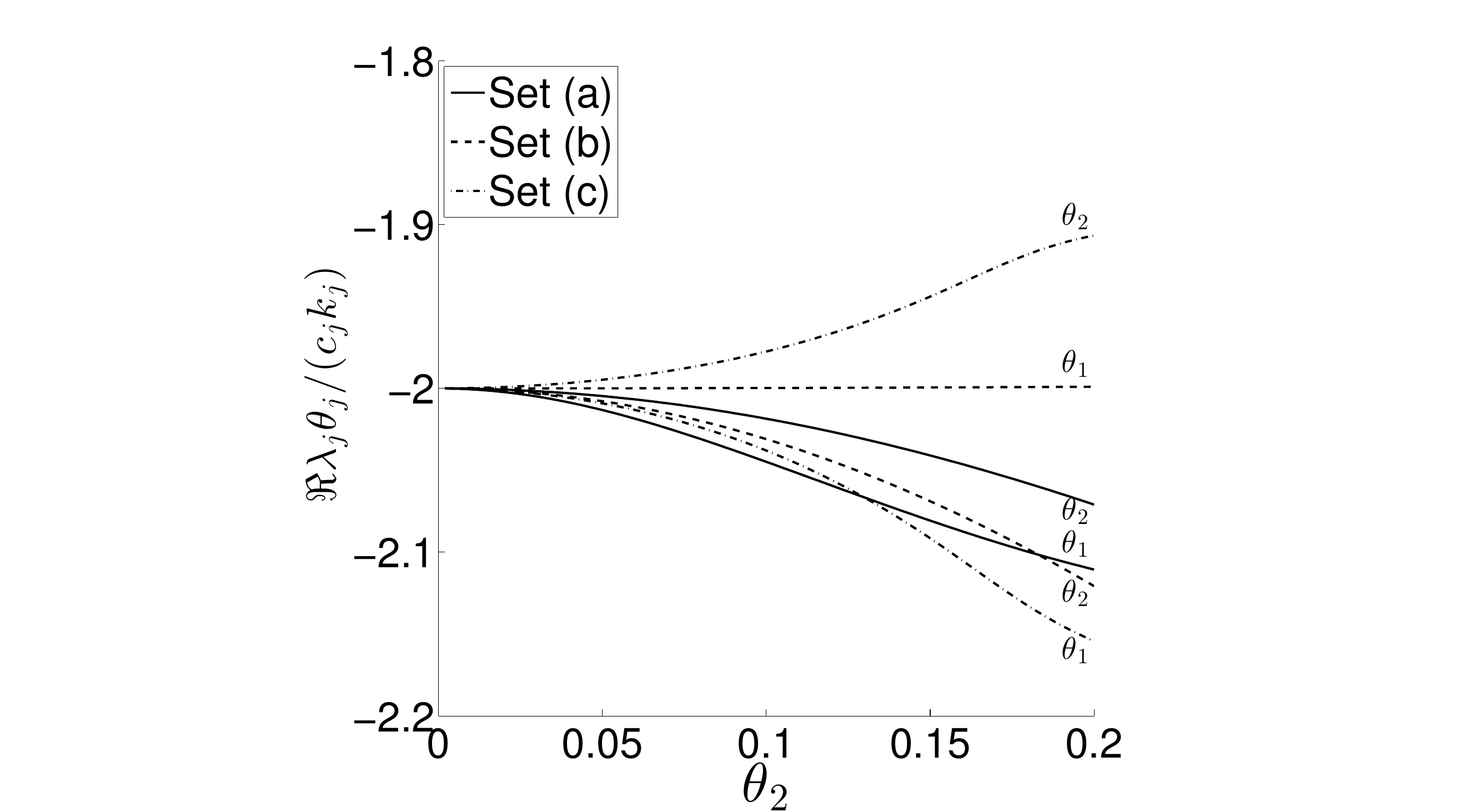}
    \caption{$\frac{\Re\lambda_j\theta_j}{c_jk_j}$ vs. $\theta_2$, $j=1,2$.}
    \label{fg:app_2nd_real_lj}
\end{subfigure}
\caption{The (\ref{fg:app_2nd_imag_lj}) imaginary and (\ref{fg:app_2nd_real_lj}) real parts of the other eigenvalue $\lambda_{1,2}$ as $\theta_2\to0$.}
\label{fg:app_2nd_lj}
\end{figure}
The asymptotic behaviors of the imaginary and real parts of $\lambda_0$ are plotted in Figure~\ref{fg:app_2nd_l0}; and those of the other two eigenvalues $\lambda_{1,2}$ are plotted in Figure~\ref{fg:app_2nd_lj}.
Note that it appears from Figure~\ref{fg:app_2nd_real_l0} that $\Re\lambda_0$ could be $o(\theta^2)$ instead of $O(\theta^2)$; 
nevertheless we stick to the latter because it is sufficient for our purpose.

Supposing that (\ref{eq:2nd_eigs}) is true, it remains to show that the solutions of (\ref{eq:2nd_ode}) are second-order in $\theta_{1,2}$.
This can be done, by explicitly constructing the formula for $A(t)$ and $N_{1,2}(t)$ similar to (\ref{eq:order_sol}).
We will only do this for $A(t)$, and the procedure for the other quantities follows similarly.
Denoting the matrix in (\ref{eq:2nd_ode}) by $\bb{C}$, and we compute its eigenvalue decomposition as:
\begin{displaymath}
\bb{C} = \bb{R}\bb{D}\bb{L},
\end{displaymath}
where $\bb{D}$ is the diagonal matrix with diagonal entries $\lambda_0, \lambda_1, \lambda_2$; and the other two matrices are given by:
\begin{displaymath}
\bb{R} = \left[\begin{array}{ccc}
1 & 1 & 1 \\
\frac{2\bsb{b}_1\tau_2}{\theta_1\lambda_0+2\bsb{a}_1+i\theta_1\bsb{b}_2} &
\frac{2\bsb{b}_1\tau_2}{\theta_1\lambda_1+2\bsb{a}_1+i\theta_1\bsb{b}_2} &
\frac{2\bsb{b}_1\tau_2}{\theta_1\lambda_2+2\bsb{a}_1+i\theta_1\bsb{b}_2} \\
\frac{2\bsb{b}_2\tau_1}{\theta_2\lambda_0+2\bsb{a}_2+i\theta_2\bsb{b}_1} &
\frac{2\bsb{b}_2\tau_1}{\theta_2\lambda_1+2\bsb{a}_2+i\theta_2\bsb{b}_1} &
\frac{2\bsb{b}_2\tau_1}{\theta_2\lambda_2+2\bsb{a}_2+i\theta_2\bsb{b}_1} 
\end{array}\right],
\end{displaymath}
and
\begin{displaymath}
\bb{L} = \left[\begin{array}{ccc}
\frac{1}{\bsb{d}_0} & \frac{-i\bsb{a}_1/\tau_2}{\bsb{d}_0(\lambda_0+2\bsb{a}_1/\theta_1+i\bsb{b}_2)} & \frac{-i\bsb{a}_2/\tau_1}{\bsb{d}_0(\lambda_0+2\bsb{a}_2/\theta_2+i\bsb{b}_1)} \\
\frac{1}{\bsb{d}_1} & \frac{-i\bsb{a}_1/\tau_2}{\bsb{d}_1(\lambda_1+2\bsb{a}_1/\theta_1+i\bsb{b}_2)} & \frac{-i\bsb{a}_2/\tau_1}{\bsb{d}_1(\lambda_1+2\bsb{a}_2/\theta_2+i\bsb{b}_1)} \\
\frac{1}{\bsb{d}_2} & \frac{-i\bsb{a}_1/\tau_2}{\bsb{d}_2(\lambda_2+2\bsb{a}_1/\theta_1+i\bsb{b}_2)} & \frac{-i\bsb{a}_2/\tau_1}{\bsb{d}_2(\lambda_2+2\bsb{a}_2/\theta_2+i\bsb{b}_1)} 
\end{array}\right],
\end{displaymath}
where $\bsb{a}_j = c_jk_j$, $\bsb{b}_j = \bsb{a}_j\tau_je^{-i\theta_j/2} = \bsb{a}_j(1+O(\theta_j))$ for $j=1,2$, and:
\begin{displaymath}
\bsb{d}_j = 1 - \frac{2i\bsb{a}_1\bsb{b}_1\theta_1}{(\lambda_j\theta_1+2\bsb{a}_1+i\bsb{b}_2\theta_1)^2} - \frac{2i\bsb{a}_2\bsb{b}_2\theta_2}{(\lambda_j\theta_2+2\bsb{a}_2+i\bsb{b}_1\theta_2)^2},\quad j = 0,1,2.
\end{displaymath}
Using~(\ref{eq:2nd_eigs}), we have:
\begin{displaymath}
\bsb{d}_0 = 1 
- \frac{2i\bsb{a}_1^2\theta_1(1+O(\theta_1))}{(2\bsb{a}_1+O(\theta))^2}
- \frac{2i\bsb{a}_2^2\theta_2(1+O(\theta_2))}{(2\bsb{a}_2+O(\theta))^2} = 
1 - \frac{1}{2}i\theta_1 - \frac{1}{2}i\theta_2 + O(\theta^2).
\end{displaymath}
Thus if we write $A(t) = A_0e^{\lambda_0t} + A_1e^{\lambda_1t} + A_2e^{\lambda_2t}$, there is:
\begin{displaymath}
\begin{array}{>{\displaystyle}r>{\displaystyle}c>{\displaystyle}l}
A_0 & = & \frac{1}{\bsb{d}_0}\left(1
        -\frac{i\bsb{a}_1\theta_1/\tau_2}{\lambda_0\theta_1+2\bsb{a}_1+i\bsb{b}_2\theta_1}
        -\frac{i\bsb{a}_2\theta_2/\tau_1}{\lambda_0\theta_2+2\bsb{a}_2+i\bsb{b}_1\theta_2}\right)
= \frac{1}{\bsb{d}_0}\left(1
        -\frac{i\bsb{a}_1\theta_1(1+O(\theta_2^2))}{2\bsb{a}_1+O(\theta)}
        -\frac{i\bsb{a}_2\theta_2(1+O(\theta_1^2))}{2\bsb{a}_2+O(\theta)}\right) \\
 & = & \frac{1}{\bsb{d}_0}\left(1-\frac{1}{2}i\theta_1-\frac{1}{2}i\theta_2+O(\theta^2)\right)
   = 1 + O(\theta^2),
\end{array}
\end{displaymath}
which confirms that $A_0e^{\lambda_0t} = e^{-i(c_1k_1+c_2k_2)t} + O(\theta^2)$.
For the remaining two terms $A_je^{\lambda_jt}$, $j=1,2$, we note that $A_j$ grows (at most) at polynomial rate of $1/\theta$; 
whereas by (\ref{eq:2nd_eigs}), $e^{\lambda_jt}$ decays exponentially in $1/\theta_j$ providing that $c_j>0$, which is assumed.
Hence we conclude that $A(t) = e^{-i(c_1k_1+c_2k_2)t}+O(\theta^2)$, and similar result holds for $N_1(t)$ and $N_2(t)$.

\end{document}